\documentclass[english]{article}
\usepackage[fontsize=11pt]{scrextend}
\usepackage[utf8]{inputenc}
\usepackage{graphicx}
\usepackage{amsmath}
\usepackage{amsfonts}
\usepackage{enumitem}   
\usepackage{amssymb}
\usepackage{comment}
\usepackage[english]{babel}
\usepackage{hyperref} 
\usepackage[top=2cm, bottom=3cm, left=2.5cm , right=2.5cm]{geometry}
\usepackage{mathrsfs}
\usepackage{dsfont}
\usepackage{color}
\usepackage{flushend}
\usepackage{mathtools}
\usepackage[ntheorem]{empheq}
\usepackage{fancyhdr}
\usepackage[thmmarks,amsmath]{ntheorem}
\usepackage{stmaryrd}
\usepackage[nottoc]{tocbibind}
\usepackage[toc,page]{appendix}
\usepackage{appendix}
\usepackage{caption}
\usepackage{xcolor}
\usepackage{float}
\usepackage{textcomp}
\usepackage{bm}  
\usepackage{lipsum}
\usepackage{titlesec}
\usepackage{braket} 
\usepackage{esint} 
\usepackage{centernot}

\titleformat{\paragraph}[block]
{\normalsize\bfseries}{}{0pt}{}

\titlespacing{\paragraph}{0pt}{*1.5}{*0.5}

\renewcommand{\textbf}[1]{\begingroup\bfseries\mathversion{bold}#1\endgroup}

\newcommand{\N}{\mathbb{N}}

\newcommand{\Q}{\mathbb{Q}}
\newcommand{\R}{\mathbb{R}}
\newcommand{\C}{\mathbb{C}}

\renewcommand{\P}{\mathbb{P}}

\newcommand{\mcal}[1]{\mathcal{#1}}
\newcommand{\mscr}[1]{\mathscr{#1}}
\newcommand{\mbb}[1]{\mathbb{#1}}

\newcommand{\msf}[1]{\mathsf{#1}}
\newcommand{\mbf}[1]{\mathbf{#1}}
\newcommand{\mfrak}[1]{\mathfrak{#1}}
\newcommand{\mrm}[1]{\mathrm{#1}}

\newcommand{\diff}{\mathop{}\mathopen{}\mathrm{d}}

\newcommand{\ii}{\mathrm{i}}

\newcommand{\tend}[1]{\underset{#1}{\longrightarrow}}

\theoremstyle{plain}
\newtheorem{theorem}{Theorem}[section]
\newtheorem{proposition}[theorem]{Proposition}
\newtheorem{lemma}[theorem]{Lemma}
\newcommand{\defi}{\overset{(\mrm{def})}{=}}
\newcommand{\vertiii}[1]{{\left\vert\kern-0.25ex\left\vert\kern-0.25ex\left\vert #1 
		\right\vert\kern-0.25ex\right\vert\kern-0.25ex\right\vert}}

{
	\theorembodyfont{\normalsize}
	\theoremsymbol{\ensuremath{\square}}
	\newtheorem*{proof}{Proof}
	
}
\newtheorem{definition}[theorem]{Definition}
\newtheorem{corollary}[theorem]{Corollary}
{	
	\theorembodyfont{}
	\newtheorem{remark}[theorem]{Remark}
}

{	
	\theorembodyfont{\normalsize}
	
}
\newtheorem{assumptions}[theorem]{Assumptions}

\begin{document}

	\title{Asymptotics of the partition function for $\beta$-ensembles at high temperature \thanks{This project has received funding from the European Research Council (ERC) under the European Union
			Horizon 2020 research and innovation program (grant agreement No. 884584).}}
	\author{Charlie Dworaczek Guera\footnote{ENSL, CNRS,  France \newline
			\textit{email:} charlie.dworaczek@ens-lyon.fr}}
	\date{}
	\maketitle
	
	\begin{abstract}We consider the real $\beta$-ensemble (or 1D log-gas) of dimension $N$ in the high-temperature regime, \textit{i.e.} where the inverse temperature $\beta$ scales as $N\beta=2P$, with $P$ a fixed positive parameter. We establish the large-$N$ asymptotic expansion at all orders of the partition function:
		\begin{equation*}
			\mcal{Z}_N[V]=\int_{\R^N}\prod_{i<j}^{N}\left |x_i-x_j\right|^{\frac{2P}{N}}\cdot\prod_{i=1}^{N}e^{-V(x_i)} \diff x_i
		\end{equation*}
		for $V(x)=x^2+\phi(x)$ with $\phi$ a bounded smooth function, and identify the first two terms of this expansion.
		
		In this regime, the energy no longer dominates the entropy, as in the fixed-$\beta$ case,  but rather scales at the same order in $N$. Consequently, at large $N$, the system is macroscopically described by the so-called\textit{ thermal equilibrium measure} which is supported on the entire real line. 
		
		Our proof relies on the loop equations method, previously applied in the fixed-$\beta$ setting in \cite{BoG1,BoG2}, and provides the first example in which this approach can be successfully implemented using the thermal equilibrium measure. This requires a detailed understanding of both the thermal equilibrium measure and the associated master operator, an unbounded differential operator, leading to several new analytical challenges.
		
		In this setting, we carry out a technically involved analysis to obtain precise estimates for the inverse of the master operator in suitable functional norms. In addition we establish, through subtle operator arguments, a crucial continuity property of the equilibrium density with respect to the potential dependence. These two results constitute the main novelties of the paper and allow us to exhibit a new class of multiple integrals for which such an expansion can be obtained, while providing a deeper understanding of the thermal equilibrium measure and its properties.
	\end{abstract}
	
	\tableofcontents
	
\section{Introduction}
\subsection{Setting of the problem}
Let $P>0$ and $V:\R\rightarrow\R$ be a function growing sufficiently fast at infinity, see Assumptions \ref{assumptions}. The real $\beta$-ensemble of dimension $N$ at high temperature is the particle system on $\R$, $\{x_i\}_{i=1}^{N}$ with the following distribution:
\begin{equation}\label{def:mesureparticules}
	\diff\P^{V}_N(\underline{x})\defi p_N^{V}(\underline{x})\diff x_1\dots \diff x_N\hspace{0,5cm}\text{ with }\hspace{0,5cm}p_N^{V}(\underline{x})\defi \dfrac{1}{\mcal{Z}_N[V]}\prod_{i<j}^{N}\left |x_i-x_j\right|^{\frac{2P}{N}}\cdot\prod_{i=1}^{N}e^{-V(x_i)}
\end{equation}
where $\mcal{Z}_N[V]>0$ is the partition function that ensures that $\P_N^{V}$ is a probability measure on $\R^N$, namely:
\begin{equation}\label{def:partfunction}
	\mathcal{Z}_N[V]\defi\int_{\R^N}\prod_{i<j}^{N}\left |x_i-x_j\right|^{\frac{2P}{N}}\cdot\prod_{i=1}^{N}e^{-V(x_i)} \diff x_i.
\end{equation}
Here, the factor $2$ in the two-body interaction is irrelevant and just makes the equations look "nicer". The main goal of this article is to establish the existence of the large-$N$ asymptotic expansion of the free energy $\log\mcal{Z}_N[V]$ under some assumptions on $V$ namely that there exists $(c_i)_{i\geq-1}$ such that for all $M\geq-1$:
\begin{equation}\label{introeq:goal}
	\log 	\mathcal{Z}_N[V]=\sum_{i=-1}^{M}\dfrac{c_i}{N^{i}}+O\left(N^{-(M+1)}\right) 
\end{equation}

Obtaining such large-$N$ all order expansions for the partition functions of Gibbs measures like $\P_N$ is interesting for several reasons. It translates into expansions for the volume of certain high-dimensional convex sets \cite[Corollary 1.9]{guera2025clt}, allows one to define certain quantities in quantum field theory \cite{borot2016asymptotic,dworaczek2025equilibrium} or to obtain topological data of the system in the expansion. Indeed, the entropy of the equilibrium measure, the Lebesgue measure \cite[Corollary 1.1]{leble2017large} or the number of connected components \cite[Theorem 1.4]{BoG2} of the support can be found in the coefficients in the expansion. Furthermore, obtaining the form of the expansion (namely in powers of $N^{-1}$ here) is interesting and is still unknown in the case of Coulomb gases, see \cite[Section 2.3]{serfaty2023gaussian}.

To do so, we rely on the analysis of a certain tower of equations indexed by integers $n\geq1$ called the \textit{loop equations}, or sometimes Dyson-Schwinger equations/Ward identities. This technique, which we describe further below, was first introduced in \cite{Ambj_rn_1992,Ambj_rn_1993} and later developed in \cite{BoG1,BoG2}.

It is well known that the leading order of this type of integral is obtained by using large deviations arguments for the law of the empirical measure $\mu_N$:
$$\mu_N\defi\dfrac{1}{N}\sum_{i=1}^{N}\delta_{x_i},$$
see for example \cite[Proof of Theorem 2.6.1]{anderson2010introduction} for $\beta$-ensembles with $\beta$ fixed, \cite[equation (1.9)]{leble2017large} in the Coulomb/Riesz case or \cite[(2.4.2)]{borot2016asymptotic} for a sinh-interaction. In the context of $\beta$-ensembles at high temperature, a large deviation principle (LDP) for the law of $\mu_N$ was shown in \cite{Garcia}. The rate function is given by $\mcal{E}_V-\inf\limits_{\mu\in\mcal{M}_1(\R)}\mcal{E}_V(\mu)$ where:
\begin{equation}\label{def:EV}
	\mcal{E}_V(\mu)\defi \int_\R V(x)\diff\mu(x)-2P\iint_{\R^2}\log|x-y|\diff\mu(x)\diff\mu(y)+\int_\R\log\left(\dfrac{\diff\mu(x)}{\diff x}\right )\diff\mu(x).
\end{equation}
This functional can be shown to admit a unique minimizer $\mu_V$ that we call the \textit{equilibrium measure} and which is sometimes referred to as the \textit{thermal equilibrium measure} in the context of Coulomb gases \cite{armstrong2021local,armstrong2022thermal,serfaty2023gaussian,padilla2025poisson}. This minimizer is absolutely continuous with respect to the Lebesgue measure, and its density $\rho_V$, called the \textit{equilibrium density}, is supported on the whole real line, see \cite[Lemma 3.2]{GMToda} for a proof and \cite{AllezBouchaudGuionnet} for plots of $\rho_V$.  As a consequence of \cite{Garcia,GMToda}, we have:
$$\lim\limits_{N\rightarrow\infty}\dfrac{1}{N}\log \mcal{Z}_N[V]=-\mcal{E}_V(\mu_V).$$
As the minimizer of the functional $\mcal{E}_V$, it is known that there exists a constant $\lambda_V\in\R$ such that $\rho_{V}$ satisfies the following equation:
\begin{equation}\label{eqmeasure:charac}V(x)-2P\int_\R\log|x-y|\rho_V(y)\diff y+\log\rho_V(x)=\lambda_V\hspace{1cm}\text{}x-ae
\end{equation}
(see \cite[Lemma 3.2]{GMToda}). In \eqref{eqmeasure:charac}, it is the term $\log\rho_{V}$ that forces $\rho_V$ to be supported on the whole real line, since the two other terms in the LHS remain locally bounded. This is a fundamental difference from the usual $\beta$-ensembles, where the equilibrium measure is compactly supported. This equation can be rewritten as
\begin{equation}\label{eqmeasure:exp}
	\forall x\in\R,\quad\rho_V(x)=\exp\Big(-V(x)-2PU^{\rho_V}(x)+\lambda_V\Big),\hspace{0,7cm}U^{\rho_V}(x)\defi -\displaystyle\int_\R\log|x-y|\rho_V(y)\diff y.
\end{equation}
One can observe that $U^{\rho_{V}}$ behaves like $-\log|x|$ at infinity, see \cite[Lemma 2.4]{DwoMemin}. Hence, assuming that $V$ grows fast enough at infinity, unlike the constant $\beta$ case where the equilibrium measure is compactly supported, the equilibrium density  $\rho_V$ here, "only" decays exponentially fast at infinity.
This measure was first studied in \cite{AllezBouchaudGuionnet} with quadratic $V$, where it was shown to interpolate between the Gaussian measure and the semi-circle law, see Figure 2 in this article. In the case $V_G(x)\defi x^2/2$ the equilibrium density $\rho_{V}$, which actually represents \textit{Askey-Wimp-Kerov} distribution in that case \cite{askey1984associated}, can also be made explicit:
$$\rho_{V_G}(x)=\dfrac{e^{-\frac{x^2}{2}}}{\sqrt{2\pi}}\dfrac{1}{|\widehat{f}_\alpha(x)|^2},\hspace{2cm}\widehat{f}_\alpha(x)\defi  \sqrt{\dfrac{P}{\Gamma(P)}}\displaystyle\int_0^{+\infty}t^{P-1}e^{-\frac{t^2}{2}+\ii xt}\diff t.$$
The LDP \cite{Garcia} establishes that $\mu_V$ is the almost-sure weak limit of $\mu_N$ from which we can deduce that for a large class of test function $\phi$ we have:
\begin{equation}\label{introeq: ftend0}
	\braket{\phi}_{\Delta\mu_N}^{V}\tend{N\rightarrow\infty}0,\hspace{1cm}\Delta\mu_N=\mu_N-\mu_V
\end{equation}
where we defined for any function $f_n$ depending on $n$ real variables and (possibly random) signed measures on $\R$ $\nu_i$:
\begin{equation}\label{eq:braket}
	\Braket{f_n}_{\nu_1\otimes\hdots\otimes\nu_n}^V\defi \mbb{E}_N^{V}\left[\int_{\R^{n}} f_n(x_1,\dots,x_n) \prod_{i=1}^n\diff\nu_i(x_i)\right].
\end{equation}
The analysis of loop equations allows one to generalize \eqref{introeq: ftend0} for smooth test-functions $\phi$ and general $V$; and obtain the whole asymptotic expansion of linear statistics up to order $M$ for all $M\geq0$:
\begin{equation}\label{introeq: fasymptotic}
	\braket{\phi}_{\Delta\mu_N}^{V}=\sum_{i=1}^{M}\dfrac{d_i^{V}(\phi)}{N^{i}}+O\left (N^{-(M+1)}\right ).
\end{equation}
This method (as will be explained in Subsection \ref{subsec:strategy}) relies on fine properties of a certain operator $\Xi$, called the \textit{master operator}, defined, for sufficiently smooth $\phi$, by:
\begin{equation}\label{def:masteroperator}
	\Xi[\phi]\defi \phi'+\big(\log\rho_V\big)'\phi+2P\left(\mcal{H}[\phi\rho_V]-\int_\R\mcal{H}[\phi\rho_V](y)\diff\mu_V(y)\right)
\end{equation}
where $\mcal{H}$ denotes the \textit{Hilbert transform}, which is defined for general functions by
\begin{equation}\label{def:hilberttransform}
	\mcal{H}[f](x)\defi\displaystyle\fint_\R\dfrac{f(y)}{y-x}\diff y
\end{equation}
where $\fint$ denotes the principal value integral. This operator was shown to be invertible in a certain functional space \cite[Theorem 6.9]{DwoMemin} and one of the main steps to establish \eqref{introeq: fasymptotic} is to show that $\Xi^{-1}$ is continuous on the functional spaces $W_n^{p}(\R)$ and  for all $n\geq1$ and $p\in\{2,\infty\}$, defined by:
\begin{equation}\label{def:sobolev spaces}
	W_n^{p}(\R^{k})\defi\left\{f\in L^p(\R^{k}),\,\partial_1^{a_1}\hdots\partial_k^{a_k}f\in L^{p}(\R^{k}),\,\sum_{l=1}^{k}a_l\leq n\;\text{with }a_l\in\N\right\}.\hspace{0,4cm}
\end{equation}
In the special case $p=2$, we use the notation $H^{n}(\R^{k})\defi W_n^{2}(\R^{k})$. These spaces (for general $p$) are endowed with the norm:
\begin{equation}\label{def:sobolev norms}	\|f\|_{W_n^{p}(\R^{k})}\defi\max\left\{\|\partial_1^{a_1}\hdots\partial_k^{a_k}f\|_{L^{p}(\R^{k})},\,a_l\in\N\;\text{with} \sum_{l=1}^{k}a_l\leq n\right\}.
\end{equation}
Combining such continuity properties of $\Xi^{-1}$ and basic concentrations arguments is actually enough to conclude the expansion \eqref{introeq: fasymptotic} of linear statistics at all orders for regular enough test functions.
Now, specializing $V=V_{G,\phi}\defi V_G+\phi$ where $\phi$ is such that \eqref{introeq: fasymptotic} holds then setting $V_{G,\phi,t}\defi V_G+t\phi$ for all $t\in[0,1]$
\begin{equation}\label{eq:introinterpolation}
	\log\mcal{Z}_N[V_{G,\phi}]=\log \mcal{Z}_N[V_G]-N\int_{0}^{1}\Braket{\phi}^{V_{G,\phi,t}}_{\mu_N}\diff t.
\end{equation}
Recalling that the expansion of $\log \mcal{Z}_N[V_G]$ is easy to obtain thanks to Mehta's explicit formula \eqref{partition,function gaussian}, obtaining the expansion of the LHS in \eqref{eq:introinterpolation} boils down to obtaining the asymptotic expansion of the integral in the RHS. Inserting the expansion \eqref{introeq: fasymptotic} inside the integral, assuming that every term as well as the remainder is $t$-integrable, one obtains the desired expansion \eqref{introeq:goal}. By looking at the expressions of the coefficients $d_i(\phi)^{V_{G,\phi,t}}$ and the form of the remainder, it can be shown that the continuity of $t\mapsto\rho_{V_{G,\phi,t}}$ in a sufficiently nice functional space is enough to justify the integrability of each term in \eqref{introeq: fasymptotic} and conclude.
\subsection{Assumptions}
Although Theorem \ref{thm:partitionfunction} is proved only for $V=V_{G,\phi}$ with $\phi$ smooth and bounded, most of our arguments extend to a general potential $V$. Unless stated otherwise, we impose the following assumptions on $V$ throughout the paper:
\begin{assumptions}\label{assumptions}
	The potential $V$ satisfies:
	\begin{enumerate}[label=(\roman*)]
		\item\label{assumption1} $V\in \mathcal{C}^\infty(\R,\R)$,
		\item \label{assumption2}$V(x)\underset{|x|\rightarrow+\infty}{\longrightarrow}+\infty$ and $|V'(x)|\underset{|x|\rightarrow+\infty}{\longrightarrow}+\infty$,
		\item\label{assumption3} The measure $\mu_V$ satisfies the Poincaré inequality, \textit{i.e.} there exists $C_{\mrm{Poinc}}>0$ (depending on $V$ and $P$) such that for all $f\in\mcal{C}_c^1(\R)$
		\begin{equation}\label{Poincaré:inequality}
			\mathrm{Var}_{\mu_V}(f)\defi \int_\R \left(f(x)-\int_\R f(y) \diff\mu_V(y)\right)^2\diff\mu_V(x) \leq C_{\mrm{Poinc}}\int_\R f'(x)^2 \diff\mu_V(x)\,.
		\end{equation}
		\item \label{assumption4}For all polynomial $Q\in\R[X]$ and $\alpha>0$, all $p\geq0$, $Q\left (V^{(p)}(x)\right )e^{-V(x)}=\underset{|x|\to \infty}{o}(x^{-\alpha})$\,.
		\item \label{assumption5} The function $x\mapsto V'(x)^{-2}$ is integrable at infinity, and $\dfrac{V^{(k)}(x)}{V'(x)}=\underset{|x|\rightarrow\infty}{O}(1)$ for $k\geq2$.
	\end{enumerate}
\end{assumptions}

Assumption \textit{(i)} is needed in order to analyze the loop equations. Obtaining the expansion of linear statistics up to $o(N^{-k})$ requires controlling building blocks appearing in the $i$-th loop equation for each $i\in\llbracket1,2k\rrbracket$. The further along these terms appear in the tower of equations, the more derivatives and iterations of $\Xi^{-1}$ ($\Xi$ being defined in \eqref{def:masteroperator}) they involve. Since the regularity of $\Xi^{-1}[\phi]$ depends both on the regularity of $\phi$ and on that of $V$, assuming $V\in\mcal{C}^{\infty}$ allows one to obtain the expansion \eqref{introeq: fasymptotic} at all orders. It also ensures that $\rho_V$ is smooth, which will be necessary to show the continuity of the different operators involved in the loop equations.

Assumption \textit{(ii)} ensures that $\mcal{Z}_N[V]$ is well-defined. Indeed, the assumption on $V'$ implies that $V$ grows faster than linearly at infinity. A further consequence, by Lemma \ref{lem:regularitedensite}, is that $\rho_V$ decays exponentially fast at infinity. The condition that $V'$ goes to infinity is also needed to ensure that $\Xi^{-1}[\phi]^{(k)}(x)\tend{|x|\rightarrow\infty}0$ for all $k\geq0$ and for bounded smooth functions $\phi$.

Assumption \textit{(iii)}, together with \textit{(i) } and \textit{(ii)}, implies that $\Xi$ defined in \eqref{def:masteroperator}, is invertible; see \cite[Theorem 6.9]{DwoMemin}. The authors showed that this condition is not overly restrictive \cite[Remark 1.3, Proposition 2.6]{DwoMemin}. Indeed, for any potential of the form $V=V_{\mrm{conv}}+\phi$, where $V_{\mrm{conv}}$ is strictly convex outside of a compact set and $\phi$ is bounded, $\mu_{V}$ satisfies the Poincaré inequality.

Assumption \textit{(iv)} ensures that the equilibrium density $\rho_V$ and all of its derivatives decay exponentially fast at infinity.

Assumption \textit{(v)} is used to prove that $\Xi^{-1}$ is continuous with respect to the norms \eqref{def:sobolev norms}. Indeed, when differentiating $\Xi^{-1}[\phi]$, for $\phi$ a smooth function, quantities behaving at infinity like $V^{(k)}(x)V'(x)^{-1}$ naturally arise. On the other hand, in our approach we need to integrate some functions that behave like $V'(x)^{-2}$ at infinity.

These conditions are satisfied, for example, for every $V$ of the form: $$x\mapsto a_nx^{2n}+\phi(x),\hspace{1cm}\text{or}\hspace{1cm}x\mapsto\dfrac{e^{\gamma x}+e^{-\gamma x}}{\alpha}$$
where $n>0,a_n>0,\phi^{(k)}\text{ bounded }\forall k\in\N$ and $\alpha>0,\gamma\in\R$. On the other hand, a potential like $V(x)=e^{x^2}$ violates assumption \textit{\ref{assumption5}}, therefore it does not fit into our analysis.

\subsection{Main results}

In order to state the next result, we recall that $\Delta\mu_N\defi \mu_N-\mu_V$ and that $H^n(\R^k)$ denotes the $n$-th Sobolev space as defined in \eqref{def:sobolev spaces}.

\begin{theorem}[Asymptotic expansion of linear statistics]\label{thm:correlators}
	Assume $V$ satisfies assumptions \ref{assumptions}, then for all smooth functions $\phi\in H^{r}(\R^k)$ for $r>0$ (depending on $K$) big enough, there exists a unique sequence $(b_i)_{i\geq\lceil k/2\rceil}$ depending on $V$, $\phi$ and $P$ such that for all $K\in\N$:
	$$\boxed{\Braket{\phi}_{\otimes^{k}\Delta\mu_N}^V=\sum_{i=\lceil k/2\rceil}^{K}\dfrac{b_i}{N^i}+O\left(N^{-(K+1)}\right).}$$
\end{theorem}

Our goal is to establish the existence of an expansion for $\mcal{Z}_N[V_{G,\phi}]$ where $V_{G,\phi}(x)\defi x^2/2+\phi(x)$ and $\phi\in H^{\infty}(\R)\defi\bigcap_{n\geq1}H^{n}(\R)$. As described in \eqref{eq:introinterpolation}, the idea is to deduce this from an integration of the expansion of $\Braket{\phi}_{\Delta\mu_N}^{V_{G,\phi,t}}$ for all  $t\in[0,1]$ where $V_{G,\phi,t}(x)= x^2/2+t\phi(x).$
Making this step rigorous requires the following continuity result with respect to $t$. We recall that $W_n^{\infty}(\R^{k})$-norm is given by $\|f\|_{W_n^\infty(\R^{k})}\defi \underset{i\in\llbracket0,n\rrbracket}{\max}\|f^{(i)}\|_{\infty}$.

\begin{theorem}\label{thm:conteqdensity}
	Assume $V$ satisfies \ref{assumptions} then for all $n\in\N$, $t'\in[0,1]$ and all $\phi\in H^{\infty}(\R)$, we have:
	$$\boxed{\|\rho_{V_{\phi,t}}-\rho_{V_{\phi,t'}}\|_{W_n^{\infty}(\R)}\underset{t\rightarrow t'}{\longrightarrow}0}$$where $V_{\phi,t}:x\mapsto V(x)+t\phi(x)$ for any $t\in[0,1]$. Furthermore, for all $x\in\R$, $t\mapsto\rho_{V_{\phi,t}}(x)\in\mcal{C}^\infty(\R)$ and satisfies the following integro-differential equation for all $t\geq0$ and $x\in\R$:
	$$\partial_t\rho_{V_{\phi,t}}(x)=\left(-\phi(x)+\int_\R\phi(s)\rho_{V_{\phi,t}}(s)\diff s \right) \rho_{V_{\phi,t}}(x).$$
\end{theorem}
Theorem \ref{thm:correlators} together with Theorem \ref{thm:conteqdensity} allow us, thanks to \eqref{interpolation} to deduce the following which constitutes the main theorem of the present article.

\begin{theorem}[Asymptotic expansion of the partition function]\label{thm:partitionfunction}
	Let $\phi\in H^{\infty}(\R)$, then there exists a unique sequence $(c_i)_{i\geq0}\in\R^\N$ depending on $\phi$ and $P$, such that for all $K\in\N$:
	$$\boxed{\dfrac{1}{N}\log \mcal{Z}_N\left[V_{G,\phi}\right]=\sum_{i=0}^{K}\dfrac{c_i}{N^i}+O\left(N^{-(K+1)}\right) .}$$
	The leading term is given by $c_{0}\defi-\mcal{E}_{ V_{G,\phi}}(\mu_{ V_{G,\phi}})$ \textit{i.e.}:
	\begin{equation}\label{eq:leadingtermin}
		c_0= -\int_\R V_{G,\phi}(x) \diff\mu_{V_{G,\phi}}(x)+P\iint_{\R^2}\log|x-y|\diff\mu_{V_{G,\phi}}(x)\diff\mu_{V_{G,\phi}}(y)+\mrm{Ent}\left[\mu_{V_{G,\phi}}\right].\nonumber
	\end{equation}
	where $$\mrm{Ent}\left[\mu_{V_{G,\phi}}\right] =-\int_{\R}\log\left(\dfrac{\diff\mu_{V_{G,\phi}}(x)}{\diff x}\right)\diff\mu_{V_{G,\phi}}(x).$$
	The subleading term $c_1$ is given in terms of $\Xi^{-1}_t$ the inverse of the master operator (defined in \eqref{def:masteroperator}) associated with the potential $V_{G,\phi,t}$. It can be written as
	\begin{equation}
		c_1\defi g_{1}(P)-P\int_{0}^{1}\left[\Braket{\widetilde{\Xi_t^{-1}}[\phi]'}_{\mu_{V_{G,\phi,t}}}+\Braket{\Theta^{(2)} \circ\widetilde{\Xi_{t,1}^{-1}}\left[\partial_2\mcal{D}\circ\widetilde{\Xi_t^{-1}}[\phi]\right] }_{\mu_{V_{G,\phi,t}}}\right] \diff t.
	\end{equation}
	where $g_1(P)$ is given in \eqref{eq:g1}.
\end{theorem}
Above, $\Theta^{(2)}$ and $\mcal{D}$ are explicit operators given in Section \ref{section:loopequations} while $$\Xi_t^{-1}[\phi](x)=\dfrac{1}{\rho_{V_{G,\phi,t}}(x)}\displaystyle\int_{x}^{+\infty}\mcal{T}_t[\phi](y)\diff\mu_{V_{G,\phi,t}}(y)$$ where $\mcal{T}_t$ is an explicit kernel operator given in \eqref{eq:inverse fredholm}. Given an operator $\mcal{S}$, the symbol $\mcal{S}_1$ is also defined in Sections \ref{section:loopequations} and \ref{section4}.

We emphasize that Theorem \ref{thm:partitionfunction} cannot be deduced from \cite[Proposition 1.2]{BoG1}. Indeed, taking $\beta=2P/N$ and then the large $N$-limit is different from taking the large $N$ expansion of \cite{BoG1} valid for fixed $\beta$ and then inserting $\beta=2P/N$ in the large $N$ expansion. Also some key points of the proof of the mentioned result are different in the high-temperature regime (different scalings, support of the equilibrium measure, form of the master operator or also the $V$-dependence of the equilibrium density).
\subsection{Connection with the literature and motivations}
\paragraph{Comparison with the fixed-temperature regime}
The fixed-temperature regime of $\beta$-ensembles is given by the following distribution:
\small{\begin{equation*}
	\diff\Q^{(\mrm{fixed})}_N(\underline{x})\defi q_N^{(\mrm{fixed})}(\underline{x})\diff x_1\dots \diff x_N\hspace{0,4cm}\text{ with }\hspace{0,4cm}q_N^{(\mrm{fixed})}(\underline{x})\defi \dfrac{1}{\mcal{Z}_N^{(\mrm{fixed})}[W]}\prod_{i<j}^{N}\left |x_i-x_j\right|^{\beta}\cdot\prod_{i=1}^{N}e^{-\beta NW(x_i)}.
\end{equation*}}
\normalsize
Above $W:\R\rightarrow\R$ is a continuous function going to infinity sufficiently fast and $\beta>0$. Setting, $\beta=2P/N$ and $V=2PW$, one recovers the distribution $\mbb{P}_N$ \eqref{def:mesureparticules}. Since in the latter case the parameter $\beta$, interpreted as the inverse temperature of the gas, is very small, it justifies why $\mbb{P}_N$ is referred to as the high-temperature regime.

A great deal is known about $\mbb{Q}_N$, owing to two decades of extensive study. First, LDPs for the empirical measure and the largest particle have been established in \cite{arous1997large,arous2001aging}, central limit theorems were proven in \cite{johansson,BoG1,Shc1,bekerman2018clt,LambertLedouxWebb,bourgade2022optimal, bekerman2018mesoscopic,peilen2024local,angst2024sharp,guera2025clt}, a full asymptotic expansion of the partition function \cite{ErMc,BoG1,BoG2} and local laws \cite{bourgade2014local,bourgade2014universality,bourgade2022optimal,peilen2024local} have been shown. For $\beta=2$ and $W$ polynomial, the asymptotic expansion of the partition function $\mcal{Z}_N^{(\mrm{fixed})}[W]$ has the form $$\dfrac{1}{N^{2}}\log\mcal{Z}_N^{(\mrm{fixed})}[W]_{|\beta=2}=\sum_{g\geq0}\dfrac{c_g}{N^{2g}}$$ where the previous equality has to be understood in the sense of an asymptotic expansion. The coefficients $\left(c_{g}\right)_{g\geq0}$ correspond to enumerations of maps and, more generally, the asymptotic expansion of $\log \mcal{Z}_N^{(\mrm{fixed})}[W]$ gives information on the enumerations of graphs embedded in surfaces \cite{mulase2005non,marino2014lectures}.

One of the main differences between $\mbb{P}_N$ and $\Q^{(\mrm{fixed})}_N$ is that, asymptotically the particles are spread over the entire real line in the high-temperature regime (the support of $\mu_V$ being $\R$) whereas the limiting measure $\mu_{(\mrm{fixed})}^{W}$ of $\mu_N$ under $\Q^{(\mrm{fixed})}_N$ is typically compactly supported (the semi-circle law in the case of the quadratic potential). This fact makes truncation arguments, where one makes the approximation $$\mbb{E}_N\left[\int_{\R} f(y) \diff\mu_N(y)\right]= \mbb{E}_N\left[\int_{\R}f(y)\psi(y)\diff\mu_N(y)\right] +\varepsilon_N\left(f,\psi\right) $$ where $f$ is an unbounded function and $\psi$ a smooth function whose support contains the asymptotic support, more difficult to implement in the setting of $\P_N$. It is one of the main obstacles to establishing Theorem \ref{thm:correlators} for general (unbounded) smooth functions $\phi$ and thus obtain Theorem \ref{thm:partitionfunction} for general smooth potential $V$.

About the master operator, establishing invertibility and controls (see \cite[Lemma 4.1]{guera2025clt} for example) in the fixed-temperature regime, \textit{i.e.} for:
$$\Xi^{(\mrm{fixed})}[\phi](x)\defi-V'(x)\phi(x)+\int_\R\dfrac{\phi(x)-\phi(y)}{x-y}\diff \mu_{(\mrm{fixed})}^{W}(y)$$ 
is also much simpler compared to $\Xi$ (defined in \eqref{def:masteroperator}). It relies on the explicit inversion formula due to Tricomi \cite{tricomi1957integral}. The inversion of $\Xi$, obtained in \cite{DwoMemin}, is considerably more involved requiring Hilbertian techniques. An explicit and tractable expression for $\Xi^{-1}$ is much harder to obtain, owing to the additional differential term. Thus, controls in $W_p^{n}$-norms (defined in \ref{def:sobolev norms}) are much more difficult to obtain. These continuity-estimates for $\Xi^{-1}$ are one of the key ingredients to make the loop equation analysis method work.

Finally, one of the main differences between the high and fixed-temperature regime is the integration step, \textit{i.e.} to deduce Theorem \ref{thm:partitionfunction} from Theorem \ref{thm:correlators} and \eqref{eq:introinterpolation}. This requires a precise understanding of the $t$-dependence of the coefficients $b_i$ as well as the remainder in Theorem \ref{thm:correlators} when considering $V_t=t V+(1-t)V_G$ for a general potential $V$. It can be shown that showing a continuity-type result for the equilibrium density in the potential $V$ dependence is enough to conclude. While this step is technically involved in the high-temperature regime, it is straightforward to show in the fixed-$\beta$ case. This is due to the "linear"-dependence of $\mu_{\mrm{(fixed)}}^{W}$ on $W$, namely: setting $W_t=t W_1+(1-t)W_2$ with $W_1$ and $W_2$ chosen such that $\mu_{(\mrm{fixed})}^{W_i}$ share the same support, then:
\begin{equation}\label{eq:interpolationfixed}
	\mu_{(\mrm{fixed})}^{W_t}=t\mu_{(\mrm{fixed})}^{W_1}+(1-t)\mu_{(\mrm{fixed})}^{W_2}.
\end{equation}
This is proven in \cite[Lemma 5.1]{BoG1} and is due to the absence of the entropy in the energy functional that $\mu_{(\mrm{fixed})}^{W}$ minimizes, see \cite[Theorem 1.1]{BoG1}. In the high-temperature regime, the analogue of \eqref{eq:interpolationfixed} is not true anymore for $\rho_{V_t}$ (as can be seen from the asymptotic behavior of $\rho_{V_t}$ at infinity for example) and the $t$ dependence is much more complicated to capture.
\paragraph{The high-temperature regime}
The study of the high-temperature regime started with the pioneering works \cite{cepa1997diffusing,bodineau1999stationary,AllezBouchaudGuionnet}. More recently, LDPs for the empirical measure \cite{Garcia}, and for the largest particle at high \cite{dworaczek2025large} and intermediary temperature \cite{pakzad2020large} have been shown. Central limit theorems for linear statistics have been established in \cite{hardy2021clt,NakanoTrinh,DwoMemin,mazzuca2023clt} at a great level of generality and the Poissonian limiting microscopic behaviors \cite{benaych2015poisson,pakzad2018poisson,NakanoTrinh,nakano2020poisson,lambert2021poisson,padilla2024emergence} were shown in this regime both in the bulk and at the edge. Finally, in \cite{forrester2021classical}, the authors analyzed the loop equations to deduce the moments of the subdominant correction of the equilibrium measures corresponding to Gaussian, Laguerre and Jacobi ensembles.

\paragraph{Asymptotic analysis of multiple integrals}
From the point of view of asymptotic analysis of multiple integrals, $\beta$-ensembles provide, \textit{via} its partition function, a good non-trivial (with an interaction between the integration variables) example of a $N$-fold multiple integrals whose AE is known at all orders \cite{BoG1,BoG2} at least in the case of a smooth potential $V$. This was later generalized in \cite{borot2015large}, to allow for general analytic $p$-interaction with $p\geq2$, between the integration variables. In the so-called varying weight situation and with a \textit{sinh}-interaction, the same authors managed to obtain the asymptotics of the partition function up to $o(1)$ in \cite{borot2016asymptotic}. Thanks to the developd tools, the authors managed to give a check on the so-called Lukyanov's conjecture in \cite{dworaczek2025equilibrium}. Generalizations of these results on contours in the complex plane were established in \cite{CourteautJohansson2025,guionnet2024asymptotic}. A motivation to extend the first AE result of \cite{BoG1} to more general domain of integration, interaction and scaling in $N$ in front of the potential is to developp tools for the so-called \textit{quantum separation of variables method}, see \cite[Section 1.5]{borot2016asymptotic}, an alternative to the \textit{algebraic Bethe Ansatz}.
In the so-called Freud weight context, \textit{i.e.} when $V=|.|^{p}$, the $\beta$-ensembles are related to the geometry of the so-called \textit{Schatten-balls}. For small $p>0$, the potential $V$ is very singular and most of the techniques break down. However, AE of the log-partition function up to $o(N)$ were established in  \cite{dadoun2023asymptotics,guera2025clt,sonnleitner2025next}. A certain low-temperature regime of $\beta$-ensembles (with slightly different conventions, \textit{i.e.} $W$ scaled as $N\sqrt{\beta}W$ instead of $\beta NW$) with $\beta\rightarrow\infty$ has been studied in \cite{MR3165804}. Theorem \ref{thm:partitionfunction} completes the picture by filling the gaps needed to perform the loop equations method.

\paragraph{Link with integrable systems}
The study of the $\beta$-ensembles at high-temperature has attracted a lot of attention recently since links were discovered with integrable systems, such as the famous \textit{classical Toda chain} \cite{toda1967vibration}. The integrable structure of this system, namely the existence of a sufficient number of conserved quantities, can be established by the existence of a so-called \textit{Lax matrix}, whose spectrum is invariant under the dynamics. At long times, the model does not \textit{thermalize}, \textit{i.e.} it does not reach thermal equilibrium but is rather described by a more sophisticated probability measure called the \textit{Generalized Gibbs Ensemble} (GGE) \cite{jaynes1957information}. This is due to the existence of a set of locally conserved quantities, which highly constrains the dynamics. In the context of the Toda chain, the GGE has been studied in \cite{Spohn1} and a link was established with the Gaussian $\beta$-ensembles. In the case of a Gaussian potential, it was shown that the distribution of the Lax matrix under the GGE was similar to the law of the tridiagonal representation of the Gaussian $\beta$-ensembles of Dumitriu and Edelman. This link was explored in \cite{GMToda,mazzuca2023large} for more general potentials via large deviation techniques. Theorem \ref{thm:partitionfunction} can be seen as a toy model to study more involved integrals appearing in this literature such as the partition function of the GGE of the Toda and Calogero fluid, see \cite[(9.34), (11.51)]{spohn2023hydrodynamic}:
$$\mcal{Z}_N^{\mrm{sys}}[V]=\dfrac{2^{N}}{N!}\int_{\R^{N}}\diff^{N}\underline{\lambda}\cdot\prod_{i=1}^{N}e^{-V(\lambda_i)}\cdot\prod_{j=1}^{N}K_0\left(e^{-N\nu}Y_j^{\mrm{sys}}(\underline{\lambda})\right) $$
where the zero-order Bessel function of the second kind is given by $K_0(x)\defi\int_0^{+\infty}e^{-x\cosh t}\diff t$ and
$$Y_j^{\mrm{sys}}(\underline{\lambda})=\begin{cases}
	\prod_{i\neq j}^{N}|\lambda_i-\lambda_j|^{-1}&\quad\text{if sys}=\mrm{Toda}
	\\\prod_{i\neq j}^{N}\left(1+|\lambda_i-\lambda_j|^{-1/2}\right)^{2}&\quad\text{if sys}=\mrm{Calogero}
\end{cases}.$$
These integrals share many similarities with $\mcal{Z}_N[V]$ in the sense that the leading order should also be given by a competition between the entropy and the energy and the minimizer, \textit{i.e.} the equilibrium measure of both systems should also be very similar to $\mu_V$ in our context.
\paragraph{The thermal equilibrium measure}
The present article shows for the first time that the loop equation method can be used with the thermal equilibrium measure, \textit{i.e.} the minimizer of $\mcal{E}_V$. For general $\beta>0$ (possibly depending on $N$), this measure (when adding a prefactor $\dfrac{2}{N\beta}$ in front of the entropy in \eqref{def:EV}) should provide a better description of the system than the usual equilibrium measure (when there is no entropy term in the functional). However, doing so implies dealing with a $N$-dependent measure and also master operator which is technically very demanding. The goal of this paper is to consider, the simplest case with $N\beta$ constant. To the knowledge of the author, while this measure is well-understood in the context of Coulomb gases \cite{armstrong2022thermal}, there is no good understanding in full generality of this measure quantitatively speaking in one-dimension.
\subsection{Outline of the proof}\label{subsec:strategy}
The strategy used to show Theorem \ref{thm:partitionfunction} is based on the following \textit{interpolation equation} of the form:
\begin{equation*}
	\log\mcal{Z}_N[V_{G,\phi}]=\log \mcal{Z}_N[V_G]-N\int_{0}^{1}\Braket{\phi}^{V_{G,\phi,t}}_{\mu_N}\diff t.
\end{equation*}
On the RHS, the expansion of $\log \mcal{Z}_N[V_G]$ follows from Mehta's formula \cite[17.6.7]{mehta2004random}, see Theorem \ref{thm:asymptfctpartgauss}. Once this identity is obtained, the expansion of the \textit{free energy} of the model $\log\mcal{Z}_N[V]$ follows from the one for the second term in the RHS in the equation above. We now explain  how to prove Theorem \ref{thm:correlators}, \textit{i.e.} how to obtain the AE of $\Braket{\phi}_{\mu_N}^{V}$ for general $V$ satisfying \ref{assumptions}.
\paragraph{Expansion of linear statistics}
The proof is based on the analysis of the loop equations. It consists of a tower of equations that link linear statistics of different orders. The simplest equation is the one at level $1$, which reads for any bounded smooth function $\phi$:
\begin{equation}\label{intro:SDlvl1}
	\Braket{\phi}_{\Delta\mu_N}^V=
	\dfrac{P}{N}\Braket{\big(\Xi^{-1}[\phi]\big)'}_{\mu_V}+\dfrac{P}{N}\Braket{\big(\Xi^{-1}[\phi]\big)'}_{\Delta\mu_N}^V-P\Braket{\mcal{D}\circ\Xi^{-1}[\phi]}_{\Delta\mu_N\otimes \Delta\mu_N}^V
\end{equation}
where $\mcal{D}$ is the operator defined for all $x\neq y$ by $\mcal{D}[\phi](x,y)\defi \dfrac{\phi(x)-\phi(y)}{x-y}$. This equation links the $1$-linear statistic and the $2$-linear statistic $\Braket{\mcal{D}\circ\Xi^{-1}[\phi]}_{\Delta\mu_N\otimes \Delta\mu_N}^V$. Note that these equations have also been derived for the resolvent $W_1(x)$ and two-point functions $W_2(x,y)$ defined as Stieltjes transforms of $\mbb{E}_N[\mu_N]$ in \cite{forrester2021classical} in the Gaussian case. For $W_1(x)$, the equation is asymptotically of Riccati type in the high-temperature, whereas it is merely quadratic in the fixed-temperature regime. The deduction of the asymptotic expansion for linear statistics from the loop equations is based on a so-called \textit{a priori bound}, which we assume for the moment, of the following form:
\begin{equation}\label{eq:borneaprioriintro}
	|\braket{\phi}_{\overset{k}{\otimes}\Delta\mu_N}|\leq C\dfrac{\|\phi\|_k}{N^{k/2}}
\end{equation}
for a norm $\|.\|_k$ that we don't make precise here. Such a bound also holds in the classical setting. Note however, that at high-temperature, one has to deal with a more complex norm which requires integrability conditions on the functions we apply this bound to. To explain how to obtain the first orders, we assume that we know that: $$\Braket{\mcal{D}\circ\Xi^{-1}[\phi]}_{\Delta\mu_N\otimes \Delta\mu_N}^V=\dfrac{\alpha_1(\phi)}{N}+o(N^{-1}) \quad\quad\quad\text{ and }\quad\quad\quad\Braket{\psi}_{\Delta\mu_N}^V=o(1),$$
for some $\alpha(\phi)\in\R$, then \eqref{intro:SDlvl1} allows one to obtain the leading order asymptotic for the 1-statistic: $$\Braket{\phi}_{\Delta\mu_N}^V=\dfrac{1}{N}\Big(P\Braket{\big(\Xi^{-1}[\phi]\big)'}_{\mu_V}-\alpha_1(\phi) P +o(1)\Big)=:\dfrac{\gamma_1(\phi)}{N}+o\left (\dfrac{1}{N}\right ).$$
Assuming now that for $n=2$
\begin{equation}\label{intro:2linear}
	\Braket{\mcal{D}\circ\Xi^{-1}[\phi]}_{\Delta\mu_N\otimes \Delta\mu_N}^V=\sum_{i=1}^{n}\dfrac{\alpha_i(\phi)}{N^i}+o(N^{-n}),
\end{equation}
it is not hard to see that one can iteratively derive the expansion of $\Braket{\phi}_{\Delta\mu_N}^V$ and get:
$$\Braket{\phi}_{\Delta\mu_N}^V=\sum_{i=1}^{n}\dfrac{\gamma_i(\phi)}{N^i}+o(N^{-n}),\hspace{2cm}\gamma_2(\phi)=P\gamma_1\left(\big(\Xi^{-1}[\phi]\big)'\right) -P\alpha_2(\phi).$$
By the same procedure, one can see that the extraction of the asymptotic expansion up to order $n>2$ of the $1$-linear statistics boils down to extracting the one for the $2$-linear statistics.

To achieve that, one needs to investigate the loop equation at level 2, which has the following form for a smooth function $\phi_2$ of $2$ variables:
\begin{equation*}
	\Braket{\phi_2}_{\Delta\mu_N\otimes \Delta\mu_N}^V=\dfrac{1}{N}\Braket{\mcal{U}[\phi_2]}_{\Delta\mu_N}^V+\Braket{\mcal{V}[\phi_2]}_{\overset{3}{\otimes}\Delta\mu_N}^V+\dfrac{1}{N}\Braket{\mcal{W}[\phi_2]}_{\Delta\mu_N\otimes\Delta\mu_N}^V+\dfrac{1}{N}\Braket{\mcal{Y}[\phi_2]}_{\mu_V}
\end{equation*}
where $\mcal{U}$, $\mcal{V}$, $\mcal{W}$ and $\mcal{Y}$ are some operators involving $\Xi^{-1}$ in their definition. From estimate \eqref{eq:borneaprioriintro} that we assumed at the beginning, we know that: $$\Braket{\mcal{U}[\phi_2]}_{\Delta\mu_N}^V=O(N^{-1/2}),\hspace{1cm}\Braket{\mcal{V}[\phi_2]}_{\overset{3}{\otimes}\Delta\mu_N}^V=O(N^{-3/2}),\hspace{1cm}\Braket{\mcal{W}[\phi_2]}_{\Delta\mu_N\otimes\Delta\mu_N}^V=O(N^{-1}).$$
It is straightforward to see that only the term $N^{-1}\Braket{\mcal{Y}[\phi_2]}_{\mu_V}$ yields a non-negligible contribution to the expansion of $\Braket{\phi_2}_{\Delta\mu_N\otimes \Delta\mu_N}^V$ at precision $o(N^{-1})$.  In order to push it up to $o(N^{-2})$, one needs to obtain the asymptotic expansion for the $3$-linear statistics and so on. Each additional order in the asymptotic expansion requires analysing a higher level loop equation. Crucially, only finitely many equations need to be analyzed at each order to capture all contributions, and the estimate allows one to neglect all the other terms and close the analysis. Finally, in order to apply the estimate to neglect the remainders, one needs to show that the operators involved in the SDE preserve enough of the regularity of the function they act on, especially for the inverse of the master operator $\Xi^{-1}$.

\paragraph{Controls on $\Xi^{-1}$}
In this setting, one has to obtain way more subtle controls compared to the constant $\beta$-setting. This is due to the fact that \eqref{eq:borneaprioriintro} involves a more complex norm than just a $L^\infty$-norm (which is the case for the fixed-temperature regime). Moreover, finding a manageable integral representation for $\Xi^{-1}[\phi]$ in order to extract controls out of it, is a highly non-trivial step. This makes the proof of the continuity of $\phi\mapsto\Xi^{-1}[\phi]$ quite technical. Finding such a representation and proving continuity estimates for $\Xi^{-1}$ constitute one of the main technical contributions of this article.

\paragraph{Potential dependence of the equilibrium density}

When integrating the asymptotic expansion of the $1$-linear statistics $\Braket{\phi}^{V_{G,\phi,t}}_{\mu_N}$, one needs to justify that the resulting integrals are finite, \textit{i.e.} that the integrands are integrable. To do so, it is enough to show that the map $t\mapsto\rho_{V_{G,\phi,t}}$, referred to as the \textit{perturbed thermal equilibrium measure} in the Coulomb gases literature \cite{serfaty2023gaussian}, is continuous with respect to the norm of uniform convergence of the function and all its derivatives. Establishing this in the high-temperature regime is non-trivial, due to the non-linearity of \eqref{eqmeasure:charac} in $V$. Theorem \ref{thm:conteqdensity} provides a partial answer,  in the case of the one-dimensional log gas, to a question raised in \cite{serfaty2023gaussian} concerning the perturbed thermal equilibrium measure. To show this result, our method is based on an application of the Banach fixed-point theorem to \eqref{eqmeasure:exp}.

\subsection{Notations and conventions}\label{section:notations}
\begin{itemize}
	\item Let $X$ be an open set of $\R^p$, we denote by $\mcal{C}^k(X)$ (resp. $L^p(X)$) the space of functions differentiable $k$-times for which the $k$-th derivative is continuous (resp. $p^{\mathrm{th}}$-power integrable functions) on $X$. $\mcal{C}_c^k(X)$ denotes the space of functions of class $k$ on $X$ with compact support. For $p\in\llbracket1,+\infty\rrbracket$, we denote by $L^p(X)$ the usual Lebesgue spaces on $X$ and by $L^p(\mu)$ the Lebesgue spaces with respect to a Borel measure $\mu$ on $\R$. Furthermore, we define $L^2_0(\mu)$ by $\left \{u\in L^2(\mu),\,\int_\R u(x)\diff\mu(x)=0\right \}$. For a function of several variables $f$, we denote the derivative operator with respect to its $i$-th variable by $\partial_i f$.
	\item The space of functions $f$ such that $f^{(k)}\in L^\infty(\R)$ for all $k=0,\dots,n$ is denoted $W_n^\infty(\R)$. Its norm is classically $\|f\|_{W_n^\infty(\R)}\defi \underset{k\in\llbracket0,n\rrbracket}{\max}\|f^{(k)}\|_{L^\infty(\R)}$.
	\item Let $f\in L^2(\R)$, we denote by $\mcal{H}[f]$ the Hilbert transform of $f$ defined by
	$$	\mcal{H}[f](x)\defi \fint_\R\dfrac{f(y)}{y-x}\diff y$$
	where $\displaystyle \fint$ stands for the Cauchy principal value integral.
	\item We denote the Fourier transform of $f\in L^1(\R)\cap L^2(\R)$ by $$\mcal{F}[f](t)\defi \displaystyle\int_\R f(x)e^{-\ii tx}\diff x.$$
	When $\mu$ is a signed measure over $\R$, we shall denote its Fourier transform by the same symbol $\mcal{F}[\mu]$.
	\item The $1/2$-norm is defined for any function $f$ which makes this quantity finite
	$$\|f\|_{1/2}^2\defi  \displaystyle\int_\R |t|\left|\mcal{F}[f](t)\right|^2d t.$$
	\item We denote by $\mcal{M}_1(\R)$ the set of probability measures over $\R$. For $\mu,\mu'\in\mcal{M}_1(\R)$ we define the distance (possibly infinite) $D$ by 
	\begin{equation}
		\label{def:distanceMesures}
		D(\mu,\mu')\defi \left( \int_0^{+\infty} \frac{1}{t}\big|\mcal{F}[\mu-\mu'](t) \big|^2\diff t \right)^{1/2}.
	\end{equation}
		\item We define the Sobolev spaces for all $m\geq0$ by
	
	$$H^m(\R^n)\defi \left\{u\in L^2(\R^n), \|u\|_{H^{m}(\R^n)}
	<+\infty\right\}$$
	where $$\|u\|_{H^{m}(\R^n)}^2\defi \int_{\R^n}\left(1+\|t\|_{2}\right)^{2m}\left|\mcal{F}[u](t_1,\dots,t_n)\right|^2d^n\underline{t}.$$
	Above, $\|.\|_2$ denotes the Euclidean norm on $\R^n$.
	If $\mu\in\mcal{M}_1(\R)$, we also define $$H^k(\mu)\defi \left\{u\in L^2(\mu),\,u^{(k)}\in L^2(\mu)\right\}.$$
	The infinite order Sobolev space is given by $H^{\infty}(\R^p)\defi\cap_{k\geq0}H^k(\R^p)$.

\end{itemize}

\textbf{Outline of the paper.} In Section \ref{section:boundapriori}, we establish an \textit{a priori} bound on the $n$-linear statistics that is crucial in order to analyze the loop equations. To prove this bound, we first prove a concentration inequality for the empirical measure following \cite{maida2014free}. In Section \ref{section:loopequations}, we establish controls on the operators that appear as building blocks of the loop equations. In Section \ref{section4}, we prove controls on the so-called master operator. These will play a crucial role in the analysis of the loop equations. We then state the loop equations and establish the large-$N$ asymptotic expansion of the linear statistics, Theorem \ref{thm:correlators}, in Section \ref{section:asymptot correlators}. In Section \ref{section7}, we establish the continuity of the equilibrium density in the interpolation parameter, Theorem \ref{thm:conteqdensity}. Section \ref{section:freeenergy} is dedicated to the expansion of the partition function and an explicit form for the free energy associated with the Gaussian potential, as well as the interpolation formula. We conclude with Theorem \ref{thm:partitionfunction}. We detail in Appendix \ref{app A} some results obtained in \cite{DwoMemin} upon which this article largely relies. In Appendix \ref{appB}, we prove the continuity and the integrability of the constants that appear in our problem.

\subsection*{Acknowledgements}
The author wishes to thank Alice Guionnet, Gaultier Lambert and Trinh Khanh Duy for interesting discussions. I also thank Karol Kozlowski for his valuable advice and his idea for showing the continuity of the perturbed equilibrium density with respect to the interpolation. The author is also grateful to two anonymous referees that helped improving the quality and readability of the current paper.

\section{A priori bound on the linear statistics}\label{section:boundapriori}
As explained in the introduction, before analyzing the loop equations, one needs a bound that quantifies how small is a function integrated $n$ times against the recentered empirical measure $\Delta\mu_N\defi \mu_N-\mu_V$. Before addressing this, let us recall certain properties enjoyed by $\mu_V$ and the concentration results established in \cite{DwoMemin}.
\subsection{Equilibrium measure}
We recall the definition of the logarithmic potential (or sometimes called Symm's operator) $U^f$ of a function $f:\R \to \R$. When it is defined, the latter is given for all $x\in\R$ by
\begin{equation}
	\label{def:potentielLogarithmique}
	U^{f}(x)\defi-\int_\R\log|x-y|f(y)\diff y\,.
\end{equation}
One can check that
$\big(U^{f}\big)'=\mathcal{H}[f]$.

We now describe the regularity of the equilibrium density $\rho_V$ characterized by \eqref{eqmeasure:charac}.

\begin{lemma}\cite[Lemma 2.2]{DwoMemin}
	\label{lem:regularitedensite}
	\begin{itemize}
		\item The support of $\mu_V$ is $\R$ and there exists a $P$-dependent constant $C_V$ such that for all $x\in \R$, 
		$$ \rho_V(x) \leq C_V (1+|x|)^{2P}e^{-V(x)}\,.$$
		\item The density $\rho_V\in\mcal{C}^\infty(\R)$ and it holds
		\begin{equation}\label{deriv1}
			\rho_V\,'=-\Big(V'+2 P\mcal{H}[\rho_V]\Big)\rho_V,
		\end{equation}
		as well as
		\begin{equation}\label{deriv2}
			\rho_V\,''=\Big(-2P\mcal{H}[\rho_V]'-V''+V'^2+4P^2\mcal{H}[\rho_V]^2+4 PV'\mcal{H}[\rho_V]\Big)\rho_V\,.
		\end{equation}
	\end{itemize}
\end{lemma}

\subsection{Concentration inequality}
We now use an idea introduced by \cite{maida2014free} and based on a comparison between a configuration $\underline{x} = (x_1,\ldots,x_N)$ sampled with $\mathbb{P}_N^{V}$ and a regularized version $\underline{y}=(y_1,\ldots,y_N)$, which we describe here.
\begin{definition}
	\label{def:regularized}
	Let $\underline{x} = (x_1,\ldots,x_N)\in \R^N$ and suppose (up to reordering) that $x_1\leq x_2\leq \ldots \leq x_N$. We define $\underline{y}\in \R^N$ by:
	$$y_1 \defi  x_1\hspace{1cm} \mrm{\text{ and }}\hspace{1cm}  \forall 0\leq k \leq N-1,\;\;\; y_{k+1}\defi y_k + \max\left \{x_{k+1}-x_k, e^{-(\log N)^2}\right \}.$$
	We denote $\mu_N^{(y)}\defi \dfrac{1}{N}\displaystyle\sum_{i=1}^{n}\delta_{y_i}$ and also define $\mu_{N,u}^{(y)}\defi \mu_N^{(y)}\ast\mcal{U}_N$ the convolution between $\mu_N^{(y)}$ and $\mcal{U}_N$ the uniform measure on $\left[0,N^{-2}e^{-(\log N)^2}\right]$. 
\end{definition}
Note that the configuration $\underline{y}$ given by the previous definition satisfies $y_{k+1}-y_k\geq e^{-(\log N)^2}$, and $\underline{y}$ is close to $\underline{x}$ in the sense that 
\begin{equation}
	\label{ineq:boundedDifferences}
	\sum_{k=1}^N |x_k-y_k| \leq N^2e^{-(\log N)^2}\,.
\end{equation}
One can note that we have $|x_k-y_k| = y_k-x_k \leq (k-1)e^{-(\log N)^2}$, and we get \eqref{ineq:boundedDifferences} by summing these inequalities. As in the proof of \cite[Theorem 1.5]{DwoMemin}, we obtain a bound on the density:

\begin{theorem}\label{thm:bounddensity}
	For all $N\geq 1$ and $\underline{x}=\big(x_1,\dots,x_N\big)\in\R^N$,
	\begin{equation}\label{ineq:density}
		p_N^{V}\big(\underline{x}\big) \leq \exp\Big(-NPD^2\left (\mu_{N,u}^{(y)},\mu_V\right )+K_V+2P(\log N)^2\Big)\cdot\prod_{i=1}^N\rho_V(x_i)
	\end{equation}
	where $K_V\defi 2P\|\mcal{H}[\rho_V]\|_{\infty}+C+P\displaystyle\Big|\iint_{\R^2}\log|x-y|\diff\mu_V(x)\diff\mu_V(y)\Big|$ for some fixed, $V$-independent constant $C$ and with $D$ as given in \eqref{def:distanceMesures}.
\end{theorem}
Note that we have to keep the dependence on $V$ in all of the constants involved in our problem.
\subsection{A priori bound on linear statistics}
Thanks to the bound given in Theorem \ref{thm:bounddensity}, we can prove the below \textit{a priori} bound on the linear statistics. This bound is \textit{a priori }in the sense that it is not optimal, namely, we will show later that for the $n$-linear statistics are $O(N^{-\lceil n/2\rceil})$ versus $O(N^{- n(1-\varepsilon)/2})$ as predicted by the \textit{a priori} bound. Anyway, this will allow us to neglect, in the loop equations, the terms that are integrated with respect to $\bigotimes^{n}\Delta\mu_N$, where $\Delta\mu_N\defi \mu_N-\mu_{V}$. The following theorem and its proof are just adaptations of \cite[Corollary 3.1.10]{borot2016asymptotic}.

\begin{proposition}[A priori bound on linear statistics]\label{a priori bound}Let $\varepsilon>0$, there exits $C>0$ (depending only on $n$ and $\varepsilon$) such that for all $f$ in $W^\infty_1(\R^n)\cap H^{n/2}(\R^n)$, it holds
	$$\boxed{\left| \Braket{f}_{\bigotimes^n\Delta\mu_N}\right|\leq\dfrac{Ce^{K_V}}{N^{\frac{n}{2}(1-\varepsilon)}}\left(\left \|f\right\|_{W_1^\infty(\R^n)}+\|f\|_{H^{n/2}(\R^n)}\right).}$$
	where $K_V$ is defined in Theorem \ref{thm:bounddensity}.
\end{proposition}
\begin{proof}
	We use the decomposition $\Delta\mu_N=\Big(\mu_N-\mu^{(y)}_{N,u}\Big)+\Delta\mu_{N,u}^{(y)}$ where  $\Delta\mu_{N,u}^{(y)}=\mu^{(y)}_{N,u}-\mu_V$ and obtain:
	\begin{multline*}
		\Braket{f}_{\otimes^n\Delta\mu_N}=\sum_{l=0}^{n-1}\sum_{i_1<\dots<i_l}^{n}\mbb{E}_N^{V}\left[\int_{\R^n}f(\xi_1,\dots,\xi_n)\prod_{j=1}^{l}\diff\Delta\mu_{N,u}^{(y)}(\xi_{i_j})\prod_{\substack{j=1\\\neq i_1,\dots,i_l}}^{n}\diff\left(\mu_N-\mu_{N,u}^{(y)}\right) (\xi_j)\right] \\+\Braket{f}_{\otimes^n\Delta\mu_{N,u}^{(y)}}.
	\end{multline*}
	Since the $x_i$'s are not far from the $y_i$'s, we have the following bound by the mean value theorem and the fact that all the involved measures are probability measures:
	\begin{equation*}
		\mbb{E}_N^{V}\left[\int_{\R^n}f(\xi_1,\dots,\xi_n)\prod_{j=1}^{l}\diff\Delta\mu_{N,u}^{(y)}(\xi_{i_j})\prod_{\substack{j=1\\\neq i_1,\dots,i_l}}^{n}\diff\left(\mu_N-\mu_{N,u}^{(y)}\right) (\xi_j)\right]
		\leq \dfrac{C_n\left\|f\right\|_{W_1^\infty(\R^n)}N}{e^{(\log N)^2}}
	\end{equation*}
	for some constant $C_n>0$ only depending on $n$.
	
	Let's focus now on $\Braket{f}_{\otimes^n\Delta\mu_{N,u}^{(y)}}$. We know by Theorem \ref{thm:bounddensity} that $$\mbb{P}_N^{V}(\Omega_N)=e^{K_V}O\left(e^{-cN^{\varepsilon}}\right)\hspace{1cm}\text{ where }\hspace{1cm}\Omega_N\defi \left\{\lambda\in\R^N,\,D^2\left [\mu_{N,u}^{(y)},\mu_V\right ]>\dfrac{1}{N^{1-\varepsilon}}\right\}$$ for some $c>0$ independent of $V$ and for a remainder controlled $V$-independently. 
	It ensures that:
	$$\left| \Braket{f}_{\otimes^n\Delta\mu_{N,u}^{(y)}}\right|\leq C e^{K_V} e^{-cN^{\varepsilon}}\left\|f\right\|_{W_0^\infty(\R^n)}+\mfrak{R}_{N}[f]$$
	where $$\mfrak{R}_{N}[f]\defi \mbb{E}_N^{V}\left[\mbf{1}_{\Omega_N^c}\int_{\R^n}f(\xi_1,\dots,\xi_n)\diff\Delta\mu_{N,u}^{(y)}\,^{\otimes^n}(\xi_1,\dots,\xi_n)\right].$$
	By Plancherel formula and Cauchy-Schwarz inequality, one gets:
	\begin{align*}
		\mfrak{R}_{N}[f]&=\mbb{E}_N^{V}\left[\mbf{1}_{\Omega_N^c}\int_{\R^n}\mcal{F}\left [f\right ](\varphi_1,\dots,\varphi_n)\cdot\prod_{j=1}^{n}\mcal{F}\left[\Delta\mu_{N,u}^{(y)}\right](-\varphi_j)\cdot\dfrac{\diff^n\underline{\varphi}}{(2\pi)^n}\right]
		\\&\leq\left(\int_{\R^n}\left| \mcal{F}\left [f\right ](\varphi_1,\dots,\varphi_n)\right|^2\cdot \prod_{j=1}^{n}|\varphi_j|\cdot\dfrac{\diff^n\underline{\varphi}}{(2\pi)^n}\right) ^{1/2}.\mbb{E}_N^{V}\left[\mbf{1}_{\Omega_N^c}2^{\frac{n}{2}}D^n\left[\mu_{N,u}^{(y)},\mu_V\right]\right]\nonumber
		\\&\leq\dfrac{2^{\frac{n}{2}}}{N^{\frac{n}{2}(1-\varepsilon)}}\left(\int_{\R^n}\left| \mcal{F}\left [f\right ](\varphi_1,\dots,\varphi_n)\right|^2\cdot \left\{1+\left( \sum_{j=1}^{n}|\varphi_j|^2\right)^{1/2}\right\}^n\cdot \dfrac{\diff^n\underline{\varphi}}{(2\pi)^n}\right) ^{1/2}\nonumber
		\\&\leq2^{\frac{n}{2}}\dfrac{\|f\|_{H^{n/2}(\R^n)} }{N^{\frac{n}{2}(1-\varepsilon)}}\nonumber
	\end{align*}
	which concludes the proof.
\end{proof}
\section{Controls of  $\mcal{D}$}\label{section:loopequations}
In this section, we will set some definitions of operators which arise as building blocks of the loop equations. After defining them, we will prove their continuity on appropriate spaces. This will ultimately allow us to apply the \textit{a priori} bounds given in Proposition \ref{a priori bound}.
\subsection{Definitions}
The operators that will appear in the loop equations at level $n\geq 2$ are constructed via the following extension procedure, allowing one to extend operators acting on univariate functions into operators acting on multivariate functions.

\begin{definition}[Extension of operators]
	Given an operator $\mcal{B}$ that acts on functions of one variable and yields a function of $l\in\llbracket1,2\rrbracket$ variables, $\phi_n$ a function of $n$ variables, we define $\mcal{B}_1$ by:
	\begin{equation}\label{def:extsub1}
		\mcal{B}_1[\phi_n](\xi_1,\dots,\xi_{n+l-1})\defi\mcal{B}\left [\phi_n(\cdot,\xi_{l+1},\dots,\xi_{n+l-1})\right ](\xi_1,\dots,\xi_l).
	\end{equation}
\end{definition}
For example $\mcal{D}_1[\phi_n](x_1,\dots,x_{n+1})=\mcal{D}[\phi_n(\cdot,x_3,\dots,x_{n+1})](x_1,x_2).$
\subsection{Control on the non-commutative derivative operator}
A first example of an operator appearing in the loop equations is the non-commutative derivative (NCD) operator.

\begin{definition}
	Let $f\in\mcal{C}^1(\R)$, we define the NCD operator $\mcal{D}[f]$ by:
	$$\forall x,y\in\R,\,\mcal{D}[f](x,y)\defi\begin{cases}
		\dfrac{f(x)-f(y)}{x-y}&\text{ if }x\neq y
		\\f'(x) &\text{ if }x=y
	\end{cases}.$$
\end{definition}
In the following, we fix $p\geq2$ and show continuity results for the operator $\mcal{D}_1$ in the $H^{n}$ and $W_n^{\infty}$-norms.
\begin{theorem}[$H^{n}$-control for $\mcal{D}_1$]\label{thm:fini diff ope reg}Let $n\geq1$, there exists $C>0$ (depending only on $n$) such that for all $f\in\mcal{C}^n(\R^{p-1})\cap H^{n+1}(\R^{p-1})$,
	$$\boxed{\left \|\mcal{D}_1[f]\right \|_{H^{n}(\R^p)}\leq C\|f\|_{H^{n+1}(\R^{p-1})}.}$$
\end{theorem}

Before showing this inequality, we need to show a general form of the derivatives of $\mcal{D}_1[f]$.
\begin{lemma}[General form for derivatives of \texorpdfstring{$\mcal{D}_1[f]$}{TEXT}]\label{lem: opD form leibniz} Let $\underline{m}=(m_1,\dots,m_p)\in\N^p$ satisfy $m_1\geq m_2$ and $\sum_{i=1}^{p}m_i\leq n$. Let $x_1,\dots,x_p\in\R$ be such that $x_1\neq x_2$, then one has:
	\begin{equation}\label{eq:operfindiff}
		\partial^{\underline{m}}\mcal{D}_1[f](x_1,\dots,x_n)=\sum_{j=0}^{m_2}C_{m_1,m_2,j}\dfrac{\left(g^{(j)}(x_2)-\displaystyle\sum_{k=0}^{m_1-j}\dfrac{g^{(k+j)}(x_1)}{k!}(x_2-x_1)^{k}\right)}{(x_2-x_1)^{m_1+m_2+1-j}}
	\end{equation}
	with $C_{m_1,m_2,j}\defi \binom{m_2}{j}(m_1+m_2-j)!(-1)^{m_2-j}$ and $g=\partial_{2}^{m_3}\dots\partial_{p-1}^{m_p}f(.,x_3,\dots,x_p)$.
\end{lemma}

\begin{proof}
	First, it is easy to verify that $\mcal{D}_1[f]\in\mcal{C}^{n}(\R^p)$ for $x_1\neq x_2$. Secondly, when $n\geq2$, by the Schwarz theorem, the order of the partial derivatives does not matter. It is only the derivatives with respect to $x_1$ and $x_2$ that are non-trivial to compute. Indeed, let $\underline{x}=(x_1,\dots,x_p)\in\R^p$ be such that $x_1\neq x_2$, then
	$$\partial_{3}^{m_3}\dots\partial_{p}^{m_p}\mcal{D}_1[f](x_1,\dots,x_p)=\dfrac{g(x_1)-g(x_2)}{x_1-x_2}$$
	with $g\defi \partial_{2}^{m_3}\dots\partial_{p}^{m_p}f(.,x_3,\dots,x_p)$. By applying the Leibniz formula when differentiating $m_1$ times with respect to $x_1$, one gets:
	$$\partial_{x_1}^{m_1}\partial_{x_3}^{m_3}\dots\partial_{x_p}^{m_p}\mcal{D}_1[f](x_1,\dots,x_p)=\dfrac{m_1!}{(x_2-x_1)^{m_1+1}}\left(g(x_2)-\sum_{k=0}^{m_1}\dfrac{g^{(k)}(x_1)}{k!}(x_2-x_1)^{k}\right).$$
	Again, we differentiate $m_2$ times with respect to $x_2$ and apply the Leibniz formula to get \eqref{eq:operfindiff}.
\end{proof}

\begin{proof}\textit{[of Theorem \ref{thm:fini diff ope reg}]}
	Let $\underline{m}\defi (m_1,\dots,m_p)\in\N^p$ be such that $m\defi \sum_{i=1}^{p}m_i\leq  n$. Without loss of generality, we can assume that $m_1\geq m_2$. Let's show that $$\left \|\partial_{1}^{m_1}\partial_{2}^{m_2}\dots\partial_{p}^{m_p}\mcal{D}_1[f]\right\|_{L^2(\R^p)}\leq C\|f\|_{H ^{m+1}(\R^{p-1})}$$ with $C>0$ independent of $f$. The idea is to prove separately the $L^2$ control on $\partial^{\underline{m}}\mcal{D}_1[f]$ close to the singularity (the diagonal) and far from it. To do so, we will use the Taylor formula with integral remainder to deal with the singularity and Lemma \ref{lem: opD form leibniz} when we are at a fixed distance from the diagonal.

	\textbf{Close to the diagonal:} We first show this inequality on the subspace $\left\{\underline{x}\in\R^n,|x_1-x_2|\leq1\right\}$. First note that $$\mcal{D}_1[f](x_1,x_2,\dots,x_p)=\int_{0}^{1}\partial_1f\left(x_1+t(x_2-x_1),x_3,\dots,x_p\right)\diff t $$
	an so by differentiating under the integral sign and by Jensen's inequality, we obtain:
	\begin{multline*}
		\left |\partial_{1}^{m_1}\partial_{2}^{m_2}\dots\partial_{p}^{m_p}\mcal{D}_1[f](x_1,x_2,\dots,x_p)\right |^2
		\\\leq\int_{0}^{1}(1-t)^{2m_1}t^{2m_2}\partial_{1}^{m_1+m_2+1}\partial_{2}^{m_3}\dots\partial_{p-1}^{m_p}f\left(x_1+t(x_2-x_1),x_3,\dots,x_p\right)^2\diff t.
	\end{multline*}
	Hence, by integrating with respect to $\underline{x}$, changing $x_2-x_1$ into $\widetilde{x_2}$, and using Fubini, we get:
	\begin{multline*}
		\int_{\R^n}\left |\partial_{1}^{m_1}\partial_{2}^{m_2}\dots\partial_{p}^{m_p}\mcal{D}_1[f](x_1,x_2,\dots,x_p)\right |^2\mbf{1}_{|x_2-x_1|<1}\diff^nx
		\\\leq\int_{0}^{1}\diff t(1-t)^{2m_1}t^{2m_2}\int_{-1}^{1}\diff\widetilde{x_2}\int_\R \diff x_1\int_\R \diff x_3\dots\int_\R \diff x_p\partial_{1}^{m_1+m_2+1}\partial_{2}^{m_3}\dots\partial_{p-1}^{m_p}f(x_1+t\widetilde{x_2},x_3,\dots,x_p)^2
		\\\leq C\|\partial_{1}^{m_1+m_2+1}\partial_{2}^{m_3}\dots\partial_{p-1}^{m_p}f\|_{L^2(\R^{p-1})}^2\leq C\|f\|_{H^{m+1}(\R^{p-1})}^2.
	\end{multline*}
	\textbf{Far from the diagonal:} Now we deal with the subset $\left\{\underline{x}\in\R^n,|x_1-x_2|\geq1\right\}$. By Jensen's inequality and Lemma \ref{lem: opD form leibniz}, we get:
	\begin{multline*}
		\int_\R \diff x_1\int_{|x_2-x_1|>1}\diff x_2\left|\partial_{1}^{m_1}\partial_{2}^{m_2}\dots\partial_{p}^{m_p}\mcal{D}_1[f](x_1,x_2,\dots,x_p)\right|^2\leq(m_2+1)\sum_{j=0}^{m_2}C_{m_1,m_2,j}^2(m_1-j+1)
		\\\times\int_\R \diff x_1\int_{|x_2-x_1|>1}\diff x_2\dfrac{g^{(j)}(x_2)^2+\displaystyle\sum_{k=0}^{m_1-j}\dfrac{g^{(k+j)}(x_1)^2}{k!^2}(x_2-x_1)^{2k}}{|x_2-x_1|^{2m_1+2m_2+2-2j}}.
	\end{multline*}
	with $g=\partial_{2}^{m_3}\dots\partial_{p-1}^{m_p}f(.,x_3,\dots,x_p)$. For all $j\in\llbracket0,m_2\rrbracket$, the double integral in the last line can be estimated with another constant $C$ depending only on $m_1$ and $m_2$. For that, we use Fubini's theorem, the fact that $m_1\geq m_2$: and obtain
	\begin{equation*}
		\int_\R \diff x_1\int_{|x_2-x_1|>1}\diff x_2\dfrac{g^{(j)}(x_2)^2+\displaystyle\sum_{k=0}^{m_1-j}\dfrac{g^{(k+j)}(x_1)^2}{k!^2}(x_2-x_1)^{2k}}{|x_2-x_1|^{2m_1+2m_2+2-2j}}\leq \big\|g\big\|_{H^{m_1}(\R)}^2.
	\end{equation*}
	Hence, after suming over $j$ and changing the constant appropriately, we integrate over $x_3,\dots,x_p$ to obtain:
	\begin{multline*}
		\int_{\R^{p-2}}\diff x_3\dots \diff x_p\int_\R \diff x_1\int_{|x_2-x_1|>1}\diff x_2\left|\partial_{1}^{m_1}\partial_{2}^{m_2}\dots\partial_{p}^{m_p}\mcal{D}_1[f](x_1,x_2,\dots,x_p)\right|^2
		\\\leq C\max_{l\in\llbracket 1,m_1\rrbracket}\left\|\partial_{1}^{l}\dots\partial_{p-1}^{m_p}f\right \|_{L^2(\R^{p-1})}^2
		\leq C\left\|f\right \|_{H^{m+1}(\R^{p-1})}^2.
	\end{multline*}
	This is enough to conclude.
\end{proof}
Since in Proposition \ref{a priori bound}, the bound on the linear statistic involves the $W^\infty_1(\R^p )$-norm, we state the following result.

\begin{proposition}\label{prop:diff fini linfini}
	There exists a $C>0$ such that for all $f\in\mcal{C}^n(\R^{p-1})\cap W^\infty_n(\R^{p-1})$,
	$$\left \|\mcal{D}_1[f]\right \|_{W^\infty_n(\R^p)}\leq C(n)\|f\|_{W^\infty_{n+1}(\R^{p-1})}.$$
\end{proposition}

\begin{proof}
	This follows from Lemma \ref{lem: opD form leibniz} together with the Taylor formula with integral remainder.
\end{proof}

\section{Control on the master operator $\Xi$}\label{section4}
In this section, we study the so-called master operator which plays an essential role in the following. Indeed, proving continuity estimates for its inverse (which is well defined, see the following subsection) is a crucial step if one wants to analyze the loop equations. We first give all the definitions of the different operators involved in the study of $\Xi^{-1}$ in Theorem \ref{thm inverse master}, then we show some appropriate decomposition for the derivative of $\Xi^{-1}[f]$ in Theorem \ref{thm: stru inv Xi} and from that we finally conclude the continuity results needed to study the loop equations, namely the control in $H^{n}$-norms in Theorems \ref{thm: cont inverse master} and \ref{thm:continversemaster thetaxi L2} and in $W_n^{\infty}$-norms in Theorem \ref{thm: cont inverse master linfini} and Corollary \ref{thm: cont theta inverse master linfini}.
\subsection{Definition}
We recall the definition of the operator $\mcal{L}$.
\begin{definition}\label{def:oper Aet W}
	We define, for a sufficiently smooth function $f$, the operator $$\mcal{L}[f]\defi \Xi\big [f'\big]=-\mcal{A}[f]-2P\mcal{W}[f]$$
	where $$\mcal{A}[f]\defi -\dfrac{\left(f'\rho_V\right) '}{\rho_V} \hspace{2cm}\text{ and } \hspace{2cm}\mcal{W}[f]\defi -\mcal{H}\left [f'\rho_V\right ]+\int_\R\mcal{H}\left [f'\rho_V\right ](y)\diff\mu_V(y).$$
\end{definition}$\mcal{L}$ is an unbounded operator on the space $\msf{H}$ defined by:
\begin{equation}\label{eq:def H}
	\mathsf{H}\defi \left\{u\in L^2\left (\mu_V\right )\ \Big|\ u'\in L^2\left (\mu_V\right ), \int_\R u(x)\diff\mu_V(x)=0\right\},\qquad \Braket{u,v}_\mathsf{H}\defi \Braket{u',v'}_{L^2\left (\mu_V\right )}.
\end{equation}
This space is indeed a Hilbert-space by the fact that $\mu_{V}$ verifies the Poincaré inequality (see assumption\textit{ \ref{assumption3}}). Its domain is defined by $\mcal{D}(\mcal{L})=\mcal{D}(\mcal{A})\defi \left\{u\in\msf{H},\,\mcal{A}[u]\in\msf{H}\right\}$ by \cite[Proposition 6.8]{DwoMemin}. $\mcal{A}$ and $\mcal{W}$ are also unbounded operators on $\msf{H}$ defined respectively on $\mcal{D}(\mcal{A})$ and $\mcal{D}(\mcal{W})\defi \left\{u\in\msf{H},\,\mcal{W}[u]\in\msf{H}\right\}$, see \cite[Section 6]{DwoMemin} for a more detailed description of these operators. We also emphasize that $\mcal{D}(\mcal{A})\subset\mcal{D}(\mcal{W})$ by \cite[Remark 6.4]{DwoMemin}.

For the next theorem, we recall that $\mcal{A}:\mcal{D}(\mcal{A})\rightarrow \msf{H}$ is an invertible and diagonalizable operator with positive countable spectrum, see \cite[Proposition 6.3]{DwoMemin}. We denote by $\lambda_1(\mcal{A})>0$ its smallest eigenvalue. This quantity has a role in our problem since for all $f\in\msf{H}$, $\|\mcal{L}^{-1}[f]\|_{\msf{H}}\leq\lambda_1(\mcal{A})^{-1}\|f\|_{\msf{H}}$ see \cite[Theorem 6.9]{DwoMemin}.

\begin{theorem}[Inversion of the master operator]\label{thm inverse master}$\Xi:\mfrak{D}(\Xi)\longrightarrow\msf{H}$ is invertible, of inverse defined for all $g\in\msf{H}$ by:
	$$\Xi^{-1}[g]\defi \big(\mcal{L}^{-1}[g]\big)'$$
	where $\mfrak{D}(\Xi)\defi \left\{f\in \mcal{C}^1(\R),\, \exists v\in\mfrak{D}(\mcal{L}),\,f=v'\right\}$. Furthermore for all $f\in\msf{H}$, 
	\begin{equation}\label{ineq:contXi}
		\left\|\Xi^{-1}[f]\right \|_{L^2\left(\mu_V\right)}\leq C_{\mcal{L}}\|f'\|_{L^2(\mu_V)}
	\end{equation}where $C_\mcal{L}\defi \lambda_1(\mcal{A})^{-1}$.
\end{theorem}
\begin{proof}
	To prove that $\Xi$ is invertible on $\mfrak{D}(\Xi)$, the only thing to prove is that for all $v\in\mcal{D}(\mcal{L})$, $v'\in \mcal{C}^1(\R)$ which is true because by definition, for all $v\in\mfrak{D}(\mcal{A})$, $(v'\rho_V)'=\rho_V\mcal{A}[v]\in\rho_V\msf{H}\subset\mcal{C}^{0}(\R)$. This means that $v'\rho_V\in\mcal{C}^1(\R)$ and thus that $v\in\mcal{C}^{2}(\R)$. Now it is obvious from the definition of $\mcal{L}$ and $\Xi^{-1}$ that $\Xi\circ\Xi^{-1}=\msf{id}_\msf{H}$ on $\mfrak{D}(\Xi)$. The estimate comes from the fact that given $f\in\msf{H}$, one has $\Xi^{-1}[f]=\left(\mcal{L}^{-1}[f]\right)'$. Then
	\begin{equation*}
		\left\|\Xi^{-1}[f]\right \|_{L^2\left(\mu_V\right)}=\left\|\mcal{L}^{-1}[f]\right \|_{\msf{H}}\leq C_\mcal{L}\|f\|_{\msf{H}}.
\end{equation*}\end{proof}

The crucial step when one wants to deduce Theorem \ref{thm:correlators} from the analysis of the loop equations, is to obtain controls on the master operator which we will show in this section. These bounds will allow us to apply the bound obtained in Proposition \ref{a priori bound} to functions like $\Xi^{-1}[\phi]$.
\subsection{Preliminaries}
We define an operator $\mcal{O}$ whose iterations will appear in the derivatives of the inverse of the master operator (which exists because of Lemma \ref{lem:reginverse}).

\begin{definition}\label{def:operateurO}
	Let $\mcal{O}$ be the operator defined on smooth enough functions by:
	\begin{equation}\label{eq:defopeO}
		\mcal{O}[f](x)\defi \left(\dfrac{f\rho_V}{\rho_V'}\right)'(x)
	\end{equation}
\end{definition}

\begin{remark} For example, with $\alpha\defi \dfrac{\rho_V}{\rho_V'}$, it holds that:
	\begin{itemize}
		\item $\mcal{O}[f]=\alpha' f+\alpha f'$.
		\item $\mcal{O}^2[f]=(\alpha\alpha')'f+3\alpha\alpha'f'+\alpha^2f''$.
		\item $\mcal{O}^3[f]=\left(\alpha\left(\alpha\alpha'\right)' \right)'f+\left(4\alpha^2\alpha''+7\alpha{\alpha'}^{2}\right)f'+\left(6\alpha'\alpha^2\right)f'' +\alpha^3f^{(3)}$.
	\end{itemize}
\end{remark}

In order to give a more precise description of $\mcal{O}^k$, which will allow us to analyse its asymptotics at infinity, we need the following definition.

\begin{definition}[Differential degree]
	Let $f$ be a function of one variable defined on $\R$ differentiable $n$ times, we define the differential degree denoted by $d_\partial^f$ with respect to $f$ by
	$$d_\partial^f\left(\displaystyle\prod_{i=0}^{n}\left(f^{(i)}\right)^{\alpha_i}\right)\defi \sum_{i=0}^{n}i\alpha_i$$
\end{definition}
For example the differential degree with respect to $f$ of $\left(f'\right)^2$ and $f''f$ is 2, while $d_\partial^f\left(\left(f^{(3)}  \right)^2\right)=6$. Using the notion of differential degree, we are now able to state the next theorem.

\begin{theorem}\label{thm: struc op O}
	Let $k\geq1$, $f\in\mcal{C}^k(\R)$, there exists a family of polynomials $(P_i^{k})_{0\leq i\leq k}$ such that
	\begin{equation}\label{eq:oper O recc}
		\boxed{\mcal{O}^{k}[f]=\sum_{i=0}^{k}f^{(k-i)}P_i^{k}(\alpha,\dots,\alpha^{(i)}),\hspace{1cm}\mrm{with} \hspace{1cm}\alpha\defi \dfrac{\rho_V}{\rho_V'}}
	\end{equation}
	In fact, $P_i^{k}\left (\alpha,\dots,\alpha^{(i)}\right )$, $i\in\llbracket0,k\rrbracket$, is the unique homogeneous polynomial in $i+1$ variables, with differential degree with respect to $\alpha$ equal to $i$, degree $k$ and with coefficients independent of $V$ satisfying the following reccurence relations:
	\begin{itemize}
		\item $P_0^{k+1}(\alpha)=\alpha P_0^{k}\left (\alpha\right )=\alpha^{k+1}$
		\item $\forall i\in\llbracket1,k\rrbracket$, $P_i^{k+1}\left(\alpha,\dots,\alpha^{(i)}\right)=\left(\alpha P_{i-1}^{k}(\alpha,\dots,\alpha^{(i-1)})\right)'+\alpha P_{i}^{k}\left (\alpha,\dots,\alpha^{(i)}\right ) $ 
		\item $P_{k+1}^{k+1}\left (\alpha,\dots,\alpha^{(k+1)}\right )=\left(\alpha P_{k}^{k}\left (\alpha,\dots,\alpha^{(k)}\right )\right)'=\left(\left(\alpha'\alpha\right)'\dots\alpha\right)'  $
	\end{itemize}
\end{theorem}

\begin{proof}
	Let's prove it by induction. For $k=1$,  $\mcal{O}[f]=\alpha' f+\alpha f'$ and so by setting $P_0^1(\alpha)=\alpha$, which is homogeneous, of degree $1$ and of differential degree $0$, and $P_1^1\left (\alpha,\alpha'\right )=\alpha'$ which is of degree $1$ and differential degree $1$, this proves the claim. Suppose that \eqref{eq:oper O recc} holds at rank $k\in\N^*$, then
	\begin{multline*}
		\mcal{O}^{k+1}[f]=\left(\alpha\mcal{O}^{k}[f]\right)'=\sum_{i=0}^{k}f^{(k-i)}\left[\alpha P_i^{k}(\alpha,\dots,\alpha^{(i)})\right]'+\sum_{i=0}^{k}f^{(k-i+1)}\alpha P_i^{k}\big(\alpha,\dots,\alpha^{(i)}\big)
		\\=\left[\alpha P_k^k\big (\alpha,\dots,\alpha^{(k)}\big )\right]'f+\sum_{i=0}^{k-1}f^{(k-i)}\left\{\left[\alpha P_i^k\big(\alpha,\dots,\alpha^{(i)}\big) \right]'+\alpha P_{i+1}^k \big(\alpha,\dots,\alpha^{(i+1)}\big)\right\}
		\\+\alpha P_0^k(\alpha)f^{(k+1)}
	\end{multline*}
	Hence by setting $P_0^{k+1}(\alpha)\defi \left(\alpha P_0^k(\alpha)\right)'$, $ P_{k+1}^{k+1}\left (\alpha,\dots,\alpha^{(k+1)}\right )\defi \left(\alpha P_{k}^{k}\left (\alpha,\dots,\alpha^{(k)}\right )\right)'$ and for all $i\in\llbracket0,k-1\rrbracket$, $P_{i+1}^{k+1}\left (\alpha,\dots,\alpha^{(j+1)}\right )\defi \left(\alpha P_{i}^{k}(\alpha,\dots,\alpha^{(i)})\right)'+\alpha P_{i}^{k}\left (\alpha,\dots,\alpha^{(i)}\right ) $, we obtain the desired form of \eqref{eq:oper O recc} and the recurrence relations. It remains to check that the homogeneity and degree conditions hold at rank $k+1$. This follows from the recurrence relations for the $P_i^k$'s.
\end{proof}

\subsection{Closed form for $\Xi^{-1}$}
Before showing a closed form for the derivatives of $\Xi_{1}^{-1}[f]$ and their  $L^2$ properties, we first prove that, if $f$ is sufficiently smooth, these derivatives indeed exist in a strong sense.

\begin{lemma}[Regularity of the inverse]\label{lem:reginverse}
	Let $f\in\msf{H}$ such that $f\rho_V\in H^n(\R)$ with $n\geq2$, then $\rho_V\Xi^{-1}[f] \in H^{n+1}(\R)$. Furthermore if $f\in\msf{H}\cap\mcal{C}^n(\R)$ is such that $f\rho_V\in H^n(\R)$, then $\Xi^{-1}[f]\in\mcal{C}^{n+1 }(\R)$.
\end{lemma}
Note that the last condition is verified whenever $f$ and its derivatives are continuous and grow slower than $e^{-V}$ at infinity. The proof uses the operators $\mcal{L}$ and $\mcal{A}$ introduced in Definition \ref{def:oper Aet W}.
\begin{proof}
	We recall that $\Xi^{-1}[f]=\left(\mcal{L}^{-1}[f]\right)'$. When $f\in\msf{H}$,  we know that $\rho_V\left(\mcal{L}^{-1}[f]\right)'\in H^2(\R)$. This is because $\mcal{L}^{-1}[f]\in\mcal{D}\left(\mcal{A}\right)$ and $\left(\rho_V\left(\mcal{L}^{-1}[f] \right)'\right)'=\rho_V\mcal{A}\left[\mcal{L}^{-1}[f]\right] \in H^1(\R)$.
	
	We now want to show that $\left(\rho_V\left(\mcal{L}^{-1}[f]\right)'\right)' \in H^{n}(\R)$, let's show first that $\mcal{H}\left[\rho_V\left(\mcal{L}^{-1}[f]\right)'\right]\in H^n(\R)$. First observe that
	
	$$\rho_V\mcal{A}\circ\mcal{L}^{-1}[f]=-\rho_Vf+2P\rho_V\mcal{H}\left[\rho_V\left(\mcal{L}^{-1}[f]\right)'\right]-2P\rho_V\int_\R\mcal{H}\left[\rho_V\left(\mcal{L}^{-1}[f]\right)'\right](y)\diff y.$$
	Hence, since $\rho_V\left(\mcal{L}^{-1}[f]\right)'\in H^2(\R)$, so is $\mcal{H}\left [\rho_V\left(\mcal{L}^{-1}[f]\right)'\right]$. The last term is a constant and clearly belongs to $H^n(\R)$ for all $n\in\N$ hence $\rho_V\mcal{A}\circ\mcal{L}^{-1}[f]=-\left(\rho_V\left(\mcal{L}^{-1}[f]\right)'\right)'\in H^2(\R)$ 
	and hence $\rho_V\left(\mcal{L}^{-1}[f]\right)'\in H^3(\R)$. By bootstraping, this shows that $\left(\mcal{L}^{-1}[f]\right)'=\Xi^{-1}[f]\in\dfrac{1}{\rho_V} H^{n+1}(\R)\subset\mcal{C}^n(\R)$ by Sobolev-Hölder embedding theorem and hence that $\mcal{H}\left[\rho_V\left(\mcal{L}^{-1}[f]\right)'\right]\in H^{n+1}(\R)$. Since
	\begin{equation}\label{eq:deriveesecondereg}
		\left(\mcal{L}^{-1}[f]\right)''=f-\dfrac{\rho_V'}{\rho_V}\left(\mcal{L}^{-1}[f]\right)'-2P\left(\mcal{H}\left[\rho_V\left(\mcal{L}^{-1}[f]\right)'\right]-\int_\R\mcal{H}\left[\rho_V\left(\mcal{L}^{-1}[f]\right)'\right](y)\diff\mu_V(y)\right) 
	\end{equation}
	and that $\dfrac{\rho_V'}{\rho_V}\in \mcal{C}^\infty(\R)$, we can then conclude that, under the assumption that $f\in\mcal{C}^n(\R)$,  $\left(\Xi^{-1}[f]\right)'= \mcal{L}^{-1}[f]\in\mcal{C}^{n}(\R)$, hence $\Xi^{-1}[f]\in\mcal{C}^{n+1}(\R)$.
\end{proof}

The following lemma will also be useful for the controls on $\Xi^{-1}$ as it is convenient to bound differently close to infinity and on a compact.

\begin{lemma}\label{lem:M_V}
	There exists $M_V>0$ such that for all $ |x|\geq M_V$, $\Big|\dfrac{\rho_{V}'}{\rho_V}(x)\Big|\geq 1$.
\end{lemma}
\begin{proof}
	From Lemma \ref{lemme:boundedHilbert}, $\mcal{H}[\rho_V]$ is bounded and by assumption\textit{ \ref{assumption2}}, $V'(x)$ goes to infinity at $\pm\infty$, the conclusion follows from the fact that $\dfrac{\rho_{V}'}{\rho_V}=-V'-2P\mcal{H}[\rho_V]$.
\end{proof}

We are now able to prove that a closed form  holds for the derivatives of $\Xi^{-1}$.
The idea is to use the resolvant formula which gives that for all $f\in\msf{H}$,
\begin{equation}\label{eq:resolvant}
	\mcal{L}^{-1}[f]=-\mcal{A}^{-1}\left [f+2P\mcal{W}\circ\mcal{L}^{-1}[f]\right ]
\end{equation}
and for all $x\in\R$,
\begin{equation}\label{eq:Ainverse'}
	\mcal{A}^{-1}[f](x)=\dfrac{1}{\rho_{V}(x)}\int_{x}^{\pm\infty}f(t)\rho_{V}(t)\diff t.
\end{equation}
It does not matter if one chooses $+\infty$ or $-\infty$ in \eqref{eq:Ainverse'} since $\int_\R f(t)\rho_V(t)\diff t=0$ but it will be convenient to make the choice $\mrm{sgn}(x)\infty$ for reasons that will appear further. We also mention that \eqref{eq:resolvant} is well-defined because $\mcal{D}(\mcal{L})=\mcal{D}\left(\mcal{A}\right)\subset\mcal{D}(\mcal{W})$.
Before establishing the continuity for $\Xi^{-1}$, we need to introduce an operator $\mcal{X}$ that takes a function in $\dfrac{1}{\rho_V}H^{n}(\R)\cap H^{1}(\mu_V)$ and produces one belonging to $\dfrac{1}{\rho_V}H^{n}(\R)\cap\msf{H}$ by means of a recentering, where:
\begin{equation}\label{def:h1mu}
	H^{1}(\mu_V)\defi\left\{u\in L^{2}(\mu_V),\quad u'\in L^{2}(\mu_V)\right\} 
\end{equation}

\begin{definition}\label{def:recentring}
	Let $\phi\in\dfrac{1}{\rho_V}H^n(\R)\cap H^{1}(\mu_V)$, we define the operator $\mcal{X}$ by
	$$\mcal{X}[\phi](\xi)\defi\phi(\xi)-\int_{\R}\phi(t)\diff\mu_V(t).$$
\end{definition}

It is obvious to see that for all $\phi\in\dfrac{1}{\rho_V}H^n(\R)\cap H^{1}(\mu_V)$, $\mcal{X}[\phi]\in\msf{H}\cap\dfrac{1}{\rho_V}H^n(\R)$. We denote  for all $u$ such that  $\mcal{X}[u]\in\mfrak{D}(\Xi)$ and all $v$ such that $\mcal{X}[v]\in\msf{H}$, $$\widetilde{\Xi}[u]\defi \Xi\circ \mcal{X}[u],\hspace{2cm}\widetilde{\Xi^{-1}}[v]\defi \Xi^{-1}\circ \mcal{X}[v]$$ and, given a general operator $\mcal{U}$, we adopt the notation $\widetilde{\mcal{U}}$ for the operator $\mcal{U}\circ\mcal{X}$ when the latter is well-defined.

\begin{theorem}\label{thm: stru inv Xi}
	Let $f\in\mcal{C}^n(\R)\cap\left(\dfrac{1}{\rho_V}H^n(\R)\right)\cap H^{1}(\mu_V)$, for all $|x|>M_V$ with $M_V$ given in Lemma \ref{lem:M_V}, for all $k\in\llbracket1,n+1\rrbracket$ it holds that	\begin{equation}\label{eq:formulainverse}
		\boxed{\left(\widetilde{\Xi^{-1}}[f]\right)^{(k)}=\sum_{i=0}^{k-1}Q_i^k\left(\theta,\dots,\theta^{(i)}\right)\beta_{k-i},\quad\quad\mrm{where}\quad \theta\defi \dfrac{\rho_V'}{\rho_V}}
	\end{equation}
	and $\widetilde{\Xi^{-1}}[f]=\beta_0$.
	The $\beta_i$'s are defined, for all $|x|>M_V$, for all $i\in\llbracket1,k\rrbracket$, by: $$\beta_i(x)\defi \dfrac{-1}{\rho_V(x)}\int\limits_x^{\mrm{sgn}(x)\infty}\diff t\rho_V(t)\mcal{O}^i\left[\mcal{X}[f]+2P\mcal{W}\circ\widetilde{\mcal{L}^{-1}}[f]\right] (t)$$ (where $\mcal{O}$ is defined in Defintion \ref{def:operateurO}). Above $Q_i^k$ denotes the unique homogeneous polynomial in $i+1$ variables with degree $k-i$, with differential degree with respect to $\theta$ equal to $a$ and with coefficients independent of $V$ satisfying the following induction relations:\begin{align}\label{eq:reccu formula}
		&Q_0^{k+1}(\theta)=\theta Q_0^k(\theta)=\theta^k
		\\&\forall i\in\llbracket1,k-1\rrbracket, Q_i^{k+1}\left (\theta,\dots,\theta^{(a)}\right )=\theta Q_i^k\left (\theta,\dots,\theta^{(a)}\right )+Q_{i-1}^k\left (\theta,\dots,\theta^{(i)}\right)'\label{eq:reccu formula2}
		\\&Q_{k}^{k+1}\left (\theta,\dots,\theta^{(k)}\right)=Q_{k-1}^k\left(\theta,\dots,\theta^{(k-1)}\right)'=\theta^{(k-1)}\label{eq:reccu formula3}
	\end{align}
\end{theorem}

\begin{proof}
	We prove this statement by induction. For $k=0$, by \eqref{eq:resolvant} and \eqref{eq:Ainverse'}, 
	by setting $$g\defi -\mcal{X}[f]-2P\mcal{W}\circ\widetilde{\mcal{L}^{-1}}[f],$$ by \eqref{eq:resolvant}, we get for all $x\in\R$, 
	\begin{equation}\label{eq:inverseXisansipp}
		\widetilde{\Xi^{-1}}[f](x)=\left(\widetilde{\mcal{L}^{-1}}[f]\right)'(x) =\left(\mcal{A}^{-1}[g]\right)'(x)=\dfrac{1}{\rho_V(x)}\int\limits_x^{\mrm{sgn}(x)\infty}\diff t\rho_V(t)g(t)=\beta_0(x).
	\end{equation}
	For $k=1$, differentiating again, which is allowed by Lemma \ref{lem:reginverse}, we get for $|x|$ large enough:
	$$(\widetilde{\mcal{L}^{-1}}[f])''(x)=(\mcal{A}^{-1}[g])''(x)=-g(x)-\dfrac{\rho_V'}{\rho_V}(x)\left(\mcal{A}^{-1}[g]\right)'(x).$$
	After integrating by parts in the last integral, we obtain:
	$$\left(\widetilde{\Xi^{-1}}[f]\right)'(x)= \left (\widetilde{\mcal{L}^{-1}}[f]\right )''(x)=\dfrac{\rho_V'}{\rho_V^2}(x)\int\limits_x^{\mrm{sgn}(x)\infty}\diff t\rho_V(t)\left(g\dfrac{\rho_V}{\rho_V'}\right)' (t)=(\theta\beta_1)(x).$$
	Since $Q_0^1(\theta)= \theta$, which is readily seen of degree $1$ and of differential degree with respect to $\theta$ equal to $0$.
	
	Now, let $k\in\llbracket1,n\rrbracket$ and suppose that \eqref{eq:formulainverse} is true at rank $k$, then by differentiating we get:
	\begin{equation}\label{eq:provingrec}
		\left(\widetilde{\Xi^{-1}}[f]\right) ^{(k+1)}=\sum_{i=0}^{k-1}Q_i^{k}\left(\theta,\dots,\theta^{(a)}\right)\beta_{k-i}'+Q_i^{k}\left(\theta,\dots,\theta^{(i)}\right)'\beta_{k-i}.
	\end{equation}
	First, for $i\in\llbracket1,k\rrbracket$ and $|x|\geq M_V$, we have:
	$$\beta_i'(x)=-\mcal{O}^i[g](x)-\dfrac{\rho_V'}{\rho_V^2}(x)\int\limits_x^{\mrm{sgn}(x)\infty}\diff t\rho_V(t)\mcal{O}^i[g](t)=\dfrac{\rho_V'}{\rho_V^2}(x)\int\limits_x^{\mrm{sgn}(x)\infty}\diff t\rho_V(t)\mcal{O}^{i+1}[g](t)=(\theta\beta_{i+1})(x).$$
	The second equality follows from an integration by parts and the fact that $\dfrac{\rho_V'}{\rho_V^2}\rho_V\dfrac{\rho_V}{\rho_V'}\mcal{O}^i[g]$ goes to zero at infinity.  Hence \eqref{eq:provingrec} becomes
	\begin{multline*}
		\left(\widetilde{\Xi^{-1}}[f]\right)^{(k+1)}=\sum_{i=0}^{k-1}\theta Q_i^{k}\left(\theta,\dots,\theta^{(i)}\right)\beta_{k+1-i}+Q_i^k\left(\theta,\dots,\theta^{(a)}\right)'\beta_{k-i}
		\\=\theta Q_0^k(\theta)\beta_{k+1}+\sum_{i=1}^{k-1}\left(\theta Q_i^k\left(\theta,\dots,\theta^{(a)}\right)+Q_{i-1}^k\left(\theta,\dots,\theta^{(i-1)}\right)'\right)\beta_{k+1-i}
		\\+Q _{k-1}^k\left (\theta,\dots,\theta^{(k-1)}\right )'\beta_1.
	\end{multline*}
	By the definitions of $(Q_i^{k+1})_i$, it is clear that \eqref{eq:formulainverse} is true at rank $k+1$. The fact that $Q_i^{k+1}$ are homogeneous and have degree $k-i$ and differential degree $i$ can be checked directly from the induction relations \eqref{eq:reccu formula}, \eqref{eq:reccu formula2}, \eqref{eq:reccu formula3}. This concludes the proof.
\end{proof}

\begin{remark}
	When $V(x)=x^m$ with $m$ even, it can be checked from \eqref{deriv1} that for every $i\in\llbracket0,n-2\rrbracket$, $Q_i^n(\theta,\dots,\theta^{(i)})$ is of the form $c_{i,n}x^{m(n-1-i)-(n-1)}+T_i^n(x)+R_i^n\left(x,\mcal{H}[\rho_V],\dots,\mcal{H}\left [\rho_V^{(i)}\right ]\right)(x)$ where $c_{i,n}$ is a real number, $T_i^n$ is polynomial of degree strictly lower than $m(n-1-i)-(n-1)$ and $R_i^n$ is also a polynomial of degree greater than $1$. Since all these Hilbert transform vanish at infinity, such a polynomial expression goes to zero at infinity. This decomposition holds as long as the degree of the monomial is non-negative, otherwise it is zero. We give the first decompositions for $\left(\widetilde{\Xi^{-1}}[f]\right)^{(k)}$ for $k\in\llbracket0,3\rrbracket$:
	$$\widetilde{\Xi^{-1}}[f]=\beta_0,\hspace{2cm}\left(\widetilde{\Xi^{-1}}[f]\right)'=\theta\beta_1,\hspace{2cm}\left(\widetilde{\Xi^{-1}}[f]\right)''=\theta'\beta_1+\theta^2\beta_2,$$
	and
	$$	\left(\widetilde{\Xi^{-1}}[f]\right)^{(3)}=\theta''\beta_1+\left(\theta\theta'+(\theta^2)'\right)\beta_2+\theta^3\beta_3.$$
	With the choice of potential $V(x)=x^m$ with $m$ even, choosing a bounded function $f$ with bounded derivatives at all orders and integrating by parts, it holds that for all $k\geq0$, there exists $\gamma_0^{(k)},\dots,\gamma_k^{(k)}\in\R,$
	$$|\beta_k(x)|\underset{|x|\rightarrow\infty}{\sim}\dfrac{\mcal{O}^k[g](x)}{x^{m-1}}\underset{|x|\rightarrow\infty}{\sim}\dfrac{1}{x^{m-1}}\sum_{j=0}^{k}g^{(j)}(x)\left(\dfrac{\gamma_{j}^{(k)}}{x^{km-j}}+\underset{|x|\rightarrow\infty}{o}\left(\dfrac{1}{x^{km-j}}\right)  \right).$$
	When $V(x)=\cosh(\alpha x)$, by the same computation, we get for different $\gamma_j^{(k)}$
	$$|\beta_k(x)|\underset{|x|\rightarrow\infty}{\sim}e^{-\alpha|x|}\sum_{j=0}^{k}\gamma_{j}^{(k)}g^{(j)}(x)\left(e^{-k\alpha|x|}+\underset{|x|\rightarrow\infty}{o}\left(e^{-k\alpha|x|}\right)  \right).$$
\end{remark}

\subsection{Controls on the inverse of the master operator}
Since we are going to use the polynomials $P_i^k$ and $Q_i^k$ defined previously in Theorem \ref{thm: struc op O} and \ref{thm: stru inv Xi}, a lot in our estimates on $\widetilde{\Xi^{-1}_1}$, we first need the following lemma which gives the asymptotics of these polynomials when $|x|\rightarrow\infty$. With  $\boxed{\alpha\defi\dfrac{\rho_V}{\rho_V'}=\theta^{-1}}$, the following result holds.
\begin{lemma}\label{lem:gammaetalpha}
	For all $k\geq1$, for all $j\in\llbracket0,k\rrbracket$
	\begin{enumerate}
		\item[(i)] $P_i^{k}\left(\alpha,\dots,\alpha^{(i)}\right)(x)=\underset{|x|\rightarrow\infty}{O}\left(V'(x)^{-k}\right)$,
		\item[(ii)] $ P_i^{k}\left(\alpha,\dots,\alpha^{(i)}\right)' (x)=\underset{|x|\rightarrow\infty}{O}\left(V'(x)^{-k}\right),$
		\item[(iii)] $Q_i^{k}\left(\theta,\dots,\theta^{(i)}\right)(x)=\underset{|x|\rightarrow\infty}{O}\left(V'(x)^{k-i}\right)$.
	\end{enumerate}
\end{lemma}

\begin{proof}
	For \textit{i)}, by the Faà di Bruno's formula, for all $n\geq0$,
	$$\alpha^{(n)}=\left(\dfrac{-1}{V'+2P\mcal{H}[\rho_{V}]}\right)^{(n)}=-\sum_{\lambda\vdash n}\dfrac{(-1)^{|\lambda|}|\lambda|!}{\left(V'+2P\mcal{H}[\rho_{V}]\right)^{|\lambda|+1} }\prod_{j=1}^{n}\dfrac{\left(V^{(j+1)}+2P\mcal{H}[\rho_V]^{(j)}\right)^{\lambda_j}}{\lambda_j!(j!)^{\lambda_j}}.  $$
	where the sum is over $\underline{\lambda}\defi (\lambda_1,\dots,\lambda_n)$ such that $\sum_{j=1}^{n}j\lambda_j=n$. From assumption \textit{v)} and Lemma \ref{lemme:boundedHilbert}, we see that
	$$|\alpha^{(n)}(x)|\leq\sum_{\lambda\vdash n}C_\lambda \underset{|x|\rightarrow\infty}{O}\left(V'(x)^{-1}\right)=\underset{|x|\rightarrow\infty}{O}\left(V'(x)^{-1}\right).$$
	Hence $P_j^{k}$, as a homogeneous polynomial in $\left(\alpha,\dots,\alpha^{(j)}\right)$ of degree $k$, is a $\underset{|x|\rightarrow\infty}{O}\left(V'(x)^{-k}\right)$. 
	
	For the point \textit{ii)}, one has to notice that for each monomial
	$$A_n\defi \left[\prod_{j=1}^{n}\left(\alpha^{(j)}\right)^{k_j} \right]'=\sum_{l=1}^{n}k_l\alpha^{(l+1)}\left(\alpha^{(l)}\right)^{k_l-1}\prod_{j\neq l}^{n}\left(\alpha^{(j)}\right)^{k_j}.$$
	But, we have proven that for all $j\in\N$, $\alpha^{(j)}=\underset{|x|\rightarrow\infty}{O}\left(V'(x)^{-1}\right)$, so by denoting $k\defi \sum_{j=1}^{n}k_j=\deg(A_n)$, $A_n(x)=\underset{|x|\rightarrow\infty}{O}\left(V'(x)^{-k}\right)$. Therefore, any homogeneous polynomial of degree $k$ such as $P_i^k$ in the variables $\left(\alpha,\dots,\alpha^{(i)}\right)$ is a $\underset{|x|\rightarrow\infty}{O}\left(V'(x)^{-k}\right)$.
	
	Finally for the point \textit{iii)}, it is clear that for all $j\geq0$,
	\begin{equation*}
		\theta^{(j)}(x)=V^{(j+1)}(x)-2P\mcal{H}\left[\rho_V\right]^{(j)}(x)= \underset{|x|\rightarrow\infty}{O}\left(V'(x)\right).
	\end{equation*}
	Thus $Q_i^{(k)}\left(\theta,\dots,\theta^{(i)}\right )$ as a homogeneous polynomial of degree $k-i$, is a $\underset{|x|\rightarrow\infty}{O}\left(V'(x)^{k-i}\right)$.
\end{proof}
Now that a convenient expression for $\widetilde{\Xi^{-1}_1}$ is available thanks to Theorem \ref{thm: stru inv Xi} and that the asymptotics polynomial $Q_i^{k}$ are well-understood thanks to Lemma \ref{lem:gammaetalpha}, we can obtain the first control on $\widetilde{\Xi^{-1}_1}$ for the $H^{n}$-norm.
\begin{theorem}[$H^{n}$-continuity of $\widetilde{\Xi_1^{-1}}$]\label{thm: cont inverse master}
	Let $n,p\geq1$, there exists an explicit constant $C_{H^{n}}(\widetilde{\Xi_1^{-1}})>0$ (depending only on $n$ and $V$) such that for all $f\in H^{n+1}(\R^{p})$,  $$\boxed{\left\|\widetilde{\Xi^{-1}_1}[f]\right\|_{H^{n}(\R^{p})}\leq C_{H^{n}}(\widetilde{\Xi_1^{-1}})\cdot\|f\|_{H^{n+1}(\R^{p})}.}$$
	Under the choice of potential $V_{\phi,t}$ defined in Theorem \ref{thm:conteqdensity}, for $\phi\in H^{\infty}(\R)$ the map $t\in[0,1]\mapsto C_{H^{n}}(\widetilde{\Xi_1^{-1}})$ is continuous.
\end{theorem}

\begin{proof}
	Let $m\leq n$ and $(m_1,\dots,m_p)\in\N^p$ be such that $\sum_{i=1}^{p}m_i=m$. Let $x_2,\dots,x_p\in\R^{p-1}$ be fixed, we define $h:x_1\mapsto\partial_2^{m_2}\dots\partial_p^{m_p}f(x_1,\dots,x_p)$ and $g=-\mcal{X}[h]-2P\mcal{W}\circ\widetilde{\mcal{L}^{-1}}[h]$. With these notations, $$\partial^{\underline{m}}\widetilde{\Xi_{1}^{-1}}[f](x_1,\dots,x_p)=\widetilde{\Xi^{-1}}\left [h\right]^{(m_1)}(x_1).$$
	The idea of the proof is to use the closed forms for $\widetilde{\Xi^{-1}}[f]^{(n)}$ and $\mcal{O}^{k}$ found respectively in Theorem \ref{thm: stru inv Xi} and Theorem \ref{thm: struc op O} to get the control for $|x_1|>M_V$. The bound for $|x_1|\leq M_V$ relies on Lemma \ref{lem:majorationderivee nieme de rho Xi} and the fact that $\rho_V$ is bounded by below on this set.
	
	\textbf{Control for $|x_1|>M_V$.} For $|x_1|\geq M_V$, we can apply Theorem \ref{thm: stru inv Xi} and Theorem \ref{thm: struc op O}, to get :
	\begin{multline}\label{eq:1er formule cont inver mast}
		\widetilde{	\Xi^{-1}}[h]^{(m_1)}(x_1)=\sum_{i=0}^{m_1-1}Q_i^{m_1}\left(\theta,\dots,\theta^{(i)}\right)(x_1)\beta_{m_1-i}(x_1)
		\\=\sum_{i=0}^{m_1-1}Q_i^{m_1}\left(\theta,\dots,\theta^{(i)}\right)(x_1)\dfrac{1}{\rho_V(x_1)}\int\limits_{x_1}^{\mrm{sgn}(x_1)\infty}\diff t\rho_V(t)\mcal{O}^{m_1-i}[g](t)
		\\=\sum_{i=0}^{m_1-1}\sum_{j=0}^{m_1-i}Q_i^{m_1}\left(\theta,\dots,\theta^{(i)}\right)(x_1)\dfrac{1}{\rho_V(x_1)}\int\limits_{x_1}^{\mrm{sgn}(x_1)\infty}\diff t\rho_V(t)g^{(m_1-i-j)}(t)P_{j}^{m_1-i}(\alpha,\dots,\alpha^{(j)})(t).
	\end{multline}
	Moreover an integration by parts yields:\begin{multline}\label{eq:2eme formule cont inver mast}
		\widetilde{\Xi^{-1}}[h]^{(m_1)}(x_1)=\sum_{i=0}^{m_1-1}\sum_{j=0}^{m_1-i}Q_i^{m_1}\left(\theta,\dots,\theta^{(i)}\right)(x_1)\Bigg(-g^{(m_1-i-j)}(x_1)\alpha(x_1)P_{j}^{m_1-i}(\alpha,\dots,\alpha^{(j)})(x_1)
		\\+\dfrac{1}{\rho_V(x_1)}\int\limits_{x_1}^{\mrm{sgn}(x_1)\infty}\diff t\rho_V(t)\Bigg[g^{(m_1-i-j+1)}(t)\alpha(t)P_{j}^{m_1-i}\left (\alpha,\dots,\alpha^{(j)}\right )(t)
		\\+g^{(m_1-i-j)}(t)\left[\alpha P_{j}^{m_1-i}\left (\alpha,\dots,\alpha^{(j)}\right)  \right] '(t) \Bigg]\Bigg) .
	\end{multline}
	From Lemma \ref{lem:M_V}, the functions $P^{i}_{j}\left(\alpha,\dots,\alpha^{(j)}\right) $ that appear above are well-defined on $[-M_V,M_V]^c$ \textit{i.e.} they don't have any singularity on this set.  Hence by integrating with respect to $x_1$, $\left ({{\Xi}^{-1}[h]^{(m_1)}}\right)^2$ on $[M_V,+\infty[$, we get by Jensen's inequality for a constant $C>0$ depending only on $m_1$
	\begin{multline*}
		\int_{M_V}^{+\infty}\diff x_1\left(	\widetilde{{\Xi}^{-1} }[h]^{(m_1)}(x_1)\right)^2
		\\\leq C\sum_{i=0}^{m_1-1}\sum_{j=0}^{m_1-i}\int\limits_{M_V}^{+\infty}\diff x_1Q_i^{m_1}\left(\theta,\dots,\theta^{(i)}\right)(x_1) ^2\Bigg\{g^{(m_1-i-j)}(x_1)^2\alpha(x_1)^2P_{j}^{m_1-i}(\alpha,\dots,\alpha^{(j)})(x_1)^2
		\\+\dfrac{1}{\rho_V(x_1)^2}\left[\int\limits_{x_1}^{\mrm{sgn}(x_1)\infty}\diff t\rho_V(t)g^{(m_1-i-j+1)}(t)\alpha(t)P_{j}^{m_1-i}\left (\alpha,\dots,\alpha^{(j)}\right )(t)\right]^2
		\\+\dfrac{1}{\rho_V(x_1)^2}\left[\int\limits_{x_1}^{\mrm{sgn}(x_1)\infty}\diff t\rho_V(t)g^{(m_1-i-j)}(t)\left[\alpha P_{j}^{m_1-i}\left (\alpha,\dots,\alpha^{(j)}\right)  \right] '(t)\right]^2\Bigg\}.
	\end{multline*}
	We want to bound this expression by $\|g\|_{H^{m_1+1}(\R)}^2$, but since $g=-h+2P\mcal{H}\left[\rho_V\Xi^{-1}[h]\right]+\mfrak{c}_h,$
	where $$\mfrak{c}_h=\int_\R h(y)\diff\mu_V(y)-2P\int_\R\mcal{H}\left[\rho_V\Xi^{-1}[h]\right](y)\diff\mu_V(y),$$ the constant terms will fail to be in $L^2(\R)$. We thus have to treat these terms separately. In the previous sum, $g$ is differentiated everywhere except in the term $j=m_1-i$ so this is the only value of $j$ where we have to deal with $\mfrak{c}$. By defining for $\bm{i}=(i_1,i_2)$:
	\begin{align}\label{eq:mfrak f1,2}
		&\mfrak{f}^{(1)}_{m,\bm{i}}:x\mapsto Q_{i_1}^{m}\left(\theta,\dots,\theta^{(i_1)}\right)(x)\alpha(x) P_{i_2}^{m-i_1}(\alpha,\dots,\alpha^{(i_2)})(x)
		\\&\mfrak{f}_{m,\bm{i}}^{(2)}:x\mapsto  \dfrac{Q_{i_1}^{m}\left(\theta,\dots,\theta^{(i_1)}\right)(x)}{\rho_V(x)}\int\limits_{x}^{\mrm{sgn}(x)\infty}\diff t\rho_V(t)\left[\alpha P_{i_2}^{m-i_1}\left (\alpha,\dots,\alpha^{(i_2)}\right)  \right] '(t),
	\end{align}
	keeping the $V$-dependence implicit in those functions, Jensen's inequality and inequality \eqref{ineq:contXi} yields:
	\begin{align*}
		\mfrak{c}_h^2C\sum_{\substack{\bm{i}\\i_1+i_2=m_1\\0\leq i_1\leq m_1-1}}\int_{M_V}^{+\infty}\diff x_1&\left[\mfrak{f}^{(1)}_{m_1,\bm{i}}(x_1)^2+\mfrak{f}_{m_1,\bm{i}}^{(2)}(x_1)^2\right]
		\\&\leq C\max_{\bm{i}}\max_{1\leq j\leq 2}\max_{m\leq n}\|\mfrak{f}^{(j)}_{m,\bm{i}}\|_{L^2([-M_V,M_V]^c)}^2 \left(\|h\|^2_{L^2(\mu_V)}+\left \|\mcal{H}\left[\rho_V\Xi^{-1}[h]\right]\right \|^2_{L^2(\mu_V)}\right) 
		\\&\leq C(V)\|\rho_{V}\|_{L^\infty(\R)}
		\left(\|h\|_{L^2(\R)}^2+\|\rho_{V}\|_{L^\infty(\R)}\pi^2\|\Xi^{-1}[h]\|_{L^2(\mu_V)}^2\right) 
		\\&\leq C(V)
		\left(\|h\|_{L^2(\R)}^2+\|\rho_{V}\|_{L^\infty(\R)}\pi^2C_{\mcal{L}}^2\|h'\|_{L^2(\mu_V)}^2\right) 
		\\&\leq C_1(V)\|\partial_2^{m_2}\dots\partial_p^{m_p}f(.,x_2,\dots,x_p)\|_{H^1(\R)}^2
	\end{align*}
	where at the end, the constant $C_1(V)$ is defined, for a constant $C>0$ only depending on $n$, by:
	\begin{equation}\label{eq:C_1(V,n_1)}
		C_1(V)\defi C\max_{\bm{i}}\max_{1\leq j\leq 2}\max_{m\leq n}\|\mfrak{f}^{(j)}_{m,\bm{i}}\|_{L^2([-M_V,M_V]^c)}^2C_{\mcal{L}}^2\left(1+\|\rho_{V}\|_{L^\infty(\R)}^3\right).
	\end{equation}
	Above, the first integral that appears is well-defined, since by Lemma \ref{lem:gammaetalpha}, one can check by assumption \textit{\ref{assumption5}} that $\mfrak{f}^{(1)}_{m_1,\bm{i}}(x)$ and $\mfrak{f}^{(2)}_{m_1,\bm{i}}(x)$ behave like $\underset{|x|\rightarrow\infty}{O}\left(V'(x)^{-2}\right)$ which is integrable by assumption \textit{\ref{assumption5}} again.
	
	In the following, we set $\mfrak{g}\defi g-\mfrak{c}_h$. We can now replace $g'$ by $\mfrak{g}'$ since we handled all the terms involving $\mfrak{c}_h$. By Cauchy-Schwarz inequality, with $C_2(V)$ defined, for a constant $C>0$ only depending on $n$ by:
	\begin{equation}\label{eq:def C_2(V,n_1)}
		C_2(V)=C\max_{\bm{i}}\max_{m\leq n}\left(\Big\|\mfrak{f}^{(1)}_{m,\bm{i}}\Big\|^2_{L^\infty([-M_V,M_V]^c)}+\max_{3\leq j\leq 4}\|\mfrak{f}^{(j)}_{m,\bm{i}}\|_{L^2([-M_V,M_V]^c)}^2\right),
	\end{equation}
	with 	\begin{align}\label{eq:def mfrak{f}}
		&\mfrak{f}^{(3)}_{m,\bm{i}}:x\mapsto  \dfrac{Q_{i_1}^{m}\left(\theta,\dots,\theta^{(i_1)}\right)(x)}{\rho_V(x)}\Bigg|\int\limits_{x}^{\mrm{sgn}(x)\infty}\diff t\rho_V(t)^2\alpha(t)^2P_{i_2}^{m-i_1}\left (\alpha,\dots,\alpha^{(i_2)}\right )(t)^2\Bigg|^{1/2},
		\\&\mfrak{f}^{(4)}_{m,\bm{i}}:x\mapsto\dfrac{Q_{i_1}^{m}\left(\theta,\dots,\theta^{(i_1)}\right)(x) }{\rho_V(x)}\Bigg|\int\limits_{x}^{\mrm{sgn}(x)\infty}\diff t\rho_V(t)^2\left[\alpha P_{i_2}^{m-i_1}\left (\alpha,\dots,\alpha^{(i_2)}\right)  \right] '(t)\Bigg|^{1/2},
	\end{align}
	we get, with $|\bm{i}|=i_1+i_2$:
	\begin{multline*}
		\sum_{i_1=0}^{m_1-1}\sum_{i_2=0}^{m_1-i_1}\int_{M_V}^{+\infty}\diff x_1Q_{i_1}^{m_1}\left(\theta,\dots,\theta^{(i_1)}\right)(x_1) ^2\Bigg\{\mfrak{g}^{(m_1-|\bm{i}|)}(x_1)^2\alpha(x_1)^2P_{i_2}^{m_1-i_1}(\alpha,\dots,\alpha^{(i_2)})(x_1)^2
		\\+\dfrac{1}{\rho_V(x_1)^2}\left[\int\limits_{x_1}^{\mrm{sgn}(x_1)\infty}\diff t\rho_V(t)\mfrak{g}^{(m_1-|\bm{i}|+1)}(t)\alpha(t)P_{i_2}^{n_1-i_1}\left (\alpha,\dots,\alpha^{(i_2)}\right )(t)\right]^2
		\\+\dfrac{1}{\rho_V(x_1)^2}\left[\int\limits_{x_1}^{\mrm{sgn}(x_1)\infty}\diff t\rho_V(t)\mfrak{g}^{(m_1-|\bm{i}|)}(t)\left[\alpha P_{i_2}^{m_1-i_1}\left (\alpha,\dots,\alpha^{(i_2)}\right)  \right] '(t)\right]^2\Bigg\}
		\\\leq C_2(V)\|\mfrak{g}\|_{H^{m_1+1}(\R)}^2.
	\end{multline*}
	Finally, by using that $\mfrak{g}=-h-2P\mcal{H}\left[\rho_V\Xi^{-1}[h]\right],$ that $\pi^{-1}\mcal{H}$ is an isometry in $L^2(\R)$ and that for all $u\in H^1(\R)$, $\mcal{H}[u]'=\mcal{H}[u']$, we obtain for a universal constant $C>0$,
	\begin{equation*}
		\left\|{\Xi}^{-1}[h]^{(m_1)}\right \|_{L^2([M_V,+\infty[)}\leq C\max_{i=1,2} C_i(V)^{1/2}\Big(\|h\|_{H^{m_1+1}(\R)}+\|\rho_V\Xi^{-1}[h] \|_{H^{m_1+1}(\R)}\Big).
	\end{equation*}
	We now use the form stated in Lemma \ref{lem:majorationderivee nieme de rho Xi}, to conclude that
	$$	\left\|{\Xi}_{1}^{-1}[h]^{(m_1)}\right \|_{L^2([M_V,+\infty[)}\leq C \max_{i=1,2} C_i(V)^{1/2}\left( C_3(V,n)+1\right) \|h\|_{H^{m_1+1}(\R)}.$$
	where $C_3(V,n)$ is given by Lemma \ref{lem:majorationderivee nieme de rho Xi}. The exact same bounds holds on $]-\infty,-M_V]$. 
	Now relaxing the dependence on $x_2,\dots,x_p\in\R$ and integrating with respect to these variables, we get
	\begin{multline*}
		\left\|\partial^{\underline{m}}{\Xi}_{1}^{-1}[f]\right \|_{L^2([-M_V,M_V]^c\times\R^{p-1})}\leq C\max_{i=1,2} C_i(V)^{1/2}\left(C_3(V,n)+1\right)
		\\\times\left(\int_{\R^{p-1}}\|\partial_2^{m_2}\dots\partial_p^{m_p}f(.,x_2,\dots,x_p)\|_{H^{m_1+1}(\R)}^2\diff x_2\dots \diff x_{p}\right)^{1/2}.
	\end{multline*}
	Thus we deduce that for a constant $C>0$ only depending on $n$ such that
	\begin{equation}\label{eq:presquelafindansXihn}
		\left\|\Xi_{1}^{-1}[f]\right \|_{H^{ n}([-M_V,M_V]^c\times\R^{p-1})}\leq C\max_{m_1\leq n}\max_{i=1,2} C_i(V)^{1/2}\left(C_3(V,n)+1\right)\|f\|_{H^{ n+1}(\R^p)}.
	\end{equation}
	\textbf{Control for $|x_1|<M_V$.}
	Now, we prove the control on $[-M_V,M_V]\times\R^{p-1}$, we fix $x_2,\dots,x_p\in\R$. By Cauchy-Schwarz inequality:
	\begin{multline*}
		\int_{-M_V}^{M_V}\left |\Xi^{-1}\left[h\right]^{(m_1)}(x)\right |^2\diff x\leq\left\|\rho_V^{-1}\right\|_{L^\infty([-M_V,M_V])}^2\left\|\rho_V \Xi^{-1}\left[h\right]^{(m_1)}\right\|_{L^2(\R)}^2
		\\\leq C_3(V,n)^2\left\|\rho_V^{-1}\right\|_{L^\infty([-M_V,M_V])}^2 \left\| h\right\|_{H^{m_1}(\R)}^2	
	\end{multline*}
	where the last inequality comes from Lemma \ref{lem:majorationderivee nieme de rho Xi}. Again relaxing the dependence on $x_2,\dots,x_p\in\R$ and integrating with respect to these variables, we get for a constant $C>0$ dependent only on $n$:
	$$\left\|\Xi_{1}^{-1}[f]\right \|_{H^{ n}([-M_V,M_V]\times\R^{p-1})}\leq C\max_{m_1\leq n}C_3(V,n)\left\|\rho_V^{-1}\right\|_{L^\infty([-M_V,M_V])}\|f\|_{H^{ n+1}(\R^p)}.$$
	Collecting the last bound and \eqref{eq:presquelafindansXihn} leads to the following definition of $	C_{H^{n}}(\widetilde{\Xi_1^{-1}})$ for a constant $C>0$ only depending on $n$:
	\begin{equation}\label{def:constanteXiHn}
		C_{H^{n}}(\widetilde{\Xi_1^{-1}})\defi C\left( C_3(V,n)+1\right)\left[\left\|\rho_V^{-1}\right\|_{L^\infty([-M_V,M_V])}+\max_{i=1,2} C_i(V)^{\frac12}\right].
	\end{equation}
	The fact that, upon choosing the potential $V_{\phi,t}$ with $\phi\in H^{\infty}(\R)$, $t\mapsto	C_{H^{n}}(\widetilde{\Xi_1^{-1}})$ is continuous is shown in Proposition \ref{prop:constantescontinues}.
\end{proof}

We now prove the analogous result for the norm $W_n^{\infty}$ defined in \eqref{def:sobolev norms} and the associated space of functions $W_n^{\infty}$ defined in \eqref{def:sobolev spaces}.

\begin{theorem}
	[$W_n^\infty$-continuity of $\widetilde{\Xi_1^{-1}}$]\label{thm: cont inverse master linfini}
	Let $n,p\geq1$, there exists an explicit constant $C_{W_n^{\infty}}(\widetilde{\Xi_1^{-1}})>0$ (depending only on $n$ and $V$), such that for all $f\in W_{n+1}^{\infty}(\R^{p})$, $$\boxed{\big\|\widetilde{\Xi^{-1}_1}[f]\big\|_{W^\infty_{n}(\R^{p})}\leq C_{W_n^{\infty}}(\widetilde{\Xi_1^{-1}})\cdot\|f\|_{W_{n+1}^\infty(\R^{^p})}.}$$
	Under the choice of potential $V_{\phi,t}$ defined in Theorem \ref{thm:conteqdensity}, for $\phi\in H^{\infty}(\R)$ the map $t\in[0,1]\mapsto C_{W_n^{\infty}}(\widetilde{\Xi_1^{-1}})$ is continuous.
\end{theorem}

\begin{proof}
	Let $f\in W_{n+1}(\R^p)$, let $m\leq n$ and $(m_1,\dots,m_p)\in\N^p$ be such that $\sum_{i=1}^{p}m_i=m$, let $x_2,\dots,x_p\in\R^{p-1}$, we set $h:x_1\mapsto\partial_2^{m_2}\dots\partial_p^{m_p}f(x_1,\dots,x_p)$ we know by theorem \ref{thm: stru inv Xi} that
	\begin{equation}\label{def:fct g}
		\widetilde{\Xi^{-1}}[h](x_1)=\dfrac{1}{\rho_V(x_1)}\int\limits_{x_1}^{\mrm{sgn}(x_1)\infty}\diff t\rho_V(t)g(t),\quad\quad\quad \text{with}\quad g\defi -\mcal{X}[h]+2P\mcal{X}\circ\mcal{H}\big[\rho_V\widetilde{\Xi^{-1}}[h]\big].
	\end{equation}
	The idea of the proof is to consider separatly the case with $m_1=0$ and $m_1\neq0$. The first case follows from straightforward bounds. For the case $m_1\neq0$ and $|x_1|>M_V$, we rely on the formula \eqref{eq:1er formule cont inver mast} obtained in the previous proof. We successively use Lemma \ref{lem:gammaetalpha} and Lemma \ref{lemme:boundedHilbert} as well as the isometry property of the Hilbert transform to deduce the desired bound. For the case $m_1\neq0$ and $|x_1|\leq M_V$, we use Leibniz formula in \eqref{eq:inverseXisansipp} to deduce a manageable expression \eqref{eq:Xi-1deriveelinfini} which allows to deduce this last bound and conclude the continuity result for $\widetilde{\Xi^{-1}_1}[f]$.
	
	\textbf{Case $m_1=0$.} For the following, we define
	\begin{equation}\label{eq:defI_a}
		\mcal{I}_1:x\mapsto \dfrac{1}{\rho_V(x)}\left|\int\limits_x^{\mrm{sgn}(x)\infty}\rho_V(t)\diff t\right |, \hspace{2cm}\mcal{I}_2:x\mapsto \dfrac{1}{\rho_V(x)}\left|\displaystyle\int\limits_x^{\mrm{sgn}(x)\infty}\rho_V(t)^2\diff t\right|^{1/2}.
	\end{equation}
	By integration by parts, one can see that $\mcal{I}_1^{V}(x)=\underset{|x|\rightarrow\infty}{O}\left (V'(x)^{-1}\right )$ is bounded on $\R$. So for the first and third term, by direct bounds:
	$$\Bigg|\dfrac{1}{\rho_V(x_1)}\int\limits_{x_1}^{\mrm{sgn}(x_1)\infty}\diff t\rho_V(t)\left(-h(t)+\int_\R h(s)\diff\mu_V(s)\right) \Bigg|\leq 2 \|\mcal{I}_1\|_{\infty}\|h\|_{\infty}\leq 2 \|\mcal{I}_1\|_{\infty}\|f\|_{W_n^\infty(\R^p)}.$$
	For the two last terms, we want to use that $\pi^{-1}\mcal{H}$ is an isometry on $L^2(\R)$, so we use Cauchy-Schwarz inequa,lity and the fact that $\mcal{I}_2(x)=\underset{|x|\rightarrow\infty}{O}\left (V'(x)^{-1/2}\right )$ is bounded on $\R$ so that:
	\begin{equation*}
		\sup_{x_1\in\R}\Bigg|\dfrac{2P}{\rho_V(x_1)}\int\limits_{x_1}^{\mrm{sgn}(x_1)\infty}\diff t\rho_V(t)\mcal{X}\circ\mcal{H}\left[\rho_V\widetilde{\Xi^{-1}}[h]\right ](t)\Bigg|\leq C_4(V)\left\|h'\right \|_{\infty}\leq C_4(V)\|f\|_{W_{n+1}^\infty(\R^p)}.
	\end{equation*}
	with a universal constant $C>0$,
	\begin{equation}\label{eq:C_4(V)C_5(V)}
		C_4(V)\defi C\max_{1\leq i \leq2}\|\mcal{I}_i\|_{\infty} \|\rho_{V}\|_{\infty}C_{\mcal{L}}.
	\end{equation}
	Thus, by taking the supremum of $x_2,\dots,x_p\in\R^{p-1}$, we conclude that for $m_1=0$,
	\begin{equation}\label{eq:bound avec C_4 et C_5}
		\Big\|\partial^{\underline{m}}\widetilde{\Xi_1^{-1}}[f]\Big\|_{L^\infty(\R^p)}\leq \left(2 \|\mcal{I}_1\|_{\infty}+C_4(V)\right)\cdot\|f\|_{W_{n+1}^\infty(\R^p)}.
	\end{equation}
	\textbf{Case $m_1\neq0$, $|x_1|>M_V$.} For $m_1\neq0$ by \eqref{eq:1er formule cont inver mast}, for all $|x_1|>M_V$, by Lemma \ref{lem:M_V} we have:
	\begin{multline*}
		\widetilde{\Xi^{-1}}[h]^{(m_1)}(x_1)=\sum_{i=0}^{m_1-1}\sum_{j=0}^{m_1-i}\dfrac{Q_{i}^{m_1}\left(\theta,\dots,\theta^{(a)}\right)(x_1)}{\rho_V(x_1)}\int\limits_{x_1}^{\mrm{sgn}(x_1)\infty}g^{(m_1-i-j)}(t)P_{j}^{m_1-i}(\alpha,\dots,\alpha^{(j)})(t)\diff\mu_V(t)
	\end{multline*}
	where $g\defi -\mcal{X}[h]-2P\mcal{W}\circ\widetilde{\mcal{L}^{-1}}[h]$.
	Furthermore, setting just as before $$\mfrak{c}_h= \displaystyle\int_\R h(t)\diff\mu_V(t)-2P\int_\R\mcal{H}\left[\rho_V\widetilde{\Xi^{-1}}[h]\right](t)\diff\mu_V(t)$$
	and for $\bm{i}=(i_1,i_2)$:
	\begin{equation}\label{def:f^5}
		\mfrak{f}^{(5)}_{m,\bm{i}}:x\mapsto \dfrac{Q_{i_1}^{m}\left(\theta,\dots,\theta^{(i_1)}\right)(x)}{\rho_V(x)}\int\limits_{x}^{\mrm{sgn}(x)\infty}\Big|P_{i_2}^{m-i_1}(\alpha,\dots,\alpha^{(i_2)})(t)\Big|\rho_V(t)\diff t,
	\end{equation}
	we can bound every term in the previous sum, involving $\mfrak{c}_h$, \textit{i.e.} for index $\bm{i}$ such that $i_2=m_1-i_1$, with a universal constant $C>0$: 
	\begin{equation*}
		|\mfrak{c}\mfrak{f}^{(5)}_{m_1,\bm{i}}(x_1)|\leq C(1+\|\rho_V\|_{\infty}C_{\mcal{L}})\cdot\max_{\bm{i}}\|\mfrak{f}^{(5)}_{m_1,\bm{i}}\|_{L^\infty([-M_V,M_V]^c)}\cdot\|f\|_{W_{n+1}^\infty(\R^{p})}.
	\end{equation*}
	We directly bound $\mfrak{c}_h$ in the LHS, while to bound $\mfrak{f}^{(5)}_{m_1,\bm{i}}$, we successfully applied Jensen's inequality, used the isometry property of $\pi^{-1}\mcal{H}$ on $L^2(\R)$ and used the inequality \eqref{ineq:contXi}. Furthermore, the fact that $\mfrak{f}^{(5)}_{m_1,\bm{i}}$ is bounded on $[-M_V,M_V]^c$ comes from Lemma \ref{lem:gammaetalpha}. Finally, by setting $\mfrak{g}\defi g-\mfrak{c}$, it only remains to establish the following bounds:
	\begin{align*}
		\Big|\sum_{i_1=0}^{m_1-1}\sum_{i_2=0}^{m_1-i_1}&Q_{i_1}^{m_1}\left(\theta,\dots,\theta^{(i_1)}\right)(x_1)\dfrac{1}{\rho_V(x_1)}\int\limits_{x_1}^{\mrm{sgn}(x_1)\infty}\mfrak{g}^{(m_1-|\bm{i}|)}(t)P_{i_2}^{m_1-i_1}(\alpha,\dots,\alpha^{(i_2)})(t)\rho_V(t)\diff t\Big|
		\\&\leq\sum_{i_1=0}^{m_1-1}\sum_{i_2=0}^{m_1-i_1}\Bigg\{\Bigg \|\dfrac{\Big|Q_{i_1}^{m_1}\left(\theta,\dots,\theta^{(i_1)}\right)\Big|}{\rho_V}\int\limits_{.}^{\mrm{sgn}(.)\infty}\Bigg|\mcal{H}\left[\left(\rho_V\widetilde{\Xi^{-1}}[h]\right)^{(m_1-|\bm{i}|)}\right] (t)\Bigg|
		\\&\times\Big|P_{i_2}^{m_1-i_1}(\alpha,\dots,\alpha^{(i_2)})(t)\Big|\rho_V(t)\diff t\Bigg \|_{L^\infty([-M_V,M_V]^c)}+\|h\|_{W^\infty_{m_1}(\R)}\|\mfrak{f}^{(5)}_{m_1,\bm{i}}\|_{L^\infty([-M_V,M_V]^c)}\Bigg\}
		\\&\leq\sum_{i_1=0}^{m_1-1}\sum_{i_2=0}^{m_1-i_1}\Bigg\{\|\mfrak{f}^{(5)}_{m_1,\bm{i}}\|_{L^\infty([-M_V,M_V]^c)}\cdot\|h\|_{W^\infty_{m_1}(\R)}
		\\&+\|\mfrak{f}^{(6)}_{m_1,\bm{i}}\|_{L^\infty([-M_V,M_V]^c)}\left\|\mcal{H}\left[\left(\rho_V\widetilde{\Xi^{-1}}[h]\right)^{(m_1-|\bm{i}|)}\right]\right\|_{L^2(\R)} \Bigg\}
	\end{align*}
	with $|\bm{i}|=i_1+i_2$ and
	\begin{equation}\label{eq:defdepleindeconstante}
		\mfrak{f}^{(6)}_{m,\bm{i}}:x\mapsto\Big|Q_{i_1}^{m}\left(\theta,\dots,\theta^{(i_1)}\right)(x)\Big|\dfrac{1}{\rho_V(x)}\sqrt{\int\limits_{x}^{\mrm{sgn}(x)\infty}\Big|P_{i_2}^{m-i_1}(\alpha,\dots,\alpha^{(i_2)})(t)\Big|^2\rho_V(t)^2\diff t}.
	\end{equation}
	For each $\bm{i}$, $\mfrak{f}^{(6)}_{m,\bm{i}}$ is bounded on $[-M_V,M_V]^c$ because of Lemma \ref{lem:gammaetalpha} and Lemma \ref{lemme:boundedHilbert}. By Cauchy-Schwarz inequality and Lemma \ref{lem:majorationderivee nieme de rho Xi}, we get:
	\begin{equation*}
		\left\|\mcal{H}\left[\left(\rho_V\widetilde{\Xi^{-1}}[h]\right)^{(m_1-|\bm{i}|)}\right]\right\|_{L^2(\R)}\leq\pi\Big\|\rho_V\widetilde{\Xi^{-1}}[h]\Big\|_{H^{m_1}(\R)}\leq \pi C_5(V,m_1)\cdot\|f\|_{W_{n+1}^\infty(\R^{p})}.
	\end{equation*}
	Finally by the same reasonnings as before, we get
	\begin{multline*}
		\Big|\sum_{i_1=0}^{m_1-1}\sum_{i_2=0}^{m_1-i_1}\dfrac{Q_{i_1}^{m_1}\left(\theta,\dots,\theta^{(i_1)}\right)(x_1)}{\rho_V(x_1)}\int\limits_{x_1}^{\mrm{sgn}(x_1)\infty}\mfrak{g}^{(m_1-|\bm{i}|)}(t)P_{i_2}^{m_1-i_1}(\alpha,\dots,\alpha^{(i_2)})(t)\rho_V(t)\diff t\Big|
		\\\leq C_6(V,m_1)\|f\|_{W_{n+1}^\infty(\R^{p})}
	\end{multline*}
	where $C_6(V)$ is defined, for a constant $C>0$ depending only on $n$, by:
	\begin{equation}\label{eq: C_6(V,l)}
		C_6(V)\defi C\max_{1\leq l \leq n}(1+C_5(V,l))\max_{5\leq j \leq 6}\max_{\bm{i}}\|\mfrak{f}^{(j)}_{l,\bm{i}}\|_{L^\infty([-M_V,M_V]^c)}.
	\end{equation}
	Thus, we deduce that for a constant $C>0$ depending only on $n$:
	\begin{multline}\label{eq:majoration sur (-m,m)}
		\|\widetilde{\Xi^{-1}}[h]^{(m_1)}\|_{L^\infty([-M_V,M_V]^c)} \\\leq C\Big[(1+\|\rho_V\|_{\infty}C_{\mcal{L}})\max_{m\leq n}\max_{\bm{i}}\|\mfrak{f}^{(5)}_{m,\bm{i}}\|_{L^\infty([-M_V,M_V]^c)}
		+C_6(V)\Big]\cdot\|f\|_{W_{n+1}^\infty(\R^{p})}.
	\end{multline}
	
	\textbf{Case $m_1\neq0$, $|x|\leq M_V$.} Now let $x\in[-M_V,M_V]$, by differentiating $m_1$ times $\eqref{eq:inverseXisansipp}$, the Leibniz formula ensures that there exists polynomials $R_{m_1-i}^{m_1}$ depending on $(\theta,\dots,\theta^{(m_1-1-i)})$ and a polynomial $S^{m_1}$ of degree $m_1-1$ depending on $(\theta,\dots,\theta^{(m_1)})$, whose coefficients are independent of $V$ such that
	\begin{multline}\label{eq:Xi-1deriveelinfini}
		\widetilde{\Xi^{-1}}[h]^{(m_1)}(x)=\dfrac{S^{m_1}(\theta,\dots,\theta^{(m_1)})(x)}{\rho_V(x)}\int\limits_{x}^{\mrm{sgn}(x)\infty}g(t)\rho_V(t)\diff t
		\\+\sum_{i=0}^{m_1-2}R_{m_1-i}^{m_1}\left(\theta,\dots,\theta^{(m_1-1-i)}\right)(x)g^{(i)}(x)-g^{(m_1-1)}(x).
	\end{multline}
	We recall that the function $g$ was defined in \eqref{def:fct g}. Then, for all, $x\in[-M_V,M_V]$, by the same bounds as before with $C_7(V)$ defined by
	\begin{multline}\label{eq:C_7(V,l)}
		C_7(V)\defi \max_{1\leq l \leq n}\Big(\|S^l(\theta,\dots,\theta^{(l)})\|_{L^\infty([-M_V,M_V])}\|\rho_V^{-1}\|_{L^\infty([-M_V,M_V])}
		\\+l\max_{0\leq i\leq l-2}\|R_{l-i}^{l}\big(\theta,\dots,\theta^{(l-1-i)}\big)\|_{L^\infty([-M_V,M_V])}+1\Big),
	\end{multline}
	we obtain:
	\begin{multline*}
		\Big|\dfrac{S^{m_1}(\theta,\dots,\theta^{(m_1)})(x_1)}{\rho_V(x_1)}\int\limits_{x_1}^{\mrm{sgn}(x_1)\infty}\mcal{X}[h](t)\rho_V(t)\diff t
		\\+\sum_{a=0}^{m_1-2}R_{m_1-a}^{m_1}\left(\theta,\dots,\theta^{(m_1-1-a)}\right)(x_1)\mcal{X}[h]^{(a)}(x_1)
		-\mcal{X}[h]^{(m_1-1)}(x_1)\Big|\leq2C_7(V)\cdot \|f\|_{W_{n+1}^{\infty}(\R^{p})}
	\end{multline*}
	and for a universal constant $C>0$:
	\begin{multline*}
		2P\Big|\dfrac{S^{m_1}(\theta,\dots,\theta^{(m_1)})(x_1)}{\rho_V(x_1)}\int\limits_{x_1}^{\mrm{sgn}(x_1)\infty}\rho_V(t)\diff t+R_{m_1}^{m_1}\left(\theta,\dots,\theta^{(m_1-1)}\right)(x_1)\Big|\cdot\int_\R |\mcal{H}\left[\rho_V\widetilde{\Xi^{-1}}[h]\right](t)|\diff\mu_V(t)
		\\\leq C\|\rho_{V}\|^{1/2}_{\infty}C_{\mcal{L}}C_7(V)\cdot\|f\|_{W_{n+1}^{\infty}(\R^{p})}.
	\end{multline*}
	It remains to bound the terms involving the Hilbert-transform. For that, we use that for all $\phi\in H^{m_1}(\R)$ and $i\in\llbracket0,m_1-1\rrbracket$, for a universal constant $C>0$:
	$$|\mcal{H}[\phi]^{(i)}(x)|=\sqrt{\mcal{H}[\phi^{(i)}](x)^2}=\sqrt{\int_{+\infty}^{x}2\mcal{H}[\phi^{(i)}](t)\mcal{H}[\phi^{(i+1)}](t)\diff t}\leq C\sqrt{\|\phi^{(i)}\|_{L^2(\R)}\|\phi^{(i+1)}\|_{L^2(\R)}}.$$
	Applying those results for $\phi=\rho_V\widetilde{\Xi^{-1}}[h]\in H^{m_1}(\R)$, Lemma \ref{lem:reginverse} allows us to conclude that $$\|\mcal{H}[\phi]\|_{W_{m_1-1}^\infty(\R)}\leq\sqrt{2}\pi\|\phi\|_{H^{m_1}(\R)}.$$ 
	We conclude by Lemma \ref{lem:majorationderivee nieme de rho Xi} that for a universal constant $C>0$:
	$$\|\mcal{H}[\phi]\|_{W_{m_1-1}^\infty(\R)}\leq\sqrt{2}\pi C_5(V,n)\|h\|_{W_{m_1}^\infty(\R)}\leq C C_5(V,n)\cdot\|f\|_{W_{n+1}^\infty(\R^{p})} $$
	and thus, with $C_7(V)$ defined in \eqref{eq:C_7(V,l)}, we get :
	\begin{multline*}
		2P\Big|\dfrac{S^{m_1}(\theta,\dots,\theta^{(m_1)})(x_1)}{\rho_V(x_1)}\int\limits_{x_1}^{\mrm{sgn}(x_1)\infty}\mcal{H}\left[\rho_V\widetilde{\Xi^{-1}}[h]\right](t)\rho_V(t)\diff t
		\\+\sum_{i=0}^{m_1-2}R_{m_1-i}^{m_1}\left(\theta,\dots,\theta^{(m_1-1-i)}\right)(x_1)\mcal{H}\left[\rho_V\widetilde{\Xi^{-1}}[h]\right]^{(i)}(x_1)-\mcal{H}\left[\rho_V\widetilde{\Xi^{-1}}[h]\right]^{(m_1-1)}(x_1)\Big|
		\\\leq C C_5(V,n)C_7(V)\|h\|_{W_{m_1-1}^{\infty}(\R)}\leq C C_5(V,n)C_7(V)\cdot\|f\|_{W_{n+1}^{\infty}(\R^{p})}.
	\end{multline*}
	All the previous bounds yield for a constant $C>0$:
	\begin{equation*}\label{ineq:almostlinfinicontrolxi-1}
		\Big\|\widetilde{\Xi^{-1}}[h]^{(m_1)}\Big\|_{L^\infty([-M_V,M_V])}\leq C \left(1+\|\rho_{V}\|^{1/2}_{\infty}C_{\mcal{L}}+ C_5(V,n)\right) C_7(V)\cdot\|f\|_{W_{n+1}^{\infty}(\R^{p})}.
	\end{equation*}
	Upon taking the supremum over $x_2,\dots,x_p\in\R^{p-1}$ in $\eqref{ineq:almostlinfinicontrolxi-1}$, and over $m\leq n$, we obtain the conclusion for the constant $C_{W_n^{\infty}}(\widetilde{\Xi_1^{-1}})$ defined for a constant $C>0$ depending only on $n$,
	\begin{multline}\label{eq:C_9(V,l)}
		C_{W_n^{\infty}}(\widetilde{\Xi_1^{-1}})\defi C\Big[(1+\|\rho_V\|_{\infty}C_{\mcal{L}})\max_{m\leq n}\max_{\bm{i}}\|\mfrak{f}^{(5)}_{m,\bm{i}}\|_{L^\infty([-M_V,M_V]^c)}+C_6(V)		\\+\left(1+\|\rho_{V}\|^{1/2}_{\infty}C_{\mcal{L}}+C_5(V,n)\right) C_7(V)+ \|\mcal{I}_1\|_{\infty}+C_4(V)\Big].
	\end{multline} 
	The fact that, upon choosing the potential $V_{\phi,t}$ with $\phi\in H^{\infty}(\R)$, $t\mapsto	C_{W_n^{\infty}}(\widetilde{\Xi_1^{-1}})$ is continuous is shown in Proposition \ref{prop:constantescontinues}.
\end{proof}

Finally, we define the variable insertion operators which will also be involved in the loop equations.
\begin{definition}\label{def:repoperator}
	If $\phi$ is a function in $n$ variables, we define the $n$-th variable insertion operator $\Theta^{(i)}$ as
	\begin{equation}\label{def:opertheta}
		\Theta^{(i)}[\phi](\xi_1,\dots,\xi_{n-1})=\phi(\xi_1,\dots,\xi_{i-1},\xi_1,\xi_i,\dots,\xi_{n-1})
	\end{equation}
\end{definition}

We can deduce from Theorem \ref{thm: cont inverse master linfini}, the following corollary which will be necessary for the analysis of loop equations.

\begin{corollary}
	\label{thm: cont theta inverse master linfini}
	Let $n,p\geq1$, $ j\in\llbracket2,p+1\rrbracket$, for all $f\in W_{n+1}^{\infty}(\R^{p+1})$, $$\boxed{\left\|\Theta^{( j)}\circ\widetilde{\Xi^{-1}_1}[f]\right\|_{W^\infty_{n}(\R^{p})}\leq2 C_{W_n^{\infty}}(\widetilde{\Xi_1^{-1}})\cdot\|f\|_{W_{n+1}^\infty(\R^{p+1})}}$$
	where the constant 	$C_{W_n^{\infty}}(\widetilde{\Xi_1^{-1}})$ was introduced in Theorem \ref{thm: cont inverse master linfini}.
\end{corollary}

\begin{proof}
	Let $f\in W_{n+1}^{\infty}(\R^{p+1})$, $x_1,\dots,x_p\in\R$, $\underline{m}\defi (m_1,\dots,m_p)\in\N^p$ such that $m\defi \sum_{i=1}^{p}m_i\leq n$.
	\begin{multline*}
		\partial^{\underline{m}}\Theta^{( j)}\circ\widetilde{\Xi^{-1}_1}\left [f\right ](x_1,\dots,x_p)=\partial_1^{m_1}\partial_2^{m_2}\dots\partial_{ j-1}^{m_{ j-1}}\partial_{ j+1}^{m_ j}\dots\partial_{p+1}^{m_p}\widetilde{\Xi_1^{-1}}\left [f\right ](x_1,\dots,x_{ j-1},x_1,x_ j\dots,x_p)
		\\+\partial_2^{m_2}\dots\partial_{ j-1}^{m_{ j-1}}\partial_j^{m_1}\partial_{ j+1}^{m_ j}\partial_{p+1}^{m_p}\widetilde{\Xi_1^{-1}}\left [f\right ](x_1,\dots,x_{ j-1},x_1,x_ j\dots,x_p).
	\end{multline*}
	Thus $\left \|\partial^{\underline{m}}\Theta^{( j)}\circ\widetilde{\Xi^{-1}_1}\left [f\right ]\right \|_{L^\infty(\R^p)}\leq2\left \|\widetilde{\Xi^{-1}_1}[f]\right \|_{W_n^{\infty}(\R^{p+1})}\leq 2C_{W_n^{\infty}}(\widetilde{\Xi_1^{-1}})\cdot\|f\|_{W_{n+1}^{\infty}(\R^{p+1})}$.
\end{proof}

The last control that we need is on $\Theta^{(j)}\circ\widetilde{\Xi_1}^{-1}$  in $H^n$-norm. We recall that $\Theta^{(j)}$ is defined in Definition \ref{def:repoperator}.

\begin{theorem}[$H^{n}$-continuity for $\Theta^{(j)}\circ\widetilde{\Xi_1^{-1}}$]\label{thm:continversemaster thetaxi L2}
	Let $n, p\geq1$, there exists an explicit constant $ C_{H_n}(\Theta^{(2)}\circ\widetilde{\Xi_1^{-1}})>0$ (depending only on $n$ and $V$), such that for all $f\in H^{n+1}(\R^{p+1})$ and $j\in\llbracket2,p+1\rrbracket$, $$\boxed{\|\Theta^{(j)}\circ\widetilde{\Xi_1^{-1}}[f]\|_{H^{n}(\R^p)}\leq C_{H_n}(\Theta^{(2)}\circ\widetilde{\Xi_1^{-1}})\cdot\|f\|_{H^{n+1}(\R^{p+1})}.}$$
	Under the choice of potential $V_{\phi,t}$ defined in Theorem \ref{thm:conteqdensity}, for $\phi\in H^{\infty}(\R)$ the map $t\in[0,1]\mapsto C_{H_n}(\Theta^{(2)}\circ\widetilde{\Xi_1^{-1}})$ is continuous.
\end{theorem}

\begin{proof}Let $f\in H^{n+1}(\R^{p+1})$, $\underline{x}\in\R^{p}$, $\underline{m}\defi (m_1,\dots,m_p)\in\N^p$ such that $m\defi \sum_{i=1}^{p}m_i\leq n$. We set $h:(x,y)\mapsto \partial_2^{m_2}\dots\partial_{j-1}^{m_{j-1}}\partial_{j+1}^{m_ j}\dots\partial_{p+1}^{m_p}f(x,x_2,\dots,x_{ j-1},y,x_j,\dots,x_p)$ and $$g(x,y)\defi -h(x,y)+2P\mcal{H}\left[\rho_{V}\widetilde{\Xi^{-1}}\left [h(.,y)\right ]\right](x)+\mfrak{c}_{h}(y)$$ where we have set $\mfrak{c}_{h}(y)=\displaystyle\int_\R h(s,y)\diff\mu_V(s)-2P\int_\R\mcal{H}\left[\rho_V\widetilde{\Xi^{-1}}\left [h(.,y)\right ]\right](s)\diff\mu_V(s)$. Let $x\in\R$, 
	\begin{equation}\label{eq:leibnizderivethetaxi-1}
		\partial^{\underline{m}}\Theta^{( j)}\circ\widetilde{\Xi^{-1}_1}\left [f\right ](x_1,\dots,x_p)=\widetilde{\Xi^{-1}}[\partial_2^{m_1}h(.,x_1)](x_1)+\widetilde{\Xi^{-1}}[h(.,x_1)]^{(m_1)}(x_1).
	\end{equation}
	The idea of the proof is to obtain first a straightforward bound on the first term in \eqref{eq:leibnizderivethetaxi-1}. To bound the second term, we treat the cases $|x_1|>M_V$ and $|x_1|<M_V$ differently. For the first case, we rely on the expression \eqref{eq:2eme formule cont inver mast} obtained above for the $m$-\textit{th} derivative of $\widetilde{\Xi^{-1}}[h]$. We then apply the controls of  Lemma \ref{lem:gammaetalpha} and \ref{lem:majorationderivee nieme de rho Xi} to deduce the bound outside of the compact. On the compact $|x_1|< M_V$, we use again \eqref{eq:Xi-1deriveelinfini} and Lemma \ref{lem:majorationderivee nieme de rho Xi} to conclude on the last bound and therefore on Theorem \eqref{thm:continversemaster thetaxi L2}.

	\textbf{First term in \eqref{eq:leibnizderivethetaxi-1}.} The first term in the RHS is easy to control by Theorem \ref{thm: stru inv Xi}, 
	\begin{multline*}
		\widetilde{\Xi^{-1}}[\partial_2^{m_1}h(.,x_1)](x_1)=\dfrac{1}{\rho_V(x_1)}\int\limits_{x_1}^{\mrm{sgn}(x_1)\infty}\diff t\rho_V(t)\Bigg\{-\partial_2^{m_1}h(t,x_1)+\int_\R\partial_2^{m_1}h(s,x_1)\diff\mu_V(s)
		\\+2P\mcal{H}\left[\rho_V\widetilde{\Xi^{-1} }\left[\partial_2^{m_1}h(.,x_1) \right]\right](t)-2P\int_\R2P\mcal{H}\left[\rho_V\widetilde{\Xi^{-1} }\left[\partial_2^{m_1}h(.,x_1) \right]\right](s)\diff\mu_V(s) \Bigg\}.
	\end{multline*}
	The same argument as in the previous controls give that there exists a universal constant $C>0$ such that:
	\begin{multline}\label{eq:controle 1erterme}
		\int_\R\widetilde{\Xi^{-1}}[\partial_2^{m_1}h(.,x_1)](x_1)^2\diff x_1\\\leq C(1+\|\rho_V\|_{L^\infty(\R)}^2C_{\mcal{L}}^2)\Big[\|\mcal{I}_2 \|_{L^\infty(\R)}^2+\|\rho_V\|_{L^\infty(\R)} \|\mcal{I}_1 \|_{L^\infty(\R)}^2\Big]\cdot\|h\|_{H^{m_1+1}(\R^2)}^2
	\end{multline}
	where $\mcal{I}_1$ and $\mcal{I}_2$ have been defined in \eqref{eq:defI_a}.
	
	\textbf{Bound for $|x_1|>M_V$.} We now deal with the second term in \eqref{eq:leibnizderivethetaxi-1}. By \eqref{eq:2eme formule cont inver mast}, we have for a constant $C>0$ depending only on $m_1$:
	\begin{align}\label{eq:1ereq dans ThetaXi-1 H_n}
		\int_{M_V}^{+\infty}&\Xi^{-1}[h(.,x_1)]^{(m_1)}(x_1)^2\diff x_1\leq C\sum_{i_1=0}^{m_1-1}\sum_{i_2=0}^{m_1-i_1}\int_{M_V}^{+\infty}\diff x_1Q_{i_1}^{m_1}\left(\theta,\dots,\theta^{(i_1)}\right)(x_1) ^2\nonumber
		\\&\times\Bigg\{[\partial_1^{m_1-|\bm{i}|}g(x_1,x_1)]^2\alpha(x_1)^2P_{i_2}^{m_1-i_1}(\alpha,\dots,\alpha^{(i_2)})(x_1)^2\nonumber
		\\&+\dfrac{1}{\rho_V(x_1)^2}\left[\int\limits_{x_1}^{\mrm{sgn}(x_1)\infty}\diff t\rho_V(t)\partial_1^{m_1-|\bm{i}|+1}g(t,x_1)\alpha(t)P_{i_2}^{m_1-i_1}\left (\alpha,\dots,\alpha^{(i_2)}\right )(t)\right]^2\nonumber
		\\&+\dfrac{1}{\rho_V(x_1)^2}\left[\int\limits_{x_1}^{\mrm{sgn}(x_1)\infty}\diff t\rho_V(t)\partial_1^{m_1-|\bm{i}|}g(t,x_1)\left[\alpha P_{i_2}^{m_1-i_1}\left (\alpha,\dots,\alpha^{(i_2)}\right)  \right] '(t)\right]^2\Bigg\}.
	\end{align}
	Above, we used the notation $\bm{i}=(i_1,i_2)$ and $|\bm{i}|=i_1+i_2$. We first deal with the presence of $\mfrak{c}_h$ in the sum. This term only arises in the sum when $i_2=m_1-i_1$. By using the functions $\mfrak{f}_{m_1,\bm{i}}^{(k)}$ for $k\in\{1,2\}$ defined in \eqref{eq:mfrak f1,2}, we can bound these terms, for all $i_1\in\llbracket0,m_1-1\rrbracket$, by using a universal constant $C>0$:
	\begin{multline*}
		\max_{\bm{i}}\int_{M_V}^{+\infty}\diff x_1\mfrak{c}_h(x_1)^2\sum_{k=1}^{2}\mfrak{f}_{m_1,\bm{i}}^{(k)}(x_1)^2 \\\leq C\max_{k\in\{1,2\}}\|\mfrak{f}_{m_1,\bm{i}}^{(k)}\|_{L^\infty([-M_V,M_V]^c)}^2
		\int_\R \diff x_1\left(\int_\R h(t,x_1)^2\diff\mu_V(t)+\int_\R \mcal{H}\left[\rho_V\widetilde{\Xi^{-1}}\left[h(.,x_1) \right] \right] (t)^2\diff\mu_V(t)\right ) 
		\\\leq C\max_{\bm{i}}\max_{k\in\{1,2\}}\|\mfrak{f}_{m_1,\bm{i}}^{(k)}\|_{L^\infty([-M_V,M_V]^c)}^2 (\|\rho_V\|_{L^\infty(\R)}+\|\rho_V\|_{L^\infty(\R)}^3C_\mcal{L}^2)\cdot\|h\|^2_{H^1(\R^2)}.
	\end{multline*}
	Since we handled all the therms involving $\mfrak{c}_h$, it just remains to bound \eqref{eq:1ereq dans ThetaXi-1 H_n} with the substitution $g(x,y)\rightsquigarrow\mfrak{g}(x,y)\defi g(x,y)-\mfrak{c}_h(y)$, namely:
	\begin{multline*}
		\int_{M_V}^{+\infty}\diff x_1Q_{i_1}^{m_1}\left(\theta,\dots,\theta^{(i_1)}\right)(x_1) ^2\Bigg\{\partial_1^{m_1-\bm{i}|}\mfrak{g}(x_1,x_1)^2\alpha(x_1)^2P_{i_2}^{m_1-i_1}(\alpha,\dots,\alpha^{(i_2)})(x_1)^2
		\\+\dfrac{1}{\rho_V(x_1)^2}\left[\int\limits_{x_1}^{\mrm{sgn}(x_1)\infty}\diff t\rho_V(t)\partial_1^{m_1-|\bm{i}|+1}\mfrak{g}(t,x_1)\alpha(t)P_{i_2}^{m_1-i_1}\left (\alpha,\dots,\alpha^{(i_2)}\right )(t)\right]^2
		\\+\dfrac{1}{\rho_V(x_1)^2}\left[\int\limits_{x_1}^{\mrm{sgn}(x_1)\infty}\diff t\rho_V(t)\partial_1^{m_1-|\bm{i}|}\mfrak{g}(t,x_1)\left[\alpha P_{i_2}^{m_1-i_1}\left (\alpha,\dots,\alpha^{(i_2)}\right)  \right] '(t)\right]^2\Bigg\}.
	\end{multline*}
	For the first term, we use the fact that for fixed $x\in\R$, $t\mapsto\partial_1^{m_1-|\bm{i}|}\mfrak{g}(t,x)$ goes to zero at infinity as an element of $H^1(\R)$ and by Cauchy-Schwarz inequality, we get:
	$$|\partial_1^{m_1-|\bm{i}|}\mfrak{g}(x_1,x_1)|=\sqrt{\int^{x_1}_{+\infty}2\partial_1^{m_1-|\bm{i}|}\mfrak{g}(t,x_1)\partial_1^{m_1-|\bm{i}|+1} \mfrak{g}(t,x_1)\diff t}\leq\sqrt{2}\|\mfrak{g}(.,x_1)\|_{H^{m_1-|\bm{i}|+1}(\R)}.$$
	Furthermore, for all $|x_1|>M_V$ and all $\bm{i}$, $\mfrak{f}^{(1)}_{m_1,\bm{i}}$ is bounded since it is continuous and a $\underset{|x|\rightarrow\infty}{O}\left(V'(x)^{-2}\right)$ by Lemma \ref{lem:gammaetalpha}. We conclude, by Lemma \ref{lem:majorationderivee nieme de rho Xi} that, with $\mfrak{f}^{(1)}_{m_1,\bm{i}}$ being given in \eqref{eq:mfrak f1,2}, there exists a constant $C$ depending only on $n$ such that:
	\begin{align*}
		\sum_{i_1=0}^{m_1-1}\sum_{i_2=0}^{m_1-i_1}&\int_{M_V}^{+\infty}\diff x_1 ^2\partial_1^{m_1-|\bm{i}|}\mfrak{g}(x_1,x_1)^2\mfrak{f}^{(1)}_{m_1,\bm{i}}(x_1)^2
		\\&\leq C\max_{\bm{i}}\|\mfrak{f}^{(1)}_{m_1,\bm{i}}\|_{L^\infty([-M_V,M_V]^c)}^2\int_{M_V}^{+\infty}\|\mfrak{g}(.,x_1)\|_{H^{m_1-|\bm{i}|+1}(\R)}^2\diff x_1
		\\&\leq C\max_{\bm{i}}\|\mfrak{f}^{(1)}_{m_1,\bm{i}}\|_{L^\infty([-M_V,M_V]^c)}^2\|\mfrak{g}\|_{H^{m_1+1}(\R^2)}^2
		\\&\leq C\left[1+\|\rho_V\|_{L^\infty(\R)}C_3(V,n)\right] \max_{\bm{i}}\|\mfrak{f}^{(1)}_{m_1,\bm{i}}\|_{L^\infty([-M_V,M_V]^c)}^2\cdot\|\partial_3^{m_3}f(.,z)\|_{H^{m_1+1}(\R^2)}^2.
	\end{align*}
	It just remains to bound 
	\begin{multline*}
		\int\limits_{M_V}^{+\infty}\diff x_1\dfrac{Q_{i_1}^{m_1}\left(\theta,\dots,\theta^{(i_1)}\right)(x_1) ^2}{\rho_V(x_1)^2}\Bigg\{\left[\int\limits_{x_1}^{\mrm{sgn}(x_1)\infty}\diff t\rho_V(t)\partial_1^{m_1-|\bm{i}|+1}\mfrak{g}(t,x_1)\alpha(t)P_{i_2}^{m_1-i_1}\left (\alpha,\dots,\alpha^{(i_2)}\right )(t)\right]^2
		\\+\left[\int\limits_{x_1}^{\mrm{sgn}(x_1)\infty}\diff t\rho_V(t)\partial_1^{m_1-|\bm{i}|}\mfrak{g}(t,x_1)\left[\alpha P_{i_2}^{m_1-i_1}\left (\alpha,\dots,\alpha^{(i_2)}\right)  \right] '(t)\right]^2\Bigg\}.
	\end{multline*}
	For the first term, we use Cauchy-Schwarz inequality, Lemma \ref{lem:gammaetalpha} and the function $\mfrak{f}^{(3)}_{m_1,\bm{i}}$ defined in \eqref{eq:def mfrak{f}} to get for a universal constant $C>0$:
	\begin{align*}
		\dfrac{|Q_{i_1}^{m_1}\left(\theta,\dots,\theta^{(i_1)}\right)(x_1)|^2}{\rho_V(x_1)^2}\Bigg|&\int\limits_{x_1}^{\mrm{sgn}(x_1)\infty}\diff t\rho_V(t)\partial_1^{m_1-|\bm{i}|+1}\mfrak{g}(t,x_1)\alpha(t)P_{i_2}^{m_1-i_1}\left (\alpha,\dots,\alpha^{(i_2)}\right )(t)\Bigg|^2
		\\&\leq\max_{\bm{i}}\|\mfrak{f}^{(3)}_{m_1,\bm{i}}\|_{L^\infty([-M_V,M_V]^c)}^2\|\mfrak{g}(.,x_1)\|_{H^{n_1+1}(\R)}^2
		\\&\leq C\max_{\bm{i}}\|\mfrak{f}^{(3)}_{m_1,\bm{i}}\|_{L^\infty([-M_V,M_V]^c)}^2 \left(1+\|\rho_V\|_{L^\infty(\R)}C_3(V,n)\right )
		\\&\times\|\partial_2^{m_2}\dots\partial_{j-1}^{m_{j-1}}\partial_{j+1}^{m_j}\dots\partial_{p+1}^{m_p}f(.,x_2,\dots,x_{j-1},x_1,x_j,\dots,x_p)\|_{H^{m_1+1}(\R)}^2.
	\end{align*}
	We proceed in the exact same way for the second term and do the same thing on $]-\infty,-M_V]$, for every term we dealt with. Finally, by integrating with respect to $x_1$, collecting all the terms and then integrating over $x_1\in[-M_V,M_V]^c$ and over $x_2,\dots,x_p\in\R$, we get
	$$\left\|\Theta^{(j)}\Big[\widetilde{\Xi^{-1}}[\partial_2^{m_2}\dots\partial_{j-1}^{m_{j-1}}\partial_{j+1}^{m_j}\dots\partial_{p+1}^{m_p}f]^{(m_1)}\Big]\right \|_{L^2([-M_V,M_V]^c\times\R^{p-1})}\leq  C_{9}(V)\cdot\|f\|_{H^{n+1}(\R^{p+1})}$$
	with $ C_{9}(V)$ defined, for a constant $C>0$ depending only on $n$, by:
	\begin{multline}\label{eq:C_10(V,m_1)}
		C_{9}(V)^2\defi C(n)\max_{m\leq n }\max_{j\in\{1,3\}}\max_{\bm{i}}\|\mfrak{f}_{m,\bm{i}}^{(j)}\|_{L^\infty([-M_V,M_V]^c)}^2
		\\\times\left ( \|\rho_V\|_{L^\infty(\R)}+\|\rho_V\|_{L^\infty(\R)}^3C_\mcal{L}^2
		+1+\|\rho_V\|_{L^\infty(\R)}C_3(V,n)\right).
	\end{multline}
	\textbf{Bound for $|x_1|<M_V$.} It just remains to bound 
	$\displaystyle\int_{-M_V}^{M_V}\widetilde{\Xi^{-1}}[h(.,x_1)]^{(m_1)}(x_1)^2\diff x_1$. For that, we use \eqref{eq:Xi-1deriveelinfini}, to get for a constant $C>0$ only depending on $n$:
	\begin{multline*}
		\displaystyle\int\limits_{-M_V}^{M_V}\widetilde{\Xi^{-1}}[h(.,x_1)]^{(m_1)}(x_1)^2\diff x_1
		\\\leq C\displaystyle\int\limits_{-M_V}^{M_V}\diff x_1\Bigg\{\dfrac{S^{m_1}(\theta,\dots,\theta^{(m_1)})(x_1)^2}{\rho_V(x_1)^2}\left(\int\limits_{x_1}^{\mrm{sgn}(x_1)\infty}g(t,x_1)\rho_V(t)\diff t\right)^2 
		\\+\sum_{i=0}^{m_1-2}R_{m_1-i}^{m_1}\left(\theta,\dots,\theta^{(m_1-1-i)}\right)(x_1)^2[\partial_1^ig(x_1,x_1)]^2+[\partial_1^{m_1-1}g(x_1,x_1)]^2\Bigg\}.
	\end{multline*}
	By the same procedure as before, we first deal with $\mfrak{c}_h$ defined at the beginning of the proof, this yields, with $\mcal{I}_i$ defined in \eqref{eq:defI_a} and $C>0$ depending only on $n$:
	\begin{multline*}
		\int\limits_{-M_V}^{M_V}\mfrak{c}_h(x_1)^2\diff x_1\Bigg[\dfrac{S^{m_1}(\theta,\dots,\theta^{(m_1)})(x_1)^2}{\rho_V(x_1)^2}\left(\int\limits_{x_1}^{\mrm{sgn}(x_1)\infty}\rho_V(t)\diff t\right)^2+R_{m_1}^{m_1}\left(\theta,\dots,\theta^{(m_1-1-i)}\right)(x_1)^2 \Bigg]
		\\\leq C\|h\|_{H^1(\R^2)}^2\|\rho_V\|_{L^\infty(\R)}\cdot\left(1+\|\rho_V\|_{L^\infty(\R)}^2C_\mcal{L} ^2\right)\cdot\Big\{ \|\mcal{I}_1\|_{\infty}^2\Big\|S^{m_1}(\theta,\dots,\theta^{(m_1)})\Big\|_{L^\infty([-M_V,M_V])}^2
		\\+\Big\|R_{m_1}^{m_1}\left(\theta,\dots,\theta^{(m_1-1)}\right)\Big\|_{L^\infty([-M_V,M_V])}^2\Big\}. 
	\end{multline*}
	Hence as before, we can replace $g$ by $\mfrak{g}$ and conclude with the last bounds:
	\begin{multline*}\displaystyle\int\limits_{-M_V}^{M_V}\dfrac{S^{m_1}(\theta,\dots,\theta^{(m_1)})(x_1)^2}{\rho_V(x_1)^2}\left(\int\limits_{x_1}^{\mrm{sgn}(x_1)\infty}\mfrak{g}(t,x_1)\rho_V(t)\diff t\right)^2\diff x_1
		\\\leq \Big\|\mcal{I}_2S^{m_1}(\theta,\dots,\theta^{(m_1)})\Big\|_{L^\infty([-M_V,M_V])}^2\cdot\|\mfrak{g}\|_{L^2(\R^2)}^2,
	\end{multline*}
	and by Cauchy-Schwarz inequality, with a constant $C>0$ only depending on $n$:
	\begin{multline*}
		\int\limits_{-M_V}^{M_V}\diff x_1\Big(\sum_{i=0}^{m_1-2}R_{m_1-i}^{m_1}\left(\theta,\dots,\theta^{(m_1-1-i)}\right)(x_1)^2\partial_1^i\mfrak{g}(x_1,x_1)^2+\partial_1^{m_1-1}\mfrak{g}(x_1,x_1)^2\Big)
		\\\leq\Bigg\{1+\max_{i\in\llbracket0,m_1-2\rrbracket}\Big\| R_{m_1-i}^{m_1}\left(\theta,\dots,\theta^{(m_1-1-i)}\right)\Big\|_{L^\infty([-M_V,M_V])}^2\Bigg\}
		\\\times\displaystyle\sum_{i=0}^{m_1-1}\int\limits_{-M_V}^{M_V}\diff x_1\int\limits_{+\infty}^{x_1}\diff t\partial_1^{i+1}\mfrak{g}(t,x_1)\partial_1^{i}\mfrak{g}(t,x_1)
		\\\leq C\Bigg\{1+\max_{i\in\llbracket0,m_1-2\rrbracket}\Big\| R_{m_1-i}^{m_1}\left(\theta,\dots,\theta^{(m_1-1-i)}\right)\Big\|_{L^\infty([-M_V,M_V])}^2\Bigg\}\cdot\|\mfrak{g}\|_{H^{m_1}(\R^2)}^2.
	\end{multline*}
	Moreover, by Lemma \ref{lem:majorationderivee nieme de rho Xi}, we have for $C>0$ a universal constant:
	$$\|\mfrak{g}\|_{H^{m_1}(\R^2)}\leq\|h\|_{H^{m_1}(\R^2)}+2P\pi\left\|\rho_V\widetilde{\Xi^{-1}_1}[h]\right\|_{H^{m_1}(\R^2)}\leq C \left[1+C_3(V,n)\right]\cdot\|h\|_{H^{m_1}(\R^2)}.$$
	We can then conclude that, for $C>0$ depending only on $n$:
	\begin{multline*}
		\int\limits_{-M_V}^{M_V}\diff x_1\Big(\sum_{i=0}^{m_1-2}R_{m_1-i}^{m_1}\left(\theta,\dots,\theta^{(m_1-1-i)}\right)(x_1)^2\partial_1^i\mfrak{g}(x_1,x_1)^2+\partial_1^{m_1-1}\mfrak{g}(x_1,x_1)^2\Big)
		\\\leq C\left[1+ C_3(V,n)\right]\Bigg\{1+\max_{i\in\llbracket0,m_1-2\rrbracket}\Big\| R_{m_1-i}^{m_1}\left(\theta,\dots,\theta^{(m_1-1-i)}\right)\Big\|_{L^\infty([-M_V,M_V])}^2\Bigg\}\cdot
		\|h\|_{H^{m_1}(\R^2)}^2.
	\end{multline*}
	Thus by integrating with respect to $z\in\R$, we get, 
	$$\left \|\Theta^{(j)}\Big[\widetilde{\Xi^{-1}}[\partial^{m_2}_2\dots\partial_{j-1}^{m_{j-1}}\partial_{j+1}^{m_j}\dots\partial_{p+1}^{m_p}f]^{(m_1)}\Big]\right \|_{L^2([-M_V,M_V]\times\R^{p-1})}\leq  C_{10}(V,n)\cdot\|\partial_3^{m_3}f\|_{H^{m_1+1}(\R^{p+1})}$$
	with $ C_{10}(V)>0$ defined, for a constant $C>0$ depending only on $n$, by:
	\begin{multline}\label{eq:C_11(V,n)}
		C_{10}(V)^2\defi C\max_{m\leq n}\Bigg\{\|\rho_V\|_{L^\infty(\R)}\left(1+\|\rho_V\|_{L^\infty(\R)}^2C_\mcal{L} ^2\right)
		\\\times\Big\{ \left[1+ C_3(V,n)\right]^2\Bigg(\max_{i\in\{1,2\}}\|\mcal{I}_i\|_{\infty}^2\Big\|S^{m}(\theta,\dots,\theta^{(m)})\Big\|_{L^\infty([-M_V,M_V])}^2
		\\+1+\max_{i\in\llbracket0,m-2\rrbracket}\Big\| R_{m-i}^{m}\left(\theta,\dots,\theta^{(m-1-i)}\right)\Big\|_{L^\infty([-M_V,M_V])}^2\Bigg)  \Bigg\}.
	\end{multline} Collecting the bounds on the $L^2$-norms of $\Theta^{(j)}\Big[\widetilde{\Xi^{-1}}[\partial^{m_2}_2\dots\partial_{j-1}^{m_{j-1}}\partial_{j+1}^{m_j}\dots\partial_{p+1}^{m_p}f]^{(m_1)}\Big]$ on $[-M_V,M_V]^c\times\R$ and $[-M_V,M_V]\times\R$, we obtain:
	$$\left \|\Theta^{(j)}\Big[\widetilde{\Xi^{-1}}[\partial^{m_2}_2\dots\partial_{j-1}^{m_{j-1}}\partial_{j+1}^{m_j}\dots\partial_{p+1}^{m_p}f]^{(m_1)}\Big]\right \|_{L^2(\R^p)}\leq 2\max_{i\in\{9,10\}}C_{i}(V)\cdot\|h\|_{H^{m_1+1}(\R^{p+1})}.$$
	By combining the bound above together with \eqref{eq:controle 1erterme} and taking the supremum over $m\leq n$, we get the desired conclusion
	with, $ C_{9}$ and $ C_{10}$ being given in \eqref{eq:C_10(V,m_1)} and \eqref{eq:C_11(V,n)}, $\mcal{I}_i$ being given in \eqref{eq:defI_a};
	\begin{multline}\label{eq:C12}
		C_{H_n}(\Theta^{(2)}\circ\widetilde{\Xi_1^{-1}})\defi C\max_{i\in\{9,10\}}C_{i}(V)
		\\+(1+\|\rho_V\|_{L^\infty(\R)}^2C_{\mcal{L}}^2)\left (1+\|\rho_V\|_{L^\infty(\R)}\right ) \max_{i\in\{1,2\}}\|\mcal{I}_i \|_{L^\infty(\R)}^2.
	\end{multline}
	Noticing that this constant does not depend on $j$ yields the conclusion. The fact that, upon choosing the potential $V_{\phi,t}$ with $\phi\in H^{\infty}(\R)$, $t\mapsto C_{H_n}(\Theta^{(2)}\circ\widetilde{\Xi_1^{-1}})$ is continuous is shown in Proposition \ref{prop:constantescontinues}.
\end{proof}
\section{Asymptotic expansion of the linear statistics}\label{section:asymptot correlators}
In this section, we prove the loop equations for general functions, using the continuity results proven in Sections \ref{section:loopequations} and \ref{section4}. We then prove Theorem \ref{thm:correlators} by using the \textit{a priori} bound proven in Proposition \ref{a priori bound} and again the continuity results proven in Sections \ref{section:loopequations} and \ref{section4}.

\subsection{Loop equations for general functions}
We now have all the needed ingredients to state the loop equations, we recall the definition of a linear statistic was defined in \eqref{eq:braket}. We follow the approach of  \cite[Prop 3.2.3]{borot2016asymptotic} and adapt it to our setting. Since the equilibium measure and thus the master operator are different in the high-temperature regime (compared to the fixed-temperature regime), the standard arguments of recentring and extension to general functions by density must be adapted.

\begin{theorem}[Loop equations]\label{thm:DSequations}
	The level $1$ loop equation holds for all $\psi_1\in H^3(\R)$ and takes the form:
	\begin{equation}\label{DSlvl1}
		\Braket{\psi_1}_{\Delta\mu_N}=
		\dfrac{P}{N}\Braket{\widetilde{\Xi^{-1}}[\psi_1]'}_{\mu_V}+\dfrac{P}{N}\Braket{\widetilde{\Xi^{-1}}[\psi_1]'}_{\Delta\mu_N}-P\Braket{\mcal{D}\circ\widetilde{\Xi^{-1}}[\psi_1]}_{\Delta\mu_N\otimes \Delta\mu_N}.
	\end{equation}
	For all $\psi_n\in H^{3}(\R^{n})$, the level $n>1$ loop equations reads:
	\begin{multline}\label{DSlvlnthm}
		\Braket{\psi_{n}}_{\overset{n}{\bigotimes}\Delta\mu_N}=\dfrac{P}{N}\Braket{\partial_1\widetilde{\Xi_1^{-1}}[\psi_{n}]}_{\mu_V\overset{n-1}{\bigotimes}\Delta\mu_N}+\dfrac{P}{N}\Braket{\partial_1\widetilde{\Xi_1^{-1}}[\psi_{n}]}_{\overset{n}{\bigotimes}\Delta\mu_N}-P\Braket{\mcal{D}_1\circ\widetilde{\Xi_1^{-1}}[\psi_{n}]}_{\overset{n+1}{\bigotimes}\Delta\mu_N}
		\\-\dfrac{1}{N}\sum_{i=2}^{n}\Bigg(\Braket{\Theta^{(i)} \circ\widetilde{\Xi_1^{-1}}\big[\partial_i\psi_{n}\big] }_{\overset{n-1}{\bigotimes}\Delta\mu_N}+\Braket{\Theta^{(i)} \circ\widetilde{\Xi_1^{-1}}\big[\partial_i\psi_{n}\big] }_{\mu_V\overset{n-2}{\bigotimes}\Delta\mu_N}\Bigg).
	\end{multline}
\end{theorem}

\begin{proof}Let $\big(\psi^{(i)}\big)_{i=1}^{n+1}\in H^3(\R)^{n+1}$ and set $\phi^{(1)}\defi\widetilde{\Xi^{-1}}[\psi^{(1)}]$
	, set $V_{\underline{\varepsilon}}(\lambda)\defi V(\lambda)+\displaystyle\sum_{i=2}^{n+1}\epsilon_i\psi^{(i)}(\lambda_i)$ and then define
	\begin{equation*}
		p_N^{(\underline{\varepsilon})}(\underline{\lambda})\defi \dfrac{1}{\mcal{Z}_N[V_{\underline{\varepsilon}}]}\cdot\prod_{i<j}^N|\lambda_i-\lambda_j|^{\frac{2P}{N}}\cdot\prod_{i=1}^{N}\cdot e^{-V_{\underline{\varepsilon}}(\lambda_i)}
	\end{equation*}
	We then define $G_t(\mu)=\mu+t\phi^{(1)}(\mu)$. We claim that $\partial_1\phi^{(1)}$ is bounded. Indeed because $\psi^{(1)}\in H^3(\R)$ and Theorem \ref{thm: cont inverse master}, $\phi^{(1)}\in H^{2}(\R)$ so in particular $\partial_1\phi^{(1)}$ is bounded.
	
	As a consequence, for $t$ small enough $G_t$ is a diffeomorphism over $\R$. By a change of variable $\lambda_i=G_t(\mu_i)$, we obtain :
	\begin{equation*}
		1=\int_{\R^N}p_N^{(\underline{\varepsilon})}(\lambda_1,\dots,\lambda_N)\diff^N\underline{\lambda}=\int_{\R^N}p_N^{(\underline{\varepsilon})}\Big(G_t(\mu_1,\dots,G_t(\mu_N)\Big)\cdot\prod_{i=1}^NG_t'(\mu_i)\diff^N\underline{\mu}.
	\end{equation*}
	Making an asymptotic expansion up to the first order in $t$ of the right-hand side of this equality yields:
	\begin{multline*}
		1=\int_{\R^N}\diff^N\lambda\cdot p_N^{(\underline{\varepsilon})}(\underline{\lambda})\cdot\Bigg\{1+t\sum_{i=1}^{N}\partial_1\phi^{(1)}(\lambda_i)\Bigg\}\cdot\Bigg\{1+t\dfrac{P}{N}\sum_{i\neq j}\dfrac{\phi^{(1)}(\lambda_i)-\phi^{(1)}(\lambda_j)}{\lambda_i-\lambda_j}\Bigg\}\\\cdot\Bigg\{1-t\sum_{i=1}^NV_{\underline{\varepsilon}}'(\lambda_i)\phi^{(1)}(\lambda_i)\Bigg\}+O(t^2).
	\end{multline*}
	Identifying the terms linear in $t$ leads to:
	\begin{equation*}
		0=-N\Braket{V_{\underline{\varepsilon}}'\phi^{(1)}}_{\mu_N}^{(\underline{\varepsilon})}+PN\Braket{\mcal{D}[\phi^{(1)}]}_{\mu_N\otimes \mu_N}^{(\underline{\varepsilon})}+(N-P)\Braket{\partial_1\phi^{(1)}}_{\mu_N}^{(\underline{\varepsilon})}
	\end{equation*}
	where we used, for any (signed) measure $\nu$, the following notation
	$$\Braket{f}_{\nu}^{(\underline{\varepsilon})}\defi\int_{\R^N}\left(\int_{\R}f(x)\diff\nu(x)\right)\cdot	p_N^{(\underline{\varepsilon})}(\underline{\lambda})\cdot\diff^N\underline{\lambda}.$$
	By definition of $V_{\underline{\varepsilon}}$, we get:
	\begin{equation*}
		\Braket{V'\phi^{(1)}}_{\mu_N}^{(\underline{\varepsilon})}+\sum_{i=2}^{n+1}\varepsilon_i\Braket{\phi^{(1)}(\psi^{(i)})'}_{\mu_N}^{(\underline{\varepsilon})}-P\Braket{\mcal{D}[\phi^{(1)}]}_{\mu_N\otimes \mu_N}^{(\underline{\varepsilon})}-\left (1-\dfrac{P}{N}\right )\Braket{(\phi^{(1)})'}_{\mu_N}^{(\underline{\varepsilon})}=0.
	\end{equation*}
	It becomes, after recentring the empirical measures against $\mu_V$ we obtain:
	\begin{multline*}
		\Braket{V'\phi^{(1)}}_{\Delta\mu_N}^{(\underline{\varepsilon})}+\Braket{V'\phi^{(1)}}_{\mu_V}-P\Braket{\mcal{D} [\phi^{(1)} ]}_{\Delta\mu_N\otimes \Delta\mu_N}^{(\underline{\varepsilon})}-2P\Braket{\mcal{D} [\phi^{(1)} ]}_{\Delta\mu_N\otimes \mu_V}^{(\underline{\varepsilon})}-P\Braket{\mcal{D} [\phi^{(1)} ]}_{\mu_V\otimes \mu_V}
		\\-\Big(1-\dfrac{P}{N}\Big)\left (\Braket{(\phi^{(1)})'}_{\Delta\mu_N}^{(\underline{\varepsilon})}+\Braket{(\phi^{(1)})'}_{\mu_V}\right)+\sum_{i=2}^{n+1}\varepsilon_i\Big(\Braket{\phi^{(1)}(\psi^{(i)})'}_{\Delta\mu_N}^{(\underline{\varepsilon})}+\Braket{\phi^{(1)}(\psi^{(i)})'}_{\mu_V}\Big)=0.
	\end{multline*}
	Using the following identities true for all signed measure $\nu$, 
	$$\Braket{\mcal{D} [\phi^{(1)} ]}_{\nu\otimes \mu_V}^{(\underline{\varepsilon})}=-\Braket{\phi^{(1)}\mcal{H}[\rho_V]}_{\nu}^{(\underline{\varepsilon})}+\Braket{\mcal{H}[\phi^{(1)}\rho_V]}_{\nu}^{(\underline{\varepsilon})}; \hspace{1cm}\Braket{\dfrac{\rho_V'}{\rho_V}\phi^{(1)}}_{\mu_V}=-\Braket{(\phi^{(1)})'}_{\mu_V},$$and $$V'+2P\mcal{H}[\rho_V]=-\dfrac{\rho_V'}{\rho_V}, \hspace{1,5cm}\Braket{\mcal{H}[\phi^{(1)}\rho_V]}_{\mu_V}=-\Braket{\phi^{(1)}\mcal{H}[\rho_V]}_{\mu_V},$$
	we get:
	\begin{multline*}
		-\Braket{\dfrac{\rho_V'}{\rho_V}\phi^{(1)}}_{\Delta\mu_N}^{(\underline{\varepsilon})}-\Braket{(\phi^{(1)})'}_{\Delta\mu_N}^{(\underline{\varepsilon})}-2P\Braket{\mcal{H}[\phi^{(1)}\rho_V]}_{\Delta\mu_N}^{(\underline{\varepsilon})}+\dfrac{P}{N}\Braket{(\phi^{(1)})'}_{\Delta\mu_N}^{(\underline{\varepsilon})}-P\Braket{\mcal{D}[\phi^{(1)}]}_{\Delta\mu_N\otimes \Delta\mu_N}^{(\underline{\varepsilon})}
		\\+\dfrac{P}{N} \Braket{(\phi^{(1)})'}_{\mu_V}+\sum_{i=2}^{n+1}\varepsilon_i\Big(\Braket{\phi^{(1)}(\psi^{(i)})'}_{\Delta\mu_N}^{(\underline{\varepsilon})}+\Braket{\phi^{(1)}(\psi^{(i)})'}_{\mu_V}\Big)=0.
	\end{multline*}
	which leads to:
	\begin{multline}\label{DSbuilding}
		\Braket{\Xi[\phi^{(1)}]}_{\Delta\mu_N}^{(\underline{\varepsilon})}=\dfrac{P}{N} \Braket{(\phi^{(1)})'}_{\mu_V}+\dfrac{P}{N}\Braket{(\phi^{(1)})'}_{\Delta\mu_N}^{(\underline{\varepsilon})}-P\Braket{\mcal{D} [\phi^{(1)} ]}_{\Delta\mu_N\otimes \Delta\mu_N}^{(\underline{\varepsilon})}
		\\+\sum_{i=2}^{n+1}\varepsilon_i\Big(\Braket{\phi^{(1)}(\psi^{(i)})'}_{\Delta\mu_N}^{(\underline{\varepsilon})}+\Braket{\phi^{(1)}(\psi^{(i)})'}_{\mu_V}\Big)
	\end{multline}
	Taking $\varepsilon_i$ go to 0 for all $i$, $\phi^{(1)}=\Xi^{-1}[\psi]$ leads to:
	\begin{equation*}
		\Braket{\Xi[\phi^{(1)}]}_{\Delta\mu_N}=\dfrac{P}{N} \Braket{(\phi^{(1)})'}_{\mu_V}+\dfrac{P}{N}\Braket{(\phi^{(1)})'}_{\Delta\mu_N}-P\Braket{\mcal{D} [\phi^{(1)} ]}_{\Delta\mu_N\otimes \Delta\mu_N}.
	\end{equation*}
	And by definition of $\phi^{(1)}$ and the fact that $\Delta\mu_N$ is a measure zero mass, we obtain:
	$$	\Braket{\psi^{(1)}}_{\Delta\mu_N}=
	\dfrac{P}{N}\Braket{\widetilde{\Xi^{-1}}[\psi^{(1)}]'}_{\mu_V}+\dfrac{P}{N}\Braket{\widetilde{\Xi^{-1}}[\psi^{(1)}]'}_{\Delta\mu_N}-P\Braket{\mcal{D}\circ\widetilde{\Xi^{-1}}[\psi^{(1)}]}_{\Delta\mu_N\otimes \Delta\mu_N}.$$
	Hence we obtain \eqref{DSlvl1}.
	
	Furthermore, in \eqref{DSbuilding} multiplying by $\dfrac{\mcal{Z}_N[V_{\underline{\varepsilon}}]}{\mcal{Z}_N[V]}$, applying $\partial_{\varepsilon_2}\dots\partial_{\varepsilon_{n+1}}$ and evaluate at $\varepsilon_i=0$ leads to
	\begin{multline}\label{DSbuildinglvl}
		(-N)^{n}\Braket{\Xi[\phi^{(1)}] (\xi_1)\cdot\prod_{i=2}^{n+1}\psi^{(i)}(\xi_i)}_{\overset{n+1}{\bigotimes}\Delta\mu_N}=(-N)^{n}\dfrac{P}{N}\Braket{(\phi^{(1)})'(\xi_1)\cdot\prod_{i=2}^{n+1}\psi^{(i)}(\xi_i)}_{\mu_V\overset{n}{\bigotimes}\Delta\mu_N}
		\\+(-N)^{n}\dfrac{P}{N}\Braket{(\phi^{(1)})'(\xi_1)\cdot\prod_{i=2}^{n+1}\psi^{(i)}(\xi_i)}_{\overset{n+1}{\bigotimes}\Delta\mu_N}-P(-N)^{n}\Braket{\mcal{D}[\phi^{(1)}](\xi_1,\xi_2)\prod_{i=2}^{n+1}\psi^{(i)}(\xi_{i+1})}_{\overset{n+2}{\bigotimes}\Delta\mu_N}
		\\+(-N)^{n-1}\sum_{i=2}^{n+1}\Bigg(\Braket{\phi^{(1)}(\xi_1)(\psi^{(i)})'(\xi_1)\prod_{\substack{j=2\\j\neq i}}^{n+1}\psi^{(j)}(\xi_j)}_{\overset{n}{\bigotimes}\Delta\mu_N}\hspace{-1cm}+\Braket{\phi^{(1)}(\xi_1)(\psi^{(i)})'(\xi_1)\prod_{\substack{j=2\\j\neq i}}^{n+1}\psi^{(j)}(\xi_j)}_{\mu_V\overset{n-1}{\bigotimes}\Delta\mu_N}\hspace{-1cm}\Bigg)
	\end{multline}
	Using again that $\phi^{(1)}=\widetilde{\Xi^{-1}}[\psi^{(1)}]$ and that $\Delta\mu_N$ has zero mass, dividing by $(-N)^n$ and defining $\psi_{n+1}(\xi_1,\dots,\xi_{n+1})\defi \prod_{i=1}^{n+1}\psi^{(i)}(\xi_i)$ and using the operators defined earlier leads to
	\begin{multline}\label{DSbuiltl}
		\Braket{\psi_{n+1}}_{\overset{n+1}{\bigotimes}\Delta\mu_N}=\dfrac{P}{N}\Braket{\partial_1\widetilde{\Xi_1^{-1}}[\psi_{n+1}]}_{\mu_V\overset{n}{\bigotimes}\Delta\mu_N}+\dfrac{P}{N}\Braket{\partial_1\widetilde{\Xi_1^{-1}}[\psi_{n+1}]}_{\overset{n+1}{\bigotimes}\Delta\mu_N}-P\Braket{\mcal{D}_1\circ\widetilde{\Xi_1^{-1}}[\psi_{n+1}]}_{\overset{n+2}{\bigotimes}\Delta\mu_N}
		\\-\dfrac{1}{N}\sum_{i=2}^{n+1}\Bigg(\Braket{\Theta^{(i)} \circ\widetilde{\Xi_1^{-1}}\big[\partial_i\psi_{n+1}\big] }_{\overset{n}{\bigotimes}\Delta\mu_N}+\Braket{\Theta^{(i)} \circ\widetilde{\Xi_1^{-1}}\big[\partial_i\psi_{n+1}\big] }_{\mu_V\overset{n-1}{\bigotimes}\Delta\mu_N}\Bigg).
	\end{multline}
	It then not hard to see that, since $\bigotimes^{n+1}H^{3}(\R)$ is dense in $H^{3}(\R^{n+1})$ for $\|.\|_{H^{3}(\R^{n+1})}$, because of Theorem \ref{thm:fini diff ope reg}, \ref{thm: cont inverse master} and \ref{thm:continversemaster thetaxi L2}, \eqref{DSbuiltl} remains valid for all $\psi_{n+1}\in H^{3}(\R^{n+1})$. This establishes \eqref{DSlvlnthm}.
\end{proof}

\subsection{Asymptotic expansion of linear statistics}
The a priori bound on linear statistics of Proposition \ref{a priori bound}, provides the main ingredient for obtaining the existence of their large-$N$ asymptotic expansion in powers of $N^{-1}$ up to any order using the loop equations Theorem \ref{thm:DSequations}.
\begin{theorem}\label{asymptotic correlators}
	Let $\psi_k\in \displaystyle H^\infty(\R^k)$, then for all integer $K$, there exists a sequence $(d_{i}^{(k)})_{i\geq\lceil k/2\rceil}\in\R^N$ (depending on $V$) such that
	$$\boxed{\Braket{\psi_k}_{\overset{k}{\bigotimes}\Delta\mu_N}=\sum_{i=\lceil k/2\rceil}^{K}\dfrac{d_{i}^{(k)}(\psi_k)}{N^i}+O\left (N^{-(K+1)}\right )}$$
	with $$d_{1}^{(1)}[\psi_1]=P\Braket{\widetilde{\Xi^{-1}}[\psi_1]'}_{\mu_V}+P\Braket{\Theta^{(2)} \circ\widetilde{\Xi_1^{-1}}\left[\partial_2\mcal{D}\circ\widetilde{\Xi^{-1}}[\psi_1]\right]  }_{\mu_V}.$$
	Furthermore, there exists a sequence of integers $\left(m_{K}\right) >0$ (depending on $K$ and $k$), increasing in $K$, such that for all $k\geq1$ and $K\geq0$, all $\psi_k\in H^\infty(\R^k)$,
	\begin{equation}\label{eq:controleraminedercorrelator}
		\Bigg|\Braket{\psi_k}_{\overset{k}{\bigotimes}\Delta\mu_N}-\sum_{i=\lceil k/2\rceil}^{K}\dfrac{d_{i}^{(k)}(\psi_k)}{N^i}\Bigg|\leq \dfrac{C_{\mrm{rem}}(V)}{N^{K+1}}\mcal{N}_{m_{K}}(\psi_k).
	\end{equation}
	Above $\mcal{N}_{m}(\psi_k)\defi \max\left(\|\psi_k\|_{W_m^{\infty}(\R^k)},\|\psi_k\|_{H^{m}(\R^k)}\right)$, while $ C_{\mrm{rem}}(V)>0$ is a constant (depending on $V$, $K$ and $k$). Finally, under the choice of potential $V_{\phi,t}$ defined in Theorem \ref{thm:conteqdensity}, for $\phi\in H^{\infty}(\R)$ the map $t\in[0,1]\mapsto C_{\mrm{rem}}(V_{\phi,t})$ is continuous.\end{theorem}
\begin{proof}Using the first loop equation given in Theorem \ref{thm:DSequations}, we get:
	\begin{equation}\label{premiere expansion}
		\Braket{\psi_1}_{\Delta\mu_N}=
		\dfrac{P}{N}\Braket{\widetilde{\Xi^{-1}}[\psi_1]'}_{\mu_V}+\dfrac{P}{N}\Braket{\widetilde{\Xi^{-1}}[\psi_1]'}_{\Delta\mu_N}-P\Braket{\mcal{D}\circ\widetilde{\Xi^{-1}}[\psi_1]}_{\Delta\mu_N\otimes \Delta\mu_N}.
	\end{equation}
	where we recall that $\widetilde{\Xi^{-1}}=\Xi^{-1}\circ\mcal{X}$ defined in Definition \ref{def:recentring}.
	The idea is to verify the hypotheses of Theorems \ref{a priori bound} for each function involved in the Dyson-Schwinger equations. By Proposition \ref{prop:diff fini linfini} and Theorem \ref{thm:fini diff ope reg}, \ref{thm: cont inverse master},  \ref{thm: cont inverse master linfini}, \ref{thm: cont theta inverse master linfini} and \ref{thm:continversemaster thetaxi L2} and the fact that $\psi_k\in H^\infty(\R^k)$, we're ensured that all the norms are finite and that a $n$-linear statistic will be a $O\left (N^{-\frac{n}{2}(1-\varepsilon)}\right )$ where $\varepsilon>0$ is fixed but can be chosen arbitrarly small.
	
	We show by induction on $K$ that there exists an asymptotic expansion up to $o(N^{-K})$ for any function $\psi_k\in H^\infty(\R^k)$ for all $k\leq 2K$.
	
	For $K=1$, since the first term in the RHS of \eqref{premiere expansion} clearly contributes to the asymptotic expansion of $	\Braket{\psi_1}_{\Delta\mu_N}$ up to $o(N^{-1})$, we focus on the two other terms. In \eqref{premiere expansion}, the second term is clearly a $o(N^{-1})$ since by Proposition \ref{a priori bound}, Theorem \ref{thm: cont inverse master} and Theorem \ref{thm: cont inverse master linfini}, there exists $C>0$ (depending only on $\varepsilon>0$)
	\begin{multline*}
		|\Braket{\partial_1\widetilde{\Xi^{-1}}[\psi_1]}_{\Delta\mu_N}|\leq Ce^{K_V}N^{-(1-\varepsilon)/2}\left(\|\widetilde{\Xi^{-1}}[\psi_1]'\|_{H^{1/2}(\R)}+\|\widetilde{\Xi^{-1}}[\psi_1]'\|_{W^{\infty}_1(\R)} \right)
		\\\leq Ce^{K_V}N^{-(1-\varepsilon)/2}\left(C_{H^2}(\widetilde{\Xi_1^{-1}})+C_{W_2^\infty}(\widetilde{\Xi_1^{-1}})\right) \mcal{N}_3(\psi_1).
	\end{multline*}
	To obtain the expansion of the $2$-linear statistic up to $o(N^{-1})$, we will need to use the loop equation at level 2 with $\psi_2\defi \mcal{D}\circ\widetilde{\Xi^{-1}}[\psi]$. Let $\psi_2\in H^\infty(\R^2)$ be arbitrary for now. The level $2$ equation reads:
	\begin{multline}\label{deuxieme expansion}
		\Braket{\psi_{2}}_{\overset{2}{\bigotimes}\Delta\mu_N}=\dfrac{P}{N}\Braket{\partial_1\widetilde{\Xi_1^{-1}}[\psi_{2}]}_{\mu_V\bigotimes\Delta\mu_N}+\dfrac{P}{N}\Braket{\partial_1\widetilde{\Xi_1^{-1}}[\psi_{2}]}_{\overset{2}{\bigotimes}\Delta\mu_N}-\dfrac{1}{N}\Braket{\Theta^{(2)} \circ\widetilde{\Xi_1^{-1}}\big[\partial_2\psi_{2}\big] }_{\Delta\mu_N}
		\\-P\Braket{\mcal{D}_1\circ\widetilde{\Xi_1^{-1}}[\psi_{2}]}_{\overset{3}{\bigotimes}\Delta\mu_N}-\dfrac{1}{N}\Braket{\Theta^{(2)} \circ\widetilde{\Xi_1^{-1}}\big[\partial_2\psi_{2}\big]}_{\mu_V}.
	\end{multline}
	The first term is a $o(N^{-1})$ as a 1-linear statistic $\Braket{\psi}_{\Delta\mu_N}$ where $\psi(x)\defi \displaystyle\int_\R\partial_1\widetilde{\Xi_1^{-1}}[\psi_{2}](x,y)\diff\mu_V(y)$. This function is indeed in $H^1(\R)$ because Theorem  \ref{thm: cont inverse master} gives:
	$$\|\psi\|_{H^1(\R)}^2\leq\|\rho_V\|_{\infty}\|\partial_1\widetilde{\Xi_1^{-1}[}\psi_{2}]\|_{H^1(\R^2)}^2\leq\|\rho_V\|_{\infty}C_{H^2}(\widetilde{\Xi_1^{-1}}) ^2\cdot\|\psi_2\|_{H^3(\R^2)}^2$$
	and in $W_1^{\infty}(\R)$ by Theorem \ref{thm: cont inverse master linfini}
	$$\|\psi\|_{W_1^{\infty}(\R)}\le\|\partial_1\widetilde{\Xi_1^{-1}[}\psi_{2}]\|_{W_1^{\infty}(\R^2)}\leq C_{W_2^{\infty}}(\widetilde{\Xi_1^{-1}})\cdot\|\psi_2\|_{W_3^{\infty}(\R^2)}.$$
	Thus by the a priori bound Proposition \ref{a priori bound}, we get for $C>0$ (depending only on $\varepsilon$)
	\begin{equation}\label{eq:controle sur les correlateurs}
		\Big|\dfrac{P}{N}\Braket{\partial_1\widetilde{\Xi_1^{-1}}[\psi_{2}]}_{\mu_V\bigotimes\Delta\mu_N}\Big|\leq\dfrac{PCe^{K_V}}{N^{1+(1-\varepsilon)/2}}\left[C_{W_2^\infty}(\widetilde{\Xi_1^{-1}})+\|\rho_V\|_{\infty}^{1/2}C_{H^2}(\widetilde{\Xi_1^{-1}})\right] \mcal{N}_3(\psi_2).
	\end{equation}
	The following two terms in \eqref{deuxieme expansion} are also a $o(N^{-1})$ by the same reasons as before. By Proposition \ref{a priori bound}, the $3$-linear statistics is a $o(N^{-1})$ for $\varepsilon>0$ small enough. Hence, we obtain the expansion:
	$$\Braket{\psi}_{\Delta\mu_N}= \dfrac{d_{1}^{(1)}(\psi_1)}{N}+o(N^{-1})\quad\quad\quad \text{ and }\quad\quad\quad \Braket{\psi_2}_{\Delta\mu_N\otimes\Delta\mu_N}=\dfrac{d_{1}^{(2)}(\psi_2)}{N}+o(N^{-1})$$
	where $$d_{1}^{(1)}(\psi_1)\defi P\Braket{\widetilde{\Xi^{-1}}[\psi_1]'}_{\mu_V}+P\Braket{\Theta^{(2)} \circ\widetilde{\Xi_1^{-1}}\left[\partial_2\mcal{D}\circ\widetilde{\Xi^{-1}}[\psi_1]\right]  }_{\mu_V}$$
	and
	$$d_{1}^{(2)}(\psi_2)\defi -\Braket{\Theta^{(2)} \circ\widetilde{\Xi_1^{-1}}\big[\partial_2\psi_{2}\big] }_{\mu_V}.$$
	More generally, suppose the desired expansion for $\Braket{\psi_{k}}_{\overset{k}{\bigotimes}\Delta\mu_N}$ holds up to $o(N^{-n})$ for all $k\in\llbracket0,2n\rrbracket$ and for any function $\psi_k\in H^\infty(\R^k)$. Additionnaly, suppose that \eqref{eq:controleraminedercorrelator} is true for all $k\in\llbracket1,K-2\rrbracket$. Then, taking a general function $\psi_{2n+2}\in H^\infty(\R^{2n+2})$, the $(2n+2)$-\textit{th} equation involves the $(2n+3)$-linear statistic $\Braket{\mcal{D}_1\circ\Xi_1^{-1}[\psi_{2n+2}]}_{\overset{2n+3}{\bigotimes}\Delta\mu_N}$ (see \eqref{DSlvlnthm}). By Proposition \ref{a priori bound}, it is of size $o\left (N^{-(n+1)}\right )$ for $\varepsilon$ small enough. The other terms will be either, $(2n+1)$-linear statistics with a $N^{-1}$ prefactor and therefore behave like $o(N^{-(n+1)})$ for $\varepsilon$ small enough, or a $2n$-linear statistics with a prefactor $N^{-1}$. For the latter, by hypothesis, we know the asymptotic expansion up to $o(N^{-n})$, thus with the prefactor $N^{-1}$, we deduce the following expansion for $\braket{\psi_{2n+2}}_{\overset{2n+2}{\bigotimes}\Delta\mu_N}$:
	$$\braket{\psi_{2n+2}}_{\overset{2n+2}{\bigotimes}\Delta\mu_N}=\dfrac{d_{n+1}^{(2n+2)}(\psi_{2n+2})}{N^{n+1}}+o\left(N^{-(n+1)}\right).$$
	We now deduce from the equation above, the expansion of $\braket{\psi_{2n+1}}_{\overset{2n+1}{\bigotimes}\Delta\mu_N}$ for a general $\psi_{2n+1}$ belonging  to $H^\infty(\R^{2n+1})$. In the $(2n+1)$\textit{-th} loop equation \eqref{DSlvlnthm}, the $(2n+2)$-linear statistic yields a non-trivial term of order $N^{-(n+1)}$ \textit{i.e.}:
	$$\Braket{\mcal{D}_1\circ\Xi_1^{-1}[\psi_{2n+1}]}_{\overset{2n+2}{\bigotimes}\Delta\mu_N}=\dfrac{d_{n+1}^{(2n+2)}\left(\mcal{D}_1\circ\widetilde{\Xi_1^{-1}}[\psi_{2n+1}]\right) }{N^{n+1}}+o\left(N^{-(n+1)}\right).$$
	Again, the $(2n+1)$-linear statistics with a prefactor $N^{-1}$ is a $o(N^{-(n+1)})$. Finally for the $2n$  and $(2n-1)$ linear statistics with the prefactor $N^{-1}$ appearing in the $(2n+1)$\textit{-th} loop equation, their known expansion up to $o(N^{-(n+1)})$ (by hypothesis) leads to the expansion of $\braket{\psi_{2n+2}}_{\overset{2n+2}{\bigotimes}\Delta\mu_N}$ up to $o(N^{-(n+1)})$.
	
	To conclude on the expansion of the $2n$-linear statisticsup to $o(N^{-(n+1)})$, notice that for each term appearing in the $(2n)$-\textit{th} equation each term will either be a $(2n+1)$-linear statistics for which we know the expansion up to $o(N^{-(n+1)})$, or a linear statistic for which we know, by hypothesis, the asymptotic expansion up to $o(N^{-n})$ (of order $2n$, $2n-1$ or $2n-2$), preceded by a factor $N^{-1}$. We can therefore conclude on the existence of the expansion of the $2n$-linear statistics up to $o(N^{-(n+1)})$. Then applying the same arguments for $2n-1$,\dots $1$-linear statistics allows us to conclude that the induction step is established.
	
	Finally, to conclude on \eqref{eq:controleraminedercorrelator}, one must simply notice that for all $n\geq 1$ and $\psi_n\in H^\infty(\R^n)$, for all $K\geq\lceil n/2\rceil$,
	$$\Braket{\psi_n}_{\overset{n}{\bigotimes}\Delta\mu_N}-\sum_{i=\lceil n/2\rceil}^{K}\dfrac{d_{i}^{(n)}(\psi_n)}{N^i}=\dfrac{d_{K+1}^{(n)}(\psi_n)}{N^{K+1}}+\mfrak{R}_{K+1}^{(n)}(\psi_n)\hspace{0,5cm}\text{ and }\hspace{0,5cm}\msf{R}_{K+1}^{(n)}(\psi_n)=o(N^{-(K+1)}).$$
	Above, the remainder $\msf{R}_{K+1}^{(n)}(\psi_n)$ contains all the negligible statistics involving the operators $\widetilde{\Xi_1^{-1}}$, $\mcal{D}_1\circ\widetilde{\Xi_1^{-1}}$ and $\Theta^{(a)} \circ\widetilde{\Xi_1^{-1}}$. Thus just as in \eqref{eq:controle sur les correlateurs}, by using continuity of the different operators involved in each of the statistics, there exists $m>0$ (depending only on $K$ and $n$), a polynomial $\msf{Q}_{K}^{(1)}$ (whose coefficients only depend on $K$ and $n$) in $e^{K_V}$, $\|\rho_V\|_{\infty}^{1/2}$, $C_{H^{m}}(\widetilde{\Xi_1^{-1}})$, $C_{W_{m}^{\infty}}(\widetilde{\Xi_1^{-1}})$ and $ C_{H^{m}}(\Theta^{(2)}\circ\widetilde{\Xi_{1}^{-1}})$ with coefficients independent of $V$ and a constant $C>0$ (depending only on $K$ and $n$) such that:
	$$\bigg|\mathsf{R}_{K+1}^{(n)}[\psi_n]\bigg|\leq \dfrac{C}{N^{K+1}}\msf{Q}^{(1)}_{K}\cdot\mcal{N}_{m}(\psi_n).$$
	
	To bound, $d_{K+1}^{(n)}(\psi_n)$ and extract the $V$-dependence, one just notices that it is a sum of linear statistics, involving as before the previous operators. By contintuity of these operators, there exists a polynomial $\msf{Q}_{K}^{(2)}$ (whose coefficients are independent of $V$) in the previous operator norms and $\|\rho_V\|_{\infty}^{1/2}$ such that, choosing $m>0$ and $C>0$ (a constant independent of $V$ and $\psi_n$) big enough, such that:
	$$\bigg|d_{K+1}^{(n)}(\psi_n)\bigg|\leq C\msf{Q}^{(2)}_{K}\cdot\mcal{N}^{(n)}(\psi_n).$$
	Thus setting $ c=C\left[\msf{Q}^{(1)}_{K}+\msf{Q}^{(2)}_{K}\right] $ allows us to conclude about \ref{eq:controleraminedercorrelator}. The fact that $t\in[0,1]\mapsto C_{\mrm{rem}}(V_{\phi,t})$ is continuous  follows from the fact it is a polynomial in building blocks which are continuous as it is shown in Appendix \ref{appB}, Lemma \ref{lem:contk_V}, Proposition \ref{prop:constantescontinues}.
\end{proof}

\section{Parameter continuity of the equilibrium measure}\label{section7}
We want to conclude about the asymptotic expansion of $\log\mcal{Z}_N\left[V_{G,\phi}\right]$ for a smooth $\phi,$ by inserting the asymptotic expansion of the linear statistics of Theorem \ref{asymptotic correlators} in Lemma \ref{interpolation}. In order to make that step rigorous, it is necessary to prove that all the linear statistics integrated with respect to the probability measure $\mbb{P}_N^{V_{G,\phi,t}}$ with $t\in[0,1]$ and then $\diff t$, yield a definite and finite integral. Since all the quantities depend on $t$ through the equilibrium measure $\mu_{V_{G,\phi,t}}$, we first prove a continuity result for $t\mapsto\rho_{V_{G,\phi,t}}$. The result that we are going to prove does not depend on the specificity of the quadratic potential so in the following, we set $V_{\phi,t}:x\mapsto V(x)+t\phi(x)$ with $V$ satisfying the assumptions \ref{assumptions}.

While from the measure point of view, it is easy to show that $t\mapsto\mu_{V_{\phi,t}}$ is continuous for the weak topology of measures, it is not sufficient to deduce the continuity of the quantities involved in our problem. Indeed, in the controls we showed in Section \ref{section4}, quantities like $L^\infty$-norm of derivatives of $\rho_{V_t}$ and $C_{\mcal{L}}$ arose. To prove the integrability of these quantities, it is necessary to show that $\|\rho_{V_{\phi,t}}-\rho_{V_{\phi,t_0}}\|_{W_i^{\infty}(\R)}\tend{t\rightarrow t_0}0$ for all $t_0\in[0,1]$.

\subsection{Setting for Banach fixed-point theorem}\label{subs:cont density}
Let $\phi\in H^{\infty}(\R)$ and let $t_0\in[0,1]$, we define the function $u_t$ by

\begin{equation}\label{def:u_t}
	\rho_{V_{\phi,t}}=(1+\delta t u_t)\rho_{V_{\phi,t_0}}\;\;\;\; \textit{ i.e.    }\;\;\;\; u_t=\dfrac{\rho_{V_{\phi,t}}-\rho_{V_{\phi,t_0}}}{\delta t}\dfrac{1}{\rho_{V_{\phi,t_0}}}
\end{equation}
where $t\neq t_0$ and $\delta t\defi t-t_0$. We will show, by Banach fixed-point theorem, that $x\mapsto u_t(x)\in\mcal{C}^\infty(\R)$, by Lemma \ref{lem:regularitedensite}, is the unique fixed-point of a $t$ continuous operator. This will allow us to deduce that $t\mapsto u_t$ is continuous for the $W^\infty_{k}$-norm for all $k\geq0$. The continuity of $t\mapsto\rho_{V_t}\in W^\infty_{k}(\R)$ will then follow.

In order to construct the operator of interest, we start with the following lemma.
\begin{lemma}\label{lem:lambda_t}
	Let $t,t_0\in[0,1]$,
	\begin{multline}
		\label{eq:lambda_t}
		\lambda_{V_{\phi,t}}=\lambda_{V_{\phi,t_0}}+\delta t\int_\R\phi(x) \diff\mu_{V_{\phi,t_0}}(x)-2P\delta t\iint_{\R^2}\log|x-y|u_t(x)\diff\mu_{V_{\phi,t_0}}(x)\diff\mu_{V_{\phi,t_0}}(y)
		\\+\int_\R\left[\log\big(1+\delta t u_t(x)\big)-\delta t u_t(x)\right]  \diff\mu_{V_{\phi,t_0}}(x).
	\end{multline}
	Here $\lambda_{V_{\phi,t}}$ denotes the constant appearing in \eqref{eqmeasure:charac} with potential $V_{\phi,t}$.
\end{lemma}

\begin{proof}
	We integrate with respect to $\mu_{V_{\phi,t_0}}$ \eqref{eqmeasure:charac} to get
	\begin{multline*}
		\lambda_{V_{\phi,t}}=\int_\R V_{\phi,t_0}(x)\diff\mu_{V_{\phi,t_0}}(x)+\delta t\int_\R\phi(y) \diff\mu_{V_{\phi,t_0}}(y)-2P\iint_{\R^2}\log|x-y|\diff\mu_{V_{\phi,t}}(x)\diff\mu_{V_{\phi,t_0}}(y)
		\\+\int_\R \log\rho_{V_{\phi,t}}(x)\diff\mu_{V_{\phi,t_0}}(x).
	\end{multline*}
	After using the fact that $\displaystyle\int_\R u_t(x)\diff\mu_{V_{\phi,t_0}}(x)=0$, that $\rho_{V_{\phi,,t}}=(1+\delta t u_t)\rho_{V_{\phi,t_0}}$ and the characterization \eqref{eqmeasure:charac} of $\mu_{V_{\phi,t_0}}$, this yields the result.
\end{proof}
To show that $u_t$ is a fixed point of a $t$-continuous operator, we need to invert and control the operator $\mcal{T}\defi -\mcal{L}\circ\mcal{A}^{-1}$ (these operators were inroduced in Definition \ref{def:oper Aet W}) which appears naturally when comparing $\rho_{V_{\phi,t}}$ to $\rho_{V_{\phi,t_0}}$.
\begin{proposition}\label{prop:invefredholm}
	We define the operator $\mcal{T}$  by
	$\mcal{T}[v]\defi v-\mcal{K}[v]$ for all $v\in L^2(\mu_{V_{\phi,t_0}})$, 
	where
	$$\mcal{K}[v](x)\defi 2P\int_\R k(x,y)v(y)\rho_{V_{\phi,t_0}}(y)\diff y$$
	and $$k(x,y)\defi \left (\log\dfrac{|x-y|}{1+|x|}-\displaystyle\int_\R\log\dfrac{|z-y|}{1+|z|}\rho_{V_{\phi,t_0}}(z)\diff z\right).$$
	The operator $\mcal{T}:L^2(\mu_{V_{\phi,t_0}})\rightarrow L^2(\mu_{V_{\phi,t_0}})$ is bijective and for all $n\geq0$, $\mcal{T}\left[W_n^\infty(\R)\right]= W_n^\infty(\R)$.
	Finally, for all $n\in\N$, there exists $C_{\mcal{T},n}>0$ such that for any $v\in W_n^\infty(\R)$,
	\begin{equation}\label{eq:controleTLinfini}
		\left \|\mcal{T}^{-1}[v]\right \|_{W_n^\infty(\R)}\leq C_{\mcal{T},n}\cdot\|v\|_{W_n^\infty(\R)}.
	\end{equation}
	
\end{proposition}

\begin{proof}
	It was shown in \cite[Theorem 6.13]{DwoMemin}\footnote{There is a misprint in the cited Theorem, in (65) it should be $k(x,y)w(y)d\mu_V(y)$.} that $\mcal{T}[v]=-\mcal{L}\circ\mcal{A}^{-1}[v]$ for all $v\in\msf{H}$, where we used  that $\int_\R v(y)\diff\mu_{V_{\phi,t_0}}(y)=0$. Since $\mcal{A},\,\mcal{L}:\mcal{D}(\mcal{A})\rightarrow\msf{H}$ are both bijective operators, so is $\mcal{T}:\msf{H}\rightarrow\msf{H}$.
	
	The fact that $k$ verifies $\|k\|_{L^2(\otimes^{2}\mu_{V_{\phi,t_0}})}<+\infty$ implies that $\mcal{K}$ considered as an operator from $L^2(\mu_{V_{\phi,t_0}})$ to itself is an Hilbert-Schmidt operator thus compact. This implies that $\mcal{T}$ is a Fredholm operator. We now show that the kernel of $\mcal{T}:L^2(\mu_{V_{\phi,t_0}})\rightarrow L^2(\mu_{V_{\phi,t_0}})$ is trivial. Let $v\in L^2(\mu_{V_{\phi,t_0}})$ such that $\mcal{T}[v]=0$ \textit{i.e.} $v=\mcal{K}[v]$, then $\mcal{K}[v]$ is in $H^1(\mu_{V_{\phi,t_0}})$. Indeed we have $$\mcal{K}[v]'(x)=-\mcal{H}[v\rho_{V_{\phi,t_0}}](x)-\dfrac{\mrm{sgn}(x)}{1+|x|}\int_{\R}v(y)\diff\mu_{V_{\phi,t_0}}(y)\in L^2(\mu_{V_{\phi,t_0}}).$$
	Moreover since $\int_\R\mcal{K}[v](x)\diff\mu_{V_{\phi,t_0}}(x)=0$, we conclude that $v\in\msf{H}$. We can now conclude that $v=0$ by bijectivity of $\mcal{T}$ on $\msf{H}$. Finally, by Fredhom alternative, $\mcal{T}$ is invertible on $L^2(\mu_{V_{\phi,t_0}})$ since it is injective. 
	
	We now prove that for all $n\in\N$, $\mcal{T}\left[W_n^\infty(\R)\right]= W_n^\infty(\R)$. We proceed by induction. For $n=0$, let $f\in L^\infty(\R)\subset L^2(\mu_{V_{\phi,t_0}})$. There exists a unique $v\in L^2(\mu_{V_{\phi,t_0}})$ such that $\mcal{T}[v]=f$ so $v=f+\mcal{K}[v]$ but since $f$ and $\mcal{K}[v]$ are bounded, so is $v\in L^\infty(\R)$. Reciprocally, if $v\in L^\infty(\R)$ so is $\mcal{T}[v]$, hence $\mcal{T}\left[L^\infty(\R)\right]= L^\infty(\R)$. Finally, for all $v\in L^\infty(\R)$,
	\begin{equation}\label{eq:t linfini}
		\Big\|\mcal{T}[v]\Big\|_{L^\infty(\R)}\leq\left(1+2P\max_{x\in\R}\int_{\R}|k(x,y)|\diff\mu_{V_{\phi,t_0}}(y)\right)\cdot\|v\|_{L^\infty(\R)}.
	\end{equation}
	Now suppose $\mcal{T}\left[W_n^\infty(\R)\right]= W_n^\infty(\R)$ is true and let's show it for $n+1$. Let $f \in W_{n+1}^\infty(\R)\subset W_n^\infty(\R)$, so by hypothesis, there exists $v\in W_n^{\infty}(\R)$ such that:
	\begin{multline}\label{eq:t(v)^n}
		f^{(n)}(x)=v^{(n)}(x)+2P\left(\log(1+|.|)\right)^{(n)}(x)\int_\R v(y)\diff\mu_{V_{\phi,t_0}}(y)-2P\int_\R\log|x-y| (v\rho_{V_{\phi,t_0}})^{(n)}(y)\diff y
		\\+2P\delta_{n,0}\iint_{\R^2}\dfrac{\log|z-y|}{1+|z|}v(y)\diff\mu_{V_{\phi,t_0}}(y)\diff\mu_{V_{\phi,t_0}}(z).
	\end{multline}
	We deduce that $v^{(n)}$ is differentiable of derivative
	\begin{equation}\label{eq:v n+1}
		v^{(n+1)}(x)=f^{(n+1)}(x)-2P\left(\log(1+|.|)\right)^{(n+1)}(x)\int_\R v(y)\diff\mu_{V_{\phi,t_0}}(y)-2P\mcal{H}\left[(v\rho_{V_{\phi,t_0}})^{(n)}\right](x),
	\end{equation}
	where $\left(\log(1+|.|)\right)^{(n)}(x)=\dfrac{n!\mrm{sgn}(x)^{(n)}}{(1+|x|)^n}$. Since the two first terms in the RHS of \eqref{eq:v n+1} are clearly bounded, we just have to show that  $\mcal{H}\big[(v\rho_{V_{\phi,t_0}})^{(n)}\big]\in H^1(\R)$. By boundedness of $f^{(n+1)}$ and $v^{(i)}$ for all $i\leq n$, we have
	\begin{multline*}
		(v\rho_{V_{\phi,t_0}})^{(n+1)}=\rho_{V_{\phi,t_0}}\left(v^{(n+1)}+\sum_{k=0}^{n}\binom{n+1}{k}\dfrac{\rho_{V_{\phi,t_0}}^{(n+1-k)}}{\rho_{V_{\phi,t_0}}}v^{(k)}\right) 
		\\=\rho_{V_{\phi,t_0}}\Bigg(f^{(n+1)}-2P\left(\log(1+|.|)\right)^{(n+1)}\int_\R v(y)\diff\mu_{V_{\phi,t_0}}(y)-2P\mcal{H}\left[(v\rho_{V_{\phi,t_0}})^{(n)}\right]
		\\+\sum_{k=0}^{n}\binom{n+1}{k}\dfrac{\rho_{V_{\phi,t_0}}^{(n+1-k)}}{\rho_{V_{\phi,t_0}}}v^{(k)}\Bigg)\in L^2(\R).
	\end{multline*}
	Thus, it holds that $\mcal{H}\left[(v\rho_{V_{\phi,t_0}})^{(n)}\right]\in H^1(\R)$ and that it is bounded. Hence it proves that $v^{(n+1)}\in L^\infty(\R)$ and hence that $W_{n+1}^\infty(\R)\subset \mcal{T}\left[W_{n+1}^\infty(\R)\right]$. Conversely, if $v\in W_{n+1}^\infty(\R)$, then $f\in W_{n}^\infty(\R)$ by hypothesis and just as before, we show that \eqref{eq:v n+1} holds. We conclude that $f^{(n+1)}\in L^\infty(\R)$ again by showing that $\mcal{H}\left[(v\rho_{V_{\phi,t_0}})^{(n)}\right]$ is bounded by the fact that $v\in W_{n+1}^\infty(\R)$. This establishes that $W_{n+1}^\infty(\R)=\mcal{T}\left[W_{n+1}^\infty(\R)\right]$.
	
	Thus for all $n\in\N$, $\mcal{T}:W_n^\infty(\R)\rightarrow W_n^\infty(\R)$ is a bijective operator. Furthermore, it is a bounded operator by the fact that for all $1\leq i\leq n$, there exists $C>0$ such that for all $v\in W_n^{\infty}(\R)$, by Leibniz formula:
	\begin{multline*}
		\big|\mcal{T}[v]^{(i)}(x)\big|\leq \|v^{(i)}\|_{L^\infty(\R)}+2P\big\|\left(\log(1+|.|)\right)^{(i)}\big\|_{L^\infty(\R)}\cdot\|v\|_{L^\infty(\R)}\\+2P\|v\|_{W_i^{\infty}(\R)}\cdot\sum_{k=0}^{i}\binom{i}{k}\sup_{z\in\R}\int_{\R}\Big|\log\dfrac{|z-y|}{1+|z|}\Big|\cdot\Bigg|\dfrac{\rho_{V_{\phi,t_0}}^{(k)}}{\rho_{V_{\phi,t_0}}}(y)\Bigg|\diff\mu_{V_{\phi,t_0}}(y)\leq C\|v\|_{W_i^{\infty}(\R)}.
	\end{multline*} 
	Above we used \eqref{eq:t(v)^n} and $\int_{\R}(v\rho_{V_{\phi,t_0}})^{(i)}(y)\diff y=0$ to deduce the following fact:
	\begin{multline*}
		\int_\R\log|x-y|(v\rho_{V_{\phi,t_0}})^{(i)}(y)\diff y=\int_\R\log|x-y|(v\rho_{V_{\phi,t_0}})^{(i)}(y)\diff y-\log(1+|x|)\int_\R(v\rho_{V_{\phi,t_0}})^{(i)}(y)\diff y
		\\=\int_\R\log\dfrac{|x-y|}{1+|x|}(v\rho_{V_{\phi,t_0}})^{(i)}(y)\diff y.
	\end{multline*} 
	Thus, we conclude that $\max_{i\in\llbracket1,n\rrbracket}\big\|\mcal{T}[v]^{(i)}\big\|_{L^\infty(\R)}\leq C\|v\|_{W_n^{\infty}(\R)}$. The bound on $\big\|\mcal{T}[v]\big\|_{L^\infty(\R)}$ was shown in \eqref{eq:t linfini}. We finally conclude that $\mcal{T}:W_n^\infty(\R)\rightarrow W_n^\infty(\R)$ is bounded bijective between Banach spaces and by Banach isomorphism theorem so is $\mcal{T}^{-1}$, this establishes \eqref{eq:controleTLinfini} and completes the proof.
\end{proof}

\begin{remark}
	An explicit expression for $\mcal{T}^{-1}$ is available using Fredholm determinant theory for invertible Hilbert-Schmidt operators, see \cite[Section XII]{gohberg2012traces}. For all $v\in L^2_0(\mu_{V_{\phi,t_0}})$,
	\begin{multline}\label{eq:inverse fredholm}
		\mcal{T}^{-1}[v](x)\defi v(x)\\+\dfrac{1}{\underset{2}{\det}(I-\mcal{K})}\sum_{n\geq1}\dfrac{(-1)^n}{n!}\int_{\R^{n+1}}\begin{vmatrix}
			k(x,s)&k(x,t_1)&\dots&k(x,t_n)
			\\k(t_1,s)&0&\dots&k(t_1,t_n)
			\\\vdots& & &\vdots
			\\k(t_n,s)&k(t_n,t_1)&\hdots&0
		\end{vmatrix}v(s)\diff\mu_{V_{\phi,t_0}}(s)\prod_{i=1}^{n}\diff\mu_{V_{\phi,t_0}}(t_i).
	\end{multline}
	Above $\underset{2}{\det}$ stands for the $2$-determinant. This formula was established in \cite[Theorem 6.11]{DwoMemin}
\end{remark}

We are now able to show that $u_t$ is a fixed point of a certain operator. We recall that $\phi$ was introduced in the beginning of subsection \ref{subs:cont density}.

\begin{proposition}\label{prop:def V_t, U_t}
	For all $t\in[0,1]$, $u_t$ defined in \eqref{def:u_t} is the unique measurable function such that $\int_\R u_t(x)\diff\mu_{V_{\phi,t_0}}(x)=0$ and which satisfies:
	$$u_t=\mcal{T}^{-1}\circ\mcal{V}_t[u_t]$$
	where $\mcal{V}_t[u]\defi -\phi+\int_\R\phi(y) \diff\mu_{V_{\phi,t_0}}(y)+\delta t\,\mcal{U}_t[u]$,   and
	\begin{multline*}
		\mcal{U}_t[v](x)\defi \left(-\phi(x)+\int_\R\phi(y) \diff\mu_{V_{\phi,t_0}}(y)+\mcal{K}[v](x)+\int_\R\dfrac{\log\Big(1+\delta t v(y)\Big)-\delta t v(y)}{\delta t} \diff\mu_{V_{\phi,t_0}}(y)\right)^2
		\\\times\int_{0}^{1}(1-s)\diff s \exp\Bigg\{s\delta t\left(-\phi(x)+\int_\R\phi(y) \diff\mu_{V_{\phi,t_0}}(y)+\mcal{K}[v](x)\right) 
		\\+\int_\R\frac{\log\Big(1+\delta t v(y)\Big)-\delta t v(y)}{\delta t} \diff\mu_{V_{\phi,t_0}}(y)\Bigg\} +\int_\R\dfrac{\log\Big(1+\delta t v(y)\Big)-\delta t v(y)}{(\delta t)^2} \diff\mu_{V_{\phi,t_0}}(y).
	\end{multline*}
\end{proposition}

\begin{proof}
	Lemma \ref{lem:lambda_t} allows one to substitute $\lambda_{V_{\phi,t}}$ in the representation for $\rho_{V_{\phi,t}}$ by \eqref{eq:lambda_t} hence leading to
	\begin{multline*}
		\rho_{V_{\phi,t}}=(1+\delta t u_t)\rho_{V_{\phi,t_0}}=\exp\Bigg(-V_{\phi,t_0}-2PU^{\rho_{V_{\phi,t_0}}}+\lambda_{V_{\phi,t_0}}-\delta t\phi+\delta t\int_\R\phi(y) \diff\mu_{V_{\phi,t_0}}(y)-2P\delta t U^{u_t\rho_{V_{\phi,t_0}}}\\-2P\delta t\iint_{\R^2}\log|
		y-z|u_t(z)\diff\mu_{V_{\phi,t_0}}(z) \diff\mu_{V_{\phi,t_0}}(y)+\int_\R\left[\log\Big(1+\delta t u_t(y)\Big)-\delta t u_t(y)\right]  \diff\mu_{V_{\phi,t_0}}(y)\Bigg).
	\end{multline*}
	Identifying $\rho_{V_{\phi,t_0}}$ via the first three terms in the exponential, $u_t$ has to satisfy the following relation for all $x\in\R$,
	\begin{multline*}
		1+\delta t u_t(x)=\exp\Bigg\{\delta t\Big(-\phi(x)+\int_\R\phi(y) \diff\mu_{V_{\phi,t_0}}(y)+2P \int_\R\log\dfrac{|x-y|}{1+|x|}u_t(y)\diff\mu_{V_{\phi,t_0}}(y)
		\\-2P\iint_{\R^2}\log\dfrac{|y-z|}{1+|y|}u_t(z)\diff\mu_{V_{\phi,t_0}}(z)\diff\mu_{V_{\phi,t_0}}(y)+\dfrac{1}{\delta t}\int_\R\left[\log\left(1+\delta t u_t(y)\right)-\delta t u_t(y)\right]  \diff\mu_{V_{\phi,t_0}}(y)\Big) \Bigg\}.
	\end{multline*}
	Above, we have used that $-U^{u_t\rho_{V_{\phi,t_0}}}(x)=\displaystyle\int_\R\log\dfrac{|x-y|}{1+|x|} u_t(y)\diff\mu_{V_{\phi,t_0}}(y)$ which is justified by the fact that $\displaystyle\int_\R u_t(y)\diff\mu_{V_{\phi,t_0}}(y)=0$. 
	Conversely, any $u$ such that $\int_\R u(y)\diff\mu_{V_{\phi,t_0}}(y)=0$ and satisfying the previous relation, verifies for all $x\in\R$,
	$$V_t(x)+2PU^{w}(x)+\log w(x)=\int_\R\Big\{V_t(y)+2PU^{w}(y)+\log w(y)\Big\}\diff\mu_{V_{\phi,t_0}}(y)$$
	where we have set $w\defi (1+\delta tu)\rho_{V_{\phi,t_0}}$. Because of this equation, $w$ can be written in exponential form as in \eqref{eqmeasure:exp}, it is thus positive and of mass $1$ which makes $\diff\mu(x)\defi w(x)\diff x$ a probability measure which satisfies the equation characterizing $\mu_{V_{\phi,t}}$. Thus, by uniqueness of the solution of \eqref{eqmeasure:charac}, $ \mu_{V_{\phi,t}}=\mu$ and thus $u=u_t$.
	
	We now expand $\exp$ into its Taylor-integral series of order 2, \textit{i.e.} $e^x=1+x+x^2\int_{0}^{1}(1-s)e^{sx}\diff s $. By using that: $$\mcal{K}[v]=2P\displaystyle\int_\R\log\dfrac{|x-y|}{1+|x|} u(y)\diff\mu_{V_{\phi,t_0}}(y)-2P\iint_{\R^2}\log\dfrac{|y-z|}{1+|y|}v(z)\diff\mu_{V_{\phi,t_0}}(z)\diff\mu_{V_{\phi,t_0}}(y)$$
	and $\mcal{T}[v]=v-\mcal{K}[v]$, we get:
	\begin{multline*}
		\mcal{T}[u_t](x)=-\phi(x)+\int_\R\phi \diff\mu_{V_{\phi,t_0}}
		\\+\delta t\Bigg[\left(-\phi(x)+\int_\R\phi \diff\mu_{V_{\phi,t_0}}+\mcal{K}[u_t](x)+\dfrac{1}{\delta t}\int_\R\left[\log\left(1+\delta t u_t\right)-\delta t u_t\right]  \diff\mu_{V_{\phi,t_0}}\right)^2
		\\\times\int_{0}^{1}\exp\left(s\delta t\left(-\phi(x)+\int_\R\phi \diff\mu_{V_{\phi,t_0}}+\mcal{K}[u_t](x) +\int_\R\frac{\log\left(1+\delta t u_t\right)-\delta t u_t}{\delta t} \diff\mu_{V_{\phi,t_0}}\right)\right) (1-s)\diff s 
		\\+\int_\R\dfrac{\log\left(1+\delta t u_t\right)-\delta t u_t}{\left(\delta t\right)^2 } \diff\mu_{V_{\phi,t_0}}\Bigg].
	\end{multline*}
	We next use the invertibility of $\mcal{T}$ to conclude.
\end{proof}

The next theorem shows that for each $t$ sufficiently close to $t_0$, $\mcal{T}^{-1}\circ\mcal{V}_t$ is contractive on a ball of fixed radius. Let $n\geq0$,
denote for all $R>0$, $\overline{\mscr{B}_n(0,R)}$ the closed ball of radius $R_n$ for the $W^{\infty}_n(\R)$-norm. 

\begin{theorem}\label{thm:contractivity}
	For all $\varepsilon_n>0$ small enough, there exists $R_n>0,$ such that for all $t\in]t_0-\varepsilon_n,t_0+\varepsilon_n[$, the operator $$\mcal{T}^{-1}\circ\mcal{V}_t:\left(\overline{\mscr{B}_n(0,R_n)},\|.\|_{W^{\infty}_n(\R)}\right) \longrightarrow \left(\overline{\mscr{B}_n(0,R_n)},\|.\|_{W^{\infty}_n(\R)}\right)$$ is well-defined and continuous. Furthermore it is contractive, \textit{i.e.} there exists $k_n\in]0,1[$, such that for all $v,w\in \overline{\mscr{B}_n(0,R_n)}$, $$\left \|\mcal{T}^{-1}\circ\mcal{V}_t[v]-\mcal{T}^{-1}\circ\mcal{V}_t[w]\right \|_{W^{\infty}_n(\R)}\leq k_n\|v-w\|_{W^{\infty}_n(\R)}.$$ Moreover $k_n$ is independent of $t$ on $]t_0-\varepsilon_n,t_0+\varepsilon_n[$ for $\varepsilon_n>0$ small enough.
	

\end{theorem}

\begin{proof}Let $\varepsilon>0$, $t\in]t_0-\varepsilon,t_0+\varepsilon[$, $v\in\overline{\mscr{B}_0(0,R)}$ for an arbitrary $0<R<\dfrac{1}{\delta t}$ (because otherwise the term $\int_\R\log(1+\delta tv)\rho_{V_{\phi,t_0}}$ in $\mcal{V}_t[v]$ might be ill-defined). We first show that $\|\mcal{V}_t[v]\|_{L^\infty(\R)}<+\infty$.
First, by Taylor-Lagrange inequality, we deduce that
	$$\Big|\int_\R\dfrac{\log\left[1+\delta t v(y)\right]-\delta t v(y)}{\delta t} \diff\mu_{V_{\phi,t_0}}(y)\Big|\leq \dfrac{\|v\|_{{L^\infty(\R)}}^2\delta t}{2}\sup_{|x|\leq\|v\|_{\infty}}\dfrac{1}{(1+\delta t x)^2}\leq \dfrac{R^2\delta t}{2(1-\delta t R)^2}.$$
	By recalling the definition of $\mcal{V}_t$ in Proposition \ref{prop:def V_t, U_t} and using the convexity of $x\mapsto x^2$, we get:
	\begin{multline*}
		\left \|\mcal{V}_t[v]\right \|_{L^\infty(\R)}\leq \delta t\Bigg[3\Bigg(4\|\phi\|_{L^\infty(\R)}^2+R^2\left \|\int_\R\left|\log\dfrac{|.-y|}{1+|.|}\right|\diff\mu_{V_{\phi,t_0}}(y) \right\|_{L^\infty(\R)}^2+\delta t^2\dfrac{R^4}{4(1-\delta t R)^4} \Bigg)
		\\\times\exp\left\{\delta t\left( 2\|\phi\|_{L^\infty(\R)}+R\left \|\int_\R\left|\log\dfrac{|.-y|}{1+|.|}\right|\diff\mu_{V_{\phi,t_0}}(y) \right\|_{L^\infty(\R)} +\delta t\dfrac{R^2}{2(1-\delta t R)^2}\right)\right\}+\dfrac{R^2}{2(1-\delta t R)^2} \Bigg]
		\\+2\|\phi\|_{L^\infty(\R)}.
	\end{multline*}
	The RHS is of the form $2\|\phi\|_{\infty}+\delta t g_0(\delta t)$ where $g_0$ is a positive function and with these notations:
	$$\big\|\mcal{T}^{-1}\circ\mcal{V}_t[v]\big\|_{L^\infty(\R)}\leq\mcal{C}_{\mcal{T},0}\left(2\|\phi\|_{L^\infty(\R)}+\delta t g_0(\delta t)\right).$$
	Therefore, by choosing $t$ such that  $|t-t_0|<\varepsilon_0$ for $\varepsilon_0>0$ small enough, there exists $R_0>0$ such that $2\|\phi\|_{\infty}C_{\mcal{T},0}<R_0<\dfrac{1}{\varepsilon_0}$ and $$\mcal{T}^{-1}\circ\mcal{V}_t\left(\overline{\mscr{B}_0(0,R_0)}\right)\subset\overline{\mscr{B}_0(0,R_0)}.$$ This makes the operator  $\mcal{T}^{-1}\circ\mcal{V}_t:\left(\overline{\mscr{B}_n(0,R_n)},\|.\|_{W^{\infty}_n(\R)}\right) \longrightarrow \left(\overline{\mscr{B}_n(0,R_n)},\|.\|_{W^{\infty}_n(\R)}\right)$ well-defined for all $|t-t_0|<\varepsilon_0$. For the contractivity, let $u,v\in\overline{\mscr{B}(0,R_0)}$, we get by Proposition \ref{prop:invefredholm},
	where $\mcal{U}_t$ was defined in Proposition \ref{prop:def V_t, U_t}. We now want to control $\big|\mcal{U}_t[u](x)-\mcal{U}_t[v](x)\big|$. We have, by decomposing the sum:
	$$\big|\mcal{U}_t[u](x)-\mcal{U}_t[v](x)\big|\leq \left(\Delta_1+\Delta_2(x)\msf{E}[u](x) +\Delta_3(x)\mfrak{h}[v](x)^2\right) $$
	where
	\begin{align*}
		\mfrak{h}[w](x)&\defi -\phi(x)+\int_\R\phi(y) \diff\mu_{V_{\phi,t_0}}(y)+\mcal{K}[w](x) +\int_\R\frac{\log\left(1+\delta t w(y)\right)-\delta t w(y)}{\delta t} \diff\mu_{V_{\phi,t_0}}(y)
		\\\msf{E}[w](x)&\defi \int_{0}^{1}\exp\Big[s\delta t\mfrak{h}[w](x)\Big] (1-s)\diff s 
		\\\Delta_1&\defi \int_\R\left|\dfrac{\log(1+\delta t u(y))-\delta tu(y)-\log(1+\delta t v(y))+\delta t v(y)}{(\delta t)^2}\right|\diff\mu_{V_{\phi,t_0}}(y)
		\\\Delta_2(x)&\defi \big|\mfrak{h}[u](x)^2-\mfrak{h}[v](x)^2\big|
		\\\Delta_3(x)&\defi \big|\msf{E}[u](x)-\msf{E}[v](x)\big|.
	\end{align*}
	First, Taylor Lagrange inequality leads to:
	
	$$\Delta_1\leq\sup_{y\in[-R_0,R_0]}\dfrac{|y|}{1+\delta t y}\|u-v\|_{L^\infty(\R)}\leq\dfrac{R_0}{1-\delta t R_0}\|u-v\|_{L^\infty(\R)}.$$
	Furthermore, by using $a^2-b^2=(a+b)(a-b)$, that $\mfrak{h}[u]$ and $\mfrak{h}[v]$ are bounded
	\begin{multline*}
		\Delta_2\leq 2 \max\Big(\|\mfrak{h}[u]\|_{L^\infty(\R)},\|\mfrak{h}[v]\|_{L^\infty(\R)}\Big)\Bigg[\Big|\mcal{K}[u-v](x)\Big|+\delta t \Delta_1\Bigg]
		\\\leq\left(4\|\phi\|_{L^\infty(\R)}+8PR_0\left \|\int_\R\left|\log\dfrac{|.-y|}{1+|.|}\right|\diff\mu_{V_{\phi,t_0}}(y) \right\|_{L^\infty(\R)}+\dfrac{\delta tR_0^2}{(1-\delta t R_0)^2}\right) \|u-v\|_{L^\infty(\R)}
		\\\times\left(\dfrac{R_0\delta t}{1-\delta t R_0}+4P\left \|\int_\R\left|\log\dfrac{|.-y|}{1+|.|}\right|.\rho_{V_{\phi,t_0}}(y)\diff y \right\|_{L^\infty(\R)}\right)\leq C(R_0)\|u-v\|_{L^\infty(\R)}
	\end{multline*}
	Similarly, there exists $C(R_0)>0$ such that:
	\begin{equation*}
		\Delta_3\|\mfrak{h}[v]\|_{L^\infty(\R)}\leq C(R_0)\delta t\|u-v\|_{L^\infty(\R)}
	\end{equation*}
	which finally leads to the existence of $C_0>0$ such that:
	\begin{equation*}
		\big\|\mcal{U}_t[u]-\mcal{U}_t[v]\big\|_{L^\infty(\R)}\leq C_0\cdot|\delta t|\cdot\|u-v\|_{L^\infty(\R)}.
	\end{equation*}
	We now choose $\varepsilon_0$ again small enough such that $k_0\defi \delta tC_0C_{\mcal{T},0}<1$, this concludes the proof that $\mcal{T}^{-1}\circ\mcal{V}_t$ is contractive on $\overline{\mscr{B}_0(0,R_0)}$ with contractivity constant $k_0$.  
	
	To get the contractivity property for $\mcal{T}^{-1}\circ\mcal{V}_t$ on $W_n^{\infty}(\R)$, we adapt a similar strategy.
	Let $u\in\overline{\mscr{B}_n(0,R)}$ with $\delta tR<1$, by Proposition \ref{prop:invefredholm},  
	$$\|\mcal{T}^{-1}[\mcal{V}_t[u]]\|_{W_n^\infty(\R)}\leq C_{\mcal{T},n}\|\mcal{V}_t[u]\|_{W_n^\infty(\R)}.$$
	Furthermore, it is clear that every term appearing in the definition of $\mcal{V}_t[u]$ belongs to $W_n(\R)$, thus by the same argument as before there exists a positive function $g_n$ such that, $$\|\mcal{T}^{-1}\circ\mcal{V}_t[u]\|_{W_n(\R)}\leq C_{\mcal{T},n}\Big(2\|\phi\|_{W_n^{\infty}(\R)}+\delta t g_n(R)\Big).$$ We conclude just as before that by taking $\delta t$ small enough, $\mcal{T}^{-1}\circ\mcal{V}_t:\overline{\mscr{B}_n(0,R)}\rightarrow \overline{\mscr{B}_n(0,R)}$ is well-defined.
	Finally, just as before since for all $u,v\in W_n^{\infty}(\R)$
	\begin{align*}
		\mcal{U}_t[u]^{(n)}(x)&=\delta t\sum_{k=0}^{n}\binom{n}{k}\left(\mfrak{h}[u]^2\right)^{(k)}(x)\msf{E}[u]^{(n-k)}(x)
		\\&=\delta t\sum_{0\leq i\leq k\leq n}\binom{n}{k}\binom{k}{i}\mfrak{h}[u]^{(i)}(x)\mfrak{h}[u]^{(k-i)}(x)\msf{E}[u]^{(n-k)}(x).
	\end{align*}
	Moreover, by the same controls as before it is easy to derive that for all $0\leq i\leq k\leq n$, for all $u,v\in \overline{\mscr{B}_n(0,R_n)}$, $$\|\mfrak{h}[u]^{(i)}\mfrak{h}[u]^{(k-i)}\msf{E}[u]^{(n-k)}-\mfrak{h}[v]^{(i)}\mfrak{h}[v]^{(k-i)}\msf{E}[v]^{(n-k)}\|_{L^\infty(\R)}\leq C(R_n,n,i,k)\|u-v\|_{W_n^{\infty}(\R)}.$$ This is enough to conclude that $$\left \|\mcal{T}^{-1}\circ\mcal{V}_t[u]-\mcal{T}^{-1}\circ\mcal{V}_t[v]\right \|_{W_n^{\infty}(\R)}\leq \delta tC_{\mcal{T},n}C(R_n,n)\|u-v\|_{W_n^{\infty}(\R)}.$$
	Finally, by taking $\delta t$ small enough, we conclude that $\mcal{T}^{-1}\circ\mcal{V}_t$ is contractive on $\overline{\mscr{B}_n(0,R_n)}$ with contractivity constant  $k_n\defi \delta tC_{\mcal{T},n}C(R_n,n)<1$.
\end{proof}

\begin{remark}Note that the definition of $u_t$ as an element of $W^{\infty}_n(\R)$ depends on $t_0$ and that we only proved the characterization of $u_t$ as a fixed point for $t\in]t_0-\varepsilon_n,t_0+\varepsilon_n[\setminus\{t_0\}$ with $\varepsilon_n>0$ small enough (we stress that we successively lowered $\delta t$ when increasing $n$.) Furthermore, since  for all $v\in L^\infty(\R)$
	$$\mcal{T}^{-1}\circ\mcal{V}_{t_0}[v]=-\phi+\displaystyle\int_\R\phi(y) \diff\mu_{V_{\phi,t_0}}(y),$$
	then we can set $u_{t_0}\defi -\phi+\displaystyle\int_\R\phi(y) \diff\mu_{V_{\phi,t_0}}(y)$ which is obviously the unique fixed point in $ \displaystyle\bigcap_{n\in\N}W_n^\infty{(\R)}$ of $\mcal{T}^{-1}\circ\mcal{V}_{t_0}$.
\end{remark}
\subsection{Regularity of the equilibrium measure}
We now prove the continuity of $t\mapsto u_t\in W_n^{\infty}(\R)$.

\begin{lemma}\label{lem:tv cont par rapport a t}Let $n\in\N$, $ \varepsilon_n>0$ and $R_n>0$ be as in Theorem \ref{thm:contractivity}. For all $v\in\overline{\mscr{B}_n(0,R_n)}$, for all $(t_p)_p\in\left(]t_0-\varepsilon_n,t_0+\varepsilon_n[\right)^{\N}$ such that $t_p\tend{p\rightarrow\infty}t\in]t_0-\varepsilon_n,t_0+\varepsilon_n[$ then $$\|\mcal{T}^{-1}\circ\mcal{V}_{t_p}[v]-\mcal{T}^{-1}\circ\mcal{V}_{t}[v]\|_{W_n^{\infty}(\R)}\underset{p\rightarrow \infty}{\longrightarrow}0.$$    
\end{lemma}

\begin{proof}
	The proof is based on the fact that there exists a neighboorhood $U_n$ of $t_0$ such that for all $v\in \overline{\mscr{B}_n(0,R_n)}$, $t\in U_n\mapsto\mcal{V}_t[v]\in W_n^{\infty}(\R)$ is continuous. Since $\mcal{T}^{-1}$ is also continuous in $W_n^{\infty}(\R)$-norm, we can conclude.
\end{proof}

\begin{corollary}[Continuity of the fixed-point]\label{cor:contfixe}Let $n\in\N$, for all $t,t'\in]t_0-\varepsilon_n,t_0+\varepsilon_n[$, $$\|u_t-u_{t'}\|_{W_n^{\infty}(\R)}\underset{t\rightarrow t'}{\longrightarrow}0.$$              
\end{corollary}

\begin{proof}Let $(t_n)_n\in\left(]t_0-\varepsilon,t_0+\varepsilon[\right)^\N $ such that $t_n\underset{n\rightarrow\infty}{\longrightarrow}t\in]t_0-\varepsilon_n,t_0+\varepsilon_n[$. First $$\left \|u_{t_n}-\mcal{T}^{-1}\circ\mcal{V}_{t_n}[u_{t}]\right \|_{W_n^{\infty}(\R)}=\left \|\mcal{T}^{-1}\circ\mcal{V}_{t_n}[u_{t_n}]-\mcal{T}^{-1}\circ\mcal{V}_{t_n}[u_{t}]\right \|_{W_n^{\infty}(\R)}\leq k_i\left \|u_{t_n}-u_t\right \|_{W_n^{\infty}(\R)}$$
	and by the triangle inequality, we obtain:
	\begin{equation*}
		\left \|u_{t}-\mcal{T}^{-1}\circ\mcal{V}_{t_n}[u_{t}]\right \|_{W_n^{\infty}(\R)}\geq \left \|u_{t_n}-u_t\right \|_{\infty}-\left \|u_{t_n}-\mcal{T}^{-1}\circ\mcal{V}_{t_n}[u_{t}]\right \|_{W_n^{\infty}(\R)}
		\geq (1-k_n)\left \|u_{t_n}-u_t\right \|_{W_n^{\infty}(\R)}.
	\end{equation*}
	Nevertheless by Lemma \ref{lem:tv cont par rapport a t}, $\left \|u_{t}-\mcal{T}^{-1}\circ\mcal{V}_{t_n}[u_{t}]\right \|_{W_n^\infty(\R)}\underset{n\rightarrow\infty}{\longrightarrow}0$ because $u_t$ is the fixed point of $\mcal{T}^{-1}\circ\mcal{V}_t$ which is a continuous operator with respect to $t$ hence $\left \|u_{t_n}-u_t\right \|_{W_n^\infty(\R)}\underset{n\rightarrow\infty}{\longrightarrow}0$.
\end{proof}

\begin{theorem}\label{thm:locallimittheorem}
	Let $t,t_0\in[0,1]$,
	$$\boxed{\left \|\rho_{V_{\phi,t}}-\rho_{V_{\phi,t_0}}\right \|_{W_n^{\infty}(\R)}\underset{t\rightarrow t_0}{\longrightarrow}0.}$$
	Furthermore, for all $x\in\R$, $k\in\N$,  $t\mapsto\partial_x^{k}\rho_{V_{\phi,t}}(x)\in\mcal{C}^\infty(\R)$ and satisfies the following partial integro-differential equation:
	$$\partial_t\partial_{x}^{k}\rho_{V_{\phi,t}}(x)=\partial_x^{k}\left[ \left(-\phi+\int_\R\phi(s)\rho_{V_{\phi,t}}(s)\diff s \right) \rho_{V_{\phi,t}}\right] (x).$$
\end{theorem}

\begin{proof} By setting $u_t\defi \dfrac{\rho_{V_{\phi,t}}-\rho_{V_{\phi,t'}}}{t-t'}\dfrac{1}{\rho_{V_{\phi,t'}}}$,
	
	$$\|\rho_{V_{\phi,t}}-\rho_{V_{\phi,t'}}\|_{W_n^{\infty}(\R)}=|t-t'|\cdot\left \|\rho_{V_{\phi,t'}}u_t \right \|_{W_n^{\infty}(\R)}\leq 2^n|t-t'|\cdot \left \|\rho_{V_{\phi,t'}}\right \|_{W_n^{\infty}(\R)}\cdot\left \|u_t \right \|_{W_n^{\infty}(\R)}.$$
	By Corollary \ref{cor:contfixe}, $\left \|u_t \right \|_{W_n^{\infty}(\R)}\tend{t\rightarrow t'}\left \|u_{t'}\right \|_{W_n^{\infty}(\R)}$, thus the right-hand side goes to zero proving the claim.
	
	For the second point, we notice that 
	$$\left \|\dfrac{\rho_{V_{\phi,t}}-\rho_{V_{\phi,t_0}}}{\delta t}-\rho_{V_{\phi,t_0}}u_{t_0}\right \|_{W_n^{\infty}(\R)}=\left \|(u_t-u_{t_0})\rho_{V_{\phi,t_0}}\right \|_{W_n^{\infty}(\R)}\leq2^n\|u_t-u_{t_0} \|_{W_n^{\infty}(\R)}\cdot\|\rho_{V_{\phi,t_0}} \|_{W_n^{\infty}(\R)}.$$
	Since the RHS goes to zero as $t\rightarrow t_0$, $u_{t_0}=-\phi+\displaystyle\int_\R\phi(y) \diff\mu_{V_{\phi,t_0}}(y)$ and $n$ is arbitrary, we conclude that, $x\in\R$, $t\mapsto\rho_{V_{\phi,t}}(x)$ is differentiable at every $t\in[0,1]$ of derivative $$\partial_t\rho_{V_{\phi,t}}^{(k)}(x)=-\left(\rho_{V_{\phi,t}}\phi\right)^{(k)} (x)+\rho_{V_{\phi,t}}^{(k)}(x)\int_\R\phi(y) \diff\mu_{V_{\phi,t}}(y).$$
	Since the above expression is again differentiable in $t$ (one deals with the integral by dominated convergence theorem with the domination $|\phi(x)\partial_t\rho_{V_{\phi,t}}(x)|\leq2\|\phi\|_\infty^2(1+\max_{s\in[0,1]}\|u_s\|_{\infty})\rho_{V_{\phi,t_0}}(x)$ for an arbitrary $t_0$), we conclude that for all $x\in\R$, $t\mapsto\rho_{V_{\phi,t}}(x)\in\mcal{C}^\infty(\R)$.
\end{proof}

\begin{corollary}
	[Convergence of moments]
	Let $h\in\N$, $t\in[0,1]$ by denoting $m_t(h)\defi \displaystyle\int_\R x^{h}\diff\mu_{V_{\phi,t}}(x)$, we have for all $t_0\in[0,1]$
	$$m_t(h)\underset{t\rightarrow t_0}{\rightarrow}m_{t_0}(h).$$
\end{corollary}

\begin{proof}
	For all $t\in[0,1]$, for all $x\in\R$, $\rho_{V_{\phi,t}}(x)\leq (1+\max_{s\in[0,1]}\|u_s\|_{\infty})\rho_{V_{\phi,t_0}}(x)$, hence by dominated convergence theorem $|m_{t}(h)-m_{t_0}(h)|\underset{t\rightarrow t'}{\rightarrow}0$.
\end{proof}

\section{Proof of Theorem \ref{thm:partitionfunction}}\label{section:freeenergy}
In this section, we prove, using Mehta's explicit formula, the expansion of $\log\mcal{Z}_N\left [V_G\right ]$. We then deduce Theorem \ref{thm:partitionfunction} from the latter and from Theorems \ref{thm:correlators} and \ref{thm:conteqdensity}.
\subsection{Expansion of the partition function for the Gaussian potential}

The asymptotic behaviour of $\mcal{Z}_N\left [V_G\right ]$ can be deduced from Mehta's formula \cite[17.6.7]{mehta2004random} \begin{equation}\label{partition,function gaussian}
	\mcal{Z}_N\left [V_G\right ]=(2\pi)^{N/2}\prod_{i=1}^N\dfrac{\Gamma\left (1+\dfrac{iP}{N}\right )}{\Gamma\left (1+\dfrac{P}{N}\right )}.
\end{equation}
This will allow us to use this formula in our interpolating integration formula to deduce the asymptotic expansion of $\log\mcal{Z}_N[V_\phi]$. From the previous equation, we can deduce the asymptotic behaviour of $\log	\mcal{Z}_N\left [V_G\right ]$. It is given by the following theorem:
\begin{theorem}\label{asymptotic free energy gaussian}
	There exists a sequence $(g_k)_{k\geq0}\in\R^N$, such that for all $K\geq0$,
	\begin{equation}\label{thm:asymptfctpartgauss}
		\boxed{	\dfrac{1}{N}\log\mcal{Z}_N\left [V_G\right ]=\sum_{i=0}^{K}\dfrac{g_i(P)}{N^{i}}+O\left(N^{-(K+1)}\right) }
	\end{equation}
	with
	\begin{equation}\label{eq:g1}
		g_1(P)\defi \gamma\dfrac{P}{2}+\dfrac{\log(1+P)}{2}+\dfrac{1}{2}\sum_{j\geq1}\left(\log\left(1+\dfrac{P+1}{j}\right)-\log\left(1+\dfrac{1}{j}\right)-\dfrac{P}{j}    \right).
	\end{equation}
	Above $\gamma$ denotes the Euler-Mascheroni constant.
\end{theorem}

\begin{proof}We first use \eqref{partition,function gaussian} to deduce 
	\begin{equation}\label{eq:logpartifuncgaussian}
		\log\mcal{Z}_N\left [V_G\right ]=\dfrac{N\log(2\pi)}{2}+\sum_{i=1}^N\Gamma\left (1+\dfrac{iP}{N}\right )-N\log\Gamma\left (1+\dfrac{P}{N}\right ).
	\end{equation}
	Let $K>0$, using the Taylor series expansion of $\log\Gamma$ around 1 (see \cite[8.342]{gradshteyn2014table}), one has
	\begin{equation}\label{taylor:loggamma}
		-N\log\Gamma\left (1+\dfrac{P}{N}\right )=\gamma P-\sum_{k=1}^{+\infty}\dfrac{\zeta(k+1)}{k+1}\dfrac{(-P)^{k+1}}{N^{k}}.
	\end{equation}
	where $\zeta$ denotes the Riemann $\zeta$ function. The second term in \eqref{eq:logpartifuncgaussian} can be estimated by using the Weierstrass product formula for $\dfrac{1}{\Gamma}$:
	\begin{equation}\label{key}
		\dfrac{1}{\Gamma(z)}=e^{\gamma z}z\prod_{j=1}^{+\infty}(1+\dfrac{z}{j})e^{-z/j}
	\end{equation}
	which is valid for any $z\in\C$. Hence we deduce that:
	\begin{multline}\label{otherterm}
		\sum_{k=1}^N\log\Gamma\left (1+\dfrac{kP}{N}\right )=-\sum_{k=1}^N\left (\gamma\Big(1+\dfrac{k P}{N}\Big)+\log\left (1+\dfrac{kP}{N}\right )-S_N(k) \right )
		\\=-\gamma N-\gamma\dfrac{(N+1)P}{2}-\sum_{k=1}^N\log\left (1+\dfrac{kP}{N}\right )+\sum_{k=1}^{N}S_N(k)
	\end{multline}where $S_N(k)\defi -\displaystyle\sum_{j=1}^{+\infty}\left[\log\left (1+\dfrac{1}{j}+\dfrac{kP}{Nj}\right )-\dfrac{1}{j}-\dfrac{kP}{Nj}\right] $.
	By the Euler-Maclaurin formula, we have the following identity for any $K>0$:
	\begin{equation}\label{eulermaclaurin1}
		\sum_{k=1}^N\log\left (1+\dfrac{kP}{N}\right )=\int_{0}^{N}f_N(t)\diff t+\sum_{k=1}^{K+2}\dfrac{B_k}{k!}\left (f_N^{(k-1)}(N)-f_N^{(k-1)}(0)\right )+R_{K+2}^{(N)}
	\end{equation}
	where $f_N(x)\defi \log \left (1+\dfrac{xP}{N}\right )$ and $B_k$ is the $k$-th Bernoulli number. The remainder $R_{K+2}^{(N)}$ is defined by $R_{K+2}^{(N)}\defi (-1)^{K+1}\displaystyle\int_{0}^{N}f_N^{(K+2)}(t)\dfrac{\widetilde{B}_{K+2}(t-\lfloor t \rfloor)}{(K+2)!}\diff t$, where $\widetilde{B}_{K+2}$ is the $(K+2)$-th Bernoulli polynomial. By using the following bound on Bernoulli polynomials,
	$$\forall x\in[0,1], \forall k>0,\hspace{1,5cm} |\widetilde{B}_k(x)|\leq2\dfrac{k!}{(2\pi)^k}\zeta(k)$$
	where $\zeta$ is the Riemann zeta function, $R_{K+2}^{(N)}$ can be controlled by the following inequalities:
	\begin{align*}
		|R_{K+2}^{(N)}|\leq\dfrac{2\zeta(K+2)}{(2\pi)^{K+2}}\int_{0}^{N}|f_N^{(K+2)}(t)|\diff t&=\dfrac{2\zeta(K+2)}{(2\pi)^{K+2}}\int_{0}^{N}\dfrac{P^{K+2}}{N^{K+2}}\dfrac{(K+1)!}{\left (1+\dfrac{P}{N}t\right )^{K+2}}\diff t
		\\&=\dfrac{2\zeta(K+2)(K+1)!}{(2\pi)^{K+2}}\dfrac{P^{K+1}}{N^{K+1}}\int_{0}^{P}\dfrac{\diff u}{\left (1+u\right )^{K+2}}=O\Big(N^{-{(K+1)}}\Big).
	\end{align*}Extracting the large-$N$ behaviour in (\ref{eulermaclaurin1}) leads to
	\begin{multline}\label{euleurmaclaurinNextracted}
		\sum_{k=1}^N\log\left (1+\dfrac{kP}{N}\right )=\dfrac{N}{P}\int_{0}^{P}\log(1+t)\diff t+B_1\log(1+P)\\+\sum_{k=2}^{K+2}\dfrac{B_k}{k!}(-1)^{k}(k-2)!\dfrac{P^{k-1}}{N^{k-1}}\left (\dfrac{1}{(1+P)^{k-1}}-1\right )+O\Big(N^{-(K+1)}\Big)
		\\=Nc_{-1}+\sum_{k=0}^{K}\dfrac{c_{k}}{N^{k}}+O\Big(N^{-(K+1)}\Big)
	\end{multline} 
	where $c_{-1}\defi (1+P^{-1})\log(1+P)-1$, $c_0\defi \dfrac{\log(1+P)}{2}$ and for all $k\in\llbracket1,K\rrbracket$,
	$$c_{k}=\dfrac{-B_{k+1}(-P)^{k}}{k(k+1)}\left (\dfrac{1}{(1+P)^{k}}-1\right ).$$ Also by Fubini's theorem, we get,
	
	\begin{align*}\label{Fubini}
		\sum_{k=1}^NS_N(k)&=-\sum_{j=1}^{+\infty}	\left\{\sum_{k=1}^N\log\left (1+\dfrac{1}{j}+\dfrac{kP}{Nj}\right)-\dfrac{1}{j}-\dfrac{kP}{Nj}\right\}
		\\&=\sum_{j=1}^{+\infty}\left\{\dfrac{N}{j}+\dfrac{(N+1)P}{2j}+\sum_{k=1}^Ng_{N,j}(k)\right\} 
	\end{align*}
	where $g_{N,j}(x)=\log\left (1+\dfrac{1}{j}+\dfrac{P}{jN}x\right )$. The first equality clearly shows that the RHS is a serie of general term bounded by $O\left (j^{-2}\right )$, so it converges and justifies the application of Fubini's theorem. Again by Euler-Maclaurin formula, we get:
	\begin{multline}
		\sum_{k=1}^Ng_{N,j}(k)=-\int_{0}^{N}\log\left (1+\dfrac{1}{j}+t\dfrac{P}{Nj}\right )\diff t-\dfrac{1}{2}\left [\log\left (1+\dfrac{P+1}{j}\right )-\log\left (1+\dfrac{1}{j}\right )\right ]\\-\sum_{k=2}^{K+2}\dfrac{(-1)^{k}B_k}{k(k-1)}\left\{\left (1+\dfrac{P+1}{j}\right )^{1-k}-\left (1+\dfrac{1}{j}\right )^{1-k}\right\}\left (\dfrac{P}{Nj}\right )^{k-1}+R_{K+2}^{(N)}(j)
	\end{multline}
	where again the new remainder $R_{K+2}^{(N)}(j)$ can be controlled \textit{via}\begin{multline*}
		|R_{K+2}^{(N)}(j)|\leq\dfrac{2\zeta(K+2)}{(2\pi)^{K+2}}\int_{0}^{N}\dfrac{(K+1)!\left(\dfrac{P}{Nj}\right)^{K+2} }{\left(1+\dfrac{1}{j}+t\dfrac{P}{jN}\right)^{K+2}}\diff t\\=\dfrac{2\zeta(K+2)}{(2\pi)^{K+2}}\left(\dfrac{P}{Nj}\right)^{K+1}\left\{\left (1+\dfrac{P+1}{j}\right )^{-(K+1)}-\left (1+\dfrac{1}{j}\right )^{-(K+1)}\right\}=O\Bigg (\dfrac{1}{(Nj)^{K+1}}\Bigg )
	\end{multline*} 
	where $O\left (\dfrac{1}{j^{K+1}}\right )$ depends on $K$ and $P$ but not on $N$. Hence we deduce that
	\begin{multline*}
		\sum_{k=1}^Ng_{N,j}(k)=-N\int_{0}^{1}\log\left (1+\dfrac{1+sP}{j}\right )\diff s -\dfrac{1}{2}\left [\log\left (1+\dfrac{P+1}{j}\right )-\log\left (1+\dfrac{1}{j}\right )\right ]\\+\sum_{k=1}^{K}\dfrac{(-P)^{k}B_{k+1}}{k(k+1)}\left\{\left (1+\dfrac{P+1}{j}\right )^{-k}-\left (1+\dfrac{1}{j}\right )^{-k}\right\}\dfrac{1}{(Nj)^k}+O\left ((jN)^{-(K+1)}\right ).
	\end{multline*}
	This leads to:
	\begin{multline}\label{}
		\sum_{k=1}^NS_N(k)=\sum_{j=1}^{+\infty}\left[\dfrac{N}{j}+\dfrac{(N+1)P}{2j}+\sum_{k=1}^Ng_{N,j}(k)\right] =\sum_{j=1}^{+\infty}\left(u_j^{(1)}N+\sum_{k=0}^{K}u_j^{(k)}N^{-k}\right) +O(N^{-(K+1)})\\=d_1N+\sum_{k=0}^{K}d_{-k}N^{-k} +O\Big(N^{-(K+1)}\Big)
	\end{multline}
	where for all $k=-1,0\dots K$, $\left( u_j^{(k)}\right )_{j>0}\in\ell^1(\N^*)$ and $d_{-k}\in\R$. This establishes the existence of the asymptotic expansion of $	\log	\mcal{Z}_N\left [V_G\right]$ up to $O\Big(N^{-(K+1)}\Big)$. Collecting the different terms, leads to the formula for $g_1(P)$.
\end{proof}

\subsection{Free energy of the model}
Only, in this subsection, since the parameter $P$ varies, we include the $P$-dependence of $\mcal{Z}_N[V]$ in the notation and write $\mcal{Z}_N^P[V]$ instead. The following result exhibits that there is a logarithmitc divergence of the free energy when $P\rightarrow0$ and hence a sort of phase transition.
\begin{theorem}[Free energy formula for Gaussian Potential]
	Let $P>0$, the free energy associated with the Gaussian potential is	
	\begin{equation}
		F(P)\defi \lim_{N\rightarrow\infty}N^{-1}\log\mcal{Z}_N^P\left [2PV_G\right ]=-\dfrac{1+P}{2}\log(2P)+\dfrac{\log(2\pi)}{2}+\int_{0}^{1}\log\Gamma\left(1+Px\right)\diff x
	\end{equation}
	
	As $P$ goes to $+\infty$, we have:
	\begin{equation}\label{free:energy formula asympto}
		F(P)=-P\left(\dfrac{3+\log 2}{2}\right) -\dfrac{1+\log 2}{2}+\log(2\pi)+\dfrac{\log P}{12 P}+O\left(P^{-1}\right) 
	\end{equation}
	
\end{theorem}
\begin{proof}By a change of variable, it holds that $\mcal{Z}_N^P\left [2PV_G\right ]=\left(\sqrt{2P}\right)^{-N-P(N-1)} \mcal{Z}_N^P\left[V_G\right]. $ Hence by Mehta's formula \eqref{partition,function gaussian},
	$$F(P)\defi \lim_{N\rightarrow\infty}N^{-1}\log\mcal{Z}_N^P\left [2PV_G\right ]=-\dfrac{1+P}{2}\log(2P)+\dfrac{\log(2\pi)}{2}+\int_{0}^{1}\log\Gamma\left(1+Px\right)\diff x.$$
	We can replace the last term by its asymptotic expansion so that
	
	$$\int_{0}^{1}\log\Gamma\left(1+Px\right)\diff x=\dfrac{(P+1)}{2}\log P-\dfrac{3P}{2}+\dfrac{\log(2\pi)-1}{2}+\dfrac{1}{12P}\log P+O\left(P^{-1}\right).$$
	We used the classic formula to conclude
	$$\log \Gamma(1+Px)=(1+Px)\log (1+Px)-1-Px-\dfrac{\log(1+Px)-\log(2\pi)}{2}+\dfrac{1}{12(1+Px)}+O\left(\dfrac{1}{(1+Px)^3}\right).$$
\end{proof}
\subsection{Interpolation with general potential}
We first establish the link between the 1-linear statistics and the partition function with general potential and the one with Gaussian potential.
\begin{lemma}\label{interpolation}
	Let $V_t(x)=tV(x)+(1-t)V_G(x)$ with $t\in[0,1]$. We have
	\begin{equation}\label{eq:interpolation}
		\log\dfrac{\mcal{Z}_N[V]}{	\mcal{Z}_N\left [V_G\right ]}=-N\int_{0}^{1}\Braket{V-V_G}_{\mu_N}^{V_t}\diff t
	\end{equation}
\end{lemma}

\begin{proof}
	By the fundamental theorem of calculus:
	$$\log\dfrac{\mcal{Z}_N[V]}{\mcal{Z}_N\left [V_G\right ]}=\int_{0}^{1}\partial_t\log \mcal{Z}_N[V_t]\diff t=-\int_{0}^{1}\diff t\int_{\R^N}p_N^{V_t}(\underline{x})\sum_{i=1}^{N}\partial_t V_t(x_i)\diff^N \underline{x}.$$
	Since $\displaystyle\sum_{i=1}^{N}\partial_t V_t(x_i)=N\displaystyle\int_\R \left[V(x)-V_G(x)\right] d\mu_N (x)$, where $\mu_N$ is the empricial measure associated to the external potential $V_t$, it concludes the proof.
\end{proof}

\begin{theorem}
	For all $\phi\in H^{\infty}(\R)$, there exists a sequence $(c_i)_{i\geq0}\in\R^{\N}$ depending on $\phi$ and $P$ such that for all $K\geq0$  
	$$\dfrac{1}{N}\log \mcal{Z}_N\left[V_{G,\phi}\right]=\sum_{i=0}^{K}\dfrac{c_i}{N^i}+O\left(N^{-(K+1)}\right) .$$
	The leading term $c_{0}$ is equal to the following expression:
	$$\int_\R V_{G,\phi}(x)\diff\mu_{V_{G,\phi}}(x)-P\iint_{\R^2}\log|x-y|\diff\mu_{V_{G,\phi}}(x)\diff\mu_{V_{G,\phi}}(y)+\int_{\R}\log\left(\dfrac{\diff\mu_{V_{G,\phi}}(x)}{\diff x}\right) \diff\mu_{V_{G,\phi}}(x).$$
	The subleading term $c_1$ can be written as
	\begin{multline}
		c_1\defi \gamma\dfrac{P}{2}+\dfrac{\log(1+P)}{2}+\dfrac{1}{2}\sum_{j\geq1}\left(\log\left(1+\dfrac{P+1}{j}\right)-\log\left(1+\dfrac{1}{j}\right)-\dfrac{P}{j}\right)
		\\-P\int_{0}^{1}\left[\Braket{\partial_1\widetilde{\Xi^{-1}}\phi}_{\mu_{V_{G,\phi,t}}}+\Braket{\Theta^{(2)} \circ\widetilde{\Xi_1^{-1}}\left[\partial_2\mcal{D}\circ\widetilde{\Xi^{-1}}\phi\right]  }_{\mu_{V_{G,\phi,t}}}\right] \diff t.
	\end{multline}
\end{theorem}

\begin{proof}
	By Lemma \ref{interpolation} and Theorem \ref{asymptotic free energy gaussian}, to establish the asymptotic expansion of $\mcal{Z}_N\left[V_{G,\phi}\right]$, it suffices to obtain the one for $\displaystyle\int_{0}^{1}\Braket{\phi}_{\mu_N}^{V_{G,\phi,t}}\diff t$. By Theorem \ref{asymptotic correlators}, we get, with $(d_i^{(1),t})$ the coefficients appearing in the expansion under the choice of potential $V_{G,\phi,t}$:
	$$\displaystyle\int_{0}^{1}\Braket{\phi}_{\mu_N}^{V_{G,\phi,t}}\diff t=\displaystyle\int_{0}^{1}\Braket{\phi}_{\mu_{V_{G,\phi,t}}}\diff t+\sum_{i=1}^{K}\dfrac{\displaystyle\int_0^1d_{i}^{(1),t}(\phi)\diff t}{N^i}+\displaystyle\int_{0}^{1}\left(\Braket{\phi}_{\Delta\mu_N}^{V_{G,\phi,t}}-\sum_{i=1}^{K}\dfrac{d_{i}^{(1),t}(\phi)}{N^i}\right)\diff t.$$
	Finally, we conclude that the last integral is a $O\left(N^{-(K+1)}\right)$ by \eqref{eq:controleraminedercorrelator} and the continuity of the map $t\mapsto C_{\mrm{rem}}(V_{G,\phi,t})$ obtained in Proposition \ref{prop:integrabilityremainder}. Furthermore, by collecting order $1$ for $\log\mcal{Z}_N\left[V_G\right]$ and $\displaystyle\int_0^1d_{i}^{(1),V_{G,\phi,t}}(\phi)\diff t$ for $i=1$ in Theorem \ref{asymptotic correlators}, we infer on the value of $c_{1}$.
\end{proof}

\section{Conclusion and open questions}

This work adapted the loop equations method to prove the existence of a $N^{-1}$ asymptotic expansion for the partition function, see Theorem \ref{thm:partitionfunction}, for a general class of potential. This class includes all potentials of the form $x^2+\phi$ where $\phi$ is a smooth bounded function. An immediate continuation of this result would be to extend it to more general confining potentials like $x^4$ for example. Another natural question would be to extend these ideas to more general interactions (Riesz, Coulomb, 2D-log in the spirit of \cite[Theorem 2]{serfaty2023gaussian}). Finally, our result shows that, contrary to \cite{BoG2}, there is no oscillatory terms in the expansion of $\log\mcal{Z}_N[V_{G,\phi,t}]$ since there is no multi-cut situation in the high temperature regime.

\appendix
{
	\section{Appendix: Lemmas and technical results}
	\label{app A}
	The first lemma in this appendix recalls several useful properties of the Hilbert transfrom used throughout this article.
	\begin{lemma}[Properties of the Hilbert transform]\ 
		\label{lemma:HilbertPropriétés}
		\begin{itemize}
			\item[i)] As a consequence, $\pi^{-1}\mathcal{H}$ is an isometry of $L^2(\R)$, and $\mathcal{H}$ satisfies on $L^2(\R)$ the identity $\mathcal{H}^2=-\pi^2I$.
			\item[ii)] Derivative: For any $f\in H^1(\R)$, $\mathcal{H}[f]$ is also $H^1(\R)$ and $\mcal{H}[f]'=\mcal{H}[f']$.
			\item[iii)] For all $p>1$, the Hilbert transform can be extended as a bounded operator $\mathcal{H}:L^p(\R)\to L^p(\R)$.
			\item[iv)] Skew-self adjointness: For any $f,g\in L^2(\R)$, $\Braket{\mcal{H}[f],g}_{L^2(\R)}=-\langle f,\mcal{H}[g] \rangle_{L^2(\R)}$.
			\item[v)] For all $\delta>0$, for all $f\in L^1(\R)$ such that $f'\in L^\infty(\R)$, $\|\mcal{H}[f]\|_{\infty}\leq(\delta^{-1}\|f\|_{1}+2\delta\|f'\|_\infty)$
		\end{itemize}
	\end{lemma}
	
	\begin{proof}
		We refer to \cite{Hilberttransforms} for the proofs of properties \textit{i)-iv)}. To prove \textit{v)}, let $f$ be such a function,
		$$|\mcal{H}[f](x)|
		\leq \lim_{\varepsilon\rightarrow0}\Big|\int_{\varepsilon\leq|x-y|
			\leq\delta}\dfrac{f(y)\diff y}{y-x}\Big|+\lim_{\varepsilon\rightarrow0}\Big|\int_{\delta\leq|x-y|\leq\varepsilon^{-1}}\dfrac{f(y)\diff y}{y-x}\Big|.$$
		The second term in the RHS can be bounded by $\delta^{-1}\|f\|_{L^1(\R)}$ while the first term verifies
		$$\lim_{\varepsilon\rightarrow0}\Big|\int_{\varepsilon\leq|x-y|\leq\delta}\dfrac{f(y)\diff y}{y-x}|\leq\lim_{\varepsilon\rightarrow0}\int_{\varepsilon\leq|x-y|\leq\delta}\Bigg|\dfrac{f(y)-f(x)}{y-x}\Bigg|\diff y+\lim_{\varepsilon\rightarrow0}\Bigg|f(x)\int_{\varepsilon\leq|x-y|\leq\delta}\dfrac{\diff y}{y-x}\Bigg|.$$
		The first term in the RHS can be bounded by $2\delta\|f'\|_{\infty}$ while the second is equal to $0$. This allows to conclude.
	\end{proof}
	
	The following result is proven in \cite[Lemma 2.3]{DwoMemin}.
	\begin{lemma}
		\label{lemme:boundedHilbert}
		Let $u\in L^2(\R)$ be such that $\int_\R u(t)\diff t$ exists and let $f:t\mapsto tu(t)\in H^1(\R)$ then
		$$\mcal{H}[u](x)\underset{|x|\rightarrow\infty}{\sim}\dfrac{-\displaystyle\int_\R u(t)\diff t}{x}.$$
		Moreover if  $\displaystyle\int_\R u(t)\diff t=0$, $\displaystyle\int_\R f(t)\diff t$ exists and $g:t\mapsto t^2u(t)\in H^1(\R)$, then $$\mcal{H}[u](x)\underset{|x|\rightarrow\infty}{\sim}\displaystyle\dfrac{-\displaystyle\int_\R tu(t)\diff t}{x^2}.$$
		As a consequence, we obtain that $\mcal{H}[\rho_V](x)\underset{|x|\rightarrow\infty}{\sim}-x^{-1}$ and the logarithmic potential $U^{\rho_V}$ is Lipschitz bounded, with bounded derivative $\mathcal{H}[\rho_V]$.
	\end{lemma}
	
	We now state and prove an inequality on the norm of the inverse of the master operator multiplied by the equlibrium density in Sobolev spaces. This inequality is crucial in the controls of Section \ref{section4}.

	\begin{lemma}\label{lem:majorationderivee nieme de rho Xi} Let $n\geq1$, and $h\in H^n(\R)$,
		$$\|\rho_V\widetilde{\Xi^{-1}}[h] \|_{H^{n}(\R)}\leq C_3(V,n)\|h\|_{H^{n}(\R)}.$$
		with a constant $C_3(V,n),$ only depending on $V$ and $n$. For the choice of potentiel $V=V_{\phi,t}$, for $\phi\in\mcal{C}^{\infty}(\R)$ with $\phi^{(k)}\in L^2(\R)$ for all $k\in\N$ and $t\in[0,1]$, $t\mapsto C_3(V_{\phi,t},n)$ is a continuous function.

		Moreover, for all $h\in H^n(\R)\cap W_n^{\infty}(\R)$,
		$$\|\rho_V\widetilde{\Xi^{-1}}[h] \|_{H^{n}(\R)}\leq C_5(V,n)\|h\|_{W_{n}^\infty(\R)}$$
		with a constant $C_5(V,n),$ only depending on $V$ and $n$. The function $t\mapsto C_5(V_{\phi,t},n)$ is also continuous.	
	\end{lemma}
	
	\begin{proof}Recall that $\theta=(\log\rho_V)'$, we first prove that for all $k\geq0$, $h\in\dfrac{1}{\rho_V}H^{k}(\R)$, there exists finite sets of indices $\mfrak{I}_{l,a}^k$, $\mfrak{J}^{k}_{l}$ and $\mfrak{K}_{l}^{k}$ independent of $V$ and polynomials 
		$p_{a,1,1}^k$ ,$p_{a,b,c,d}^{k}$, $q_{i,b,c}^{k}$ in $\theta,\dots,\theta^{(k-1)}$, with coefficients independent of $V$ and of degree at most $k$, such that for all $x\in\R$,
		\begin{multline}\label{eq:l-1enfonctiondehilbert}
			\widetilde{\Xi^{-1}}[h]^{(k)}=\sum_{i=0}^{k-1}p_{i,1,1}^{k}h^{(i)}+\sum_{i=0}^{k-1}\sum_{j\in\mfrak{I}^k_{2,i}}p_{i,j,2,1}^{k}\mcal{H}\left[\rho_Vp_{i,j,2,2}^kh^{(i)}\right]+\hdots
			\\+\sum_{i=0}^{k-1}\sum_{j\in\mfrak{I}^k_{k,i}}p_{i,j,k,1}^{k}\mcal{H}\Bigg[\rho_Vp_{i,j,k,2}^{k}\mcal{H}\Big[\rho_Vp_{i,j,k,3}^{k}\mcal{H}\Big[ \hdots\mcal{H}\Big[\rho_V p_{i,j,k,k}^{k}h^{(i)} \underbrace{\Big]\Big]\hdots\Bigg ]}_{k-1}+q^k_{1,1}\widetilde{\Xi^{-1}}[h] \\+\sum_{j\in\mfrak{J}_{2}^k}q^k_{j,2,1}\mcal{H}\left[\rho_V q^k_{j,2,2}\widetilde{\Xi^{-1}}[h] \right]+\hdots+\sum_{j\in\mfrak{J}_{k+1}^k}q^k_{j,k+1,1}\mcal{H}\Bigg[\rho_Vq^k_{j,k+1,2}\mcal{H}\Big[\hdots\mcal{H}\Big[\rho_Vq^k_{j,k+1,k+1}\widetilde{\Xi^{-1}}[h]\underbrace{\Big]\hdots\Bigg]}_{k}
			\\+\Bigg(r^k_{1,1}+\sum_{j\in\mfrak{K}^k_{2}}r_{j,2,1}^{k}\mcal{H}\left[\rho_Vr_{j,2,2}^k\right]+\hdots+\sum_{j\in\mfrak{K}^k_{k}}r^k_{j,k,1}\mcal{H}\Bigg[\rho_Vr^k_{j,k,2}\mcal{H}\Big[\hdots\mcal{H}\Big[\rho_Vr^k_{j,k,k}\underbrace{\Big]\hdots\Bigg]}_{k-1}\Bigg)
			\\\times\left(2P\int_\R\mcal{H}\left[\rho_V \widetilde{\Xi^{-1}}[h]\right]\diff\mu_V-\int_\R h\diff\mu_V\right).
		\end{multline}
		We prove it by induction, where for $n=1$ one just uses the definition of $\Xi$ for the initial case \textit{i.e.}
		$$		\left(\widetilde{\Xi^{-1}}[h]\right)'=h-\int_\R h(t)\diff\mu_V(t)-\dfrac{\rho_V'}{\rho_V}\left(\widetilde{\Xi^{-1}}[h]\right)'-2P\mcal{H}\left[\rho_V\widetilde{\Xi^{-1}}[h]\right]+2P\int_\R\mcal{H}\left[\rho_V\widetilde{\Xi^{-1}}[h]\right](t)\diff\mu_V(t).$$ For the induction step, we use a bootstrap argument. Suppose \eqref{eq:l-1enfonctiondehilbert} holds at rank $k$, then differentiate and replace $\left(\widetilde{\Xi^{-1}}[h]\right)'$ by the RHS of the above relation to show that \eqref{eq:l-1enfonctiondehilbert} holds at rank $k+1$.

		Now, by the Leibniz formula, for all $k\in\llbracket0,n\rrbracket$, it holds that
		$$\left(\rho_V\Xi^{-1}[h]\right) ^{(k)}=\sum_{i=0}^{k}\binom{k}{i}\rho_{V}^{(k-i)}\Xi^{-1}[h]^{(i)}.$$
		Furthermore by \eqref{eq:l-1enfonctiondehilbert}, by using successively that $\pi^{-1}\mcal{H}$ is an isometry of $L^2(\R)$, inequality \eqref{ineq:contXi} and Jensen's inequality, we obtain:
		$$\max_{0\leq k\leq n}\|\left(\rho_V\Xi^{-1}[h]\right) ^{(k)}(x)\|_{L^2(\R)}\leq C_3(V,n)\cdot\|h\|_{H^k(\R)}$$
		with $C_3(V,n)$ given by
		\begin{multline}\label{eq:def C_3(V,n)}
			C_3(V,n)\defi C(n)
			\\\times\max_{i\leq k\leq n}\Bigg\{i\max_{0\leq a<i}\left\|\rho_{V}^{(k-i)} p_{a,b,1,1}^{i}\right\|_{\infty}+i\pi\max_{0\leq a<i}|\mfrak{I}^i_{2,a}|\max_{b\in\mfrak{I}^i_{2,a}}\left (\left\|\rho_{V}^{(k-i)}p_{a,b,2,1}^{i}\right\|_{\infty}\|\rho_{V} p_{a,b,2,2}^i\|_{\infty}	\right)+\dots
			\\+i\pi^{i-1}\max_{0\leq a<i}|\mfrak{I}^i_{i,a}|\max_{b\in\mfrak{I}^i_{i,a}}\left(\left\|\rho_{V}^{(k-i)} p_{a,b,i,1}^{i}\right\|_{\infty}.\prod_{l=2}^{i}\|\rho_{V} p_{a,b,i,l}^i\|_{\infty}\right) +C_{\mcal{L}}\left \|\rho_{V}^{(k-i)} q_{1,1}^{i}\right\|_{\infty}\|\rho_V^{1/2}\|_\infty
			\\+C_{\mcal{L}}\pi|\mfrak{J}^i_{2}|\max_{b\in\mfrak{J}^i_{2}}\left \|\rho_{V}^{(k-i)} q_{b,2,1}^{i}\right\|_{\infty}\|\rho_V^{1/2}q_{b,2,2}^i\|_\infty\|\rho_V^{1/2}\|_\infty+\dots
			\\+C_{\mcal{L}}\pi^{i}|\mfrak{J}^i_{i+1}|\max_{b\in\mfrak{J}^i_{i+1}}\left \|\rho_{V}^{(k-i)} q_{b,i+1,1}^{i}\right\|_{\infty}\prod_{l=2}^{i}\|\rho_V q_{b,i+1,l}^{i}\|_{\infty}\|\rho_V^{1/2} q_{b,i+1,i+1}^{i}\|_{\infty}\|\rho_V^{1/2}\|_{\infty}
			\\+2\|\rho_V\|_{\infty}\left(1+2P\pi C_{\mcal{L}}\right) \Big[\|r_{1,1}^{i}\dfrac{\rho_{V}^{(k-i)}}{\rho_V}\sqrt{\rho_V}\|_{\infty}+\pi|\mfrak{K}^i_{2}|\max_{b\in\mfrak{K}^i_{2}}\|\rho_V^{(k-i)}r_{b,2,1}^{i}\|_{\infty}\|\sqrt{\rho_V}r_{b,2,1}^{i}\|_{\infty}+\dots
			\\+\pi^{i-1}|\mfrak{K}^i_{i}|\max_{b\in\mfrak{K}^i_{i}}\|\rho_V^{(k-i)}r_{b,i,1}^{i}\|_{\infty}\prod_{l=2}^{i-1}\|\rho_V r_{b,i,l}^{i}\|_{\infty}\|\sqrt{\rho_V}r_{b,i,i}^{i}\|_{\infty}\Big] \Bigg\}
		\end{multline}
		For the second inequality, if $h\in H^n(\R)\cap W_n^{\infty}(\R)$, we use the fact the same inequalities but we use the following integrals at the end $\displaystyle\int_\R \left(p_{a,b,i,l}^{i} h^{(a)} \rho_V(t)\right)^2 \diff t\leq \|h\|_{W_a^\infty(\R)}^2\cdot\int_\R (p_{a,b,i,l}^{i}\rho_{V_{t}})^2\diff t$. This leads to
		$$\|\rho_V\widetilde{\Xi^{-1}}[h] \|_{H^{n}(\R)}\leq C_5(V,n)\|h\|_{W_{n}^\infty(\R)}$$
		with $C_5(V,n)$ given by
		\begin{multline}\label{eq:def C_5(V,n)}
			C_5(V,n)\defi C(n)
			\\\times\max_{i\leq k\leq n}\Bigg\{i\max_{0\leq a<i}\left\|\dfrac{\rho_{V}^{(k-i)}}{\rho_{V}} p_{a,b,1,1}^{i}\sqrt{\rho_V}\right\|_{\infty}+i\pi\max_{0\leq a<i}|\mfrak{I}^i_{2,a}|\max_{b\in\mfrak{I}^i_{2,a}}\left (\left\|\rho_{V}^{(k-i)}p_{a,b,2,1}^{i}\right\|_{\infty}\|\sqrt{\rho_{V}} p_{a,b,2,2}^i\|_{\infty}	\right)+\dots
			\\+i\pi^{i-1}\max_{0\leq a<i}|\mfrak{I}^i_{i,a}|\max_{b\in\mfrak{I}^i_{i,a}}\left(\left\|\rho_{V}^{(k-i)} p_{a,i,1}^{i}\right\|_{\infty}.\prod_{l=2}^{i-1}\|\rho_{V} p_{a,i,l}^i\|_{\infty}\|\sqrt{\rho_{V} }p_{a,i,i}^i\|_{\infty}\right)
			\\ +C_{\mcal{L}}\left \|\rho_{V}^{(k-i)} q_{1,1}^{i}\right\|_{\infty}\|\rho_V\|_\infty^{1/2}+C_{\mcal{L}}\pi|\mfrak{J}^i_{2}|\max_{b\in\mfrak{J}^i_{2}}\left \|\rho_{V}^{(k-i)} q_{b,2,1}^{i}\right\|_{\infty}\|\rho_V^{1/2}q_{b,2,2}^i\|_\infty+\dots
			\\+C_{\mcal{L}}\pi^{i}|\mfrak{J}^i_{i+1}|\max_{b\in\mfrak{J}^i_{i+1}}\left \|\rho_{V}^{(k-i)} q_{b,i+1,1}^{i}\right\|_{\infty}\prod_{l=2}^{i}\|\rho_V q_{b,i+1,l}^{i}\|_{\infty}\|\rho_V^{1/2} q_{b,i+1,i+1}^{i}\|_{\infty}
			\\+2\left(1+2P\pi\|\rho_V\|_\infty^{1/2}C_{\mcal{L}}\right) \Big[\|r_{1,1}^{i}\dfrac{\rho_{V}^{(k-i)}}{\rho_{V}}\sqrt{\rho_V}\|_{L^2(\R)}+\pi|\mfrak{K}^i_{2}|\max_{b\in\mfrak{K}^i_{2}}\|\rho_V^{(k-i)}r_{b,2,1}^{i}\|_{\infty}\|\sqrt{\rho_V}r_{b,2,1}^{i}\|_{\infty}+\dots
			\\+\pi^{i-1}|\mfrak{K}^i_{i}|\max_{b\in\mfrak{K}^i_{i}}\|\rho_V^{(k-i)}r_{b,i,1}^{i}\|_{\infty}\prod_{l=2}^{i-1}\|\rho_V r_{b,i,l}^{i}\|_{\infty}\|\sqrt{\rho_V}r_{b,i,i}^{i}\|_{\infty}\Big] \Bigg\}.
		\end{multline}
		The fact that $t\mapsto C_i(V_{\phi,t},n)$ for $i=3$ or $5$ is shown in Appendix \ref{appB}.
	\end{proof}
	
	\begin{remark}
		With $\theta=\dfrac{\rho_V'}{\rho_V}$, $g\defi \widetilde{\Xi^{-1}}[f]$ and $c\defi \displaystyle\int_\R\left( 2P\mcal{H}\left[\rho_Vg\right]-f\right) \diff\mu_V$ we have
		\begin{align*}
			g'&=f-\theta g-2P\mcal{H}\left[\rho_Vg\right] +c
			\\g''&=-\theta f+f'-2P\mcal{H}\left[\rho_Vf\right]+\left(\theta^2-\theta'\right)g+2P\theta\mcal{H}\left[\rho_Vg\right] +4P^2\mcal{H}\left[\rho_V\mcal{H}\left[\rho_V g\right]\right]
			\\&+\left(-\theta-2P\mcal{H}[\rho_V]\right)c.
			\\g^{(3)}&=\left(\theta^2-2\theta'\right)f-\theta f'+f''+2P\theta\mcal{H}\left[\rho_Vf\right]-2P\mcal{H}\left[\rho_V\theta f\right]-2P\mcal{H}\left[\rho_Vf'\right]+4P^2\mcal{H}\big[\rho_V\mcal{H}\left[\rho_V f\right]\Big]
			\\&+\left[\left(\theta^2-\theta'\right)'-\left(\theta^3-\theta\theta'\right)\right]g +2P\left(\theta^2-2\theta'\right)\mcal{H}[\rho_Vg]-4P^2\theta\mcal{H}\left[\rho_V\mcal{H}\left[\rho_V g\right]\right ]+4P^2\mcal{H}\left[\rho_V\theta\mcal{H}\left[\rho_V g\right]\right]
			\\&+4P^2\mcal{H}\left[\rho_V\mcal{H}\left[\rho_V\theta g\right]\right]-8P^3\mcal{H}\Big[\rho_V\mcal{H}\big[\rho_V\mcal{H}\left[\rho_V g \right]\big]\Big]
			\\&+\left((\theta^2-2\theta')+2P\theta\mcal{H}[\rho_V]-2P\mcal{H}[\rho_V\theta]+4P^2\mcal{H}\Big [\rho_V\mcal{H}[\rho_V]\right] \Big)
		\end{align*}
	\end{remark}
	\section{Integrability of the constants}\label{appB}
	\subsection{Parameter continuity of norms of certain functions}
	In this appendix, we work with $V=V_{G,\phi,t}$, $t\in[0,1]$ and $\phi\in H^{\infty}(\R)$.
	We will show, that the constant $C_i(V_{G,\phi,t})$ appearing in our problem, see Theorem \ref{thm: cont inverse master}, Theorem \ref{thm: cont inverse master linfini} and Theorem \ref{thm:continversemaster thetaxi L2} will be continuous in $t$ hence integrable on $[0,1]$. By using the main result from Section \ref{section7}, \textit{i.e.} Theorem \ref{thm:locallimittheorem} stating that the map $t\mapsto\|\rho_{t}\|_{W_n^{\infty}(\R)}$ (setting $\rho_{t}\defi\rho_{V_{G,\phi,t}}$ for this appendix) is continuous for all $n\in\N$, we will show the following lemma.
	
	\begin{lemma}\label{lem:leshilbertssontbornéesent}
		Let $t,t_0\in[0,1]$, for all $n\in\N$,
		$$\left \|\mcal{H}[\rho_{t}-\rho_{t_0}]\right \|_{W_n^{\infty}(\R)}\underset{t\rightarrow t_0}{\longrightarrow}0.$$
	\end{lemma}
	\begin{proof}
		We prove it by induction and use Lemma \ref{lemma:HilbertPropriétés} and Theorem \ref{thm:locallimittheorem}. For $n=0$, we know that there exists $C>0$, such that:
		$$\left \|\mcal{H}[\rho_{t}-\rho_{t_0}]\right \|_{\infty}\leq C\left(\|\rho_{t}-\rho_{t_0}\right\|_{L^1(\R)}+\|\rho_{t}'-\rho_{t_0}'\|_{\infty}).$$
		By Scheffé's lemma, the $L^1$ norm goes to zero and by Theorem \ref{thm:locallimittheorem} goes also to zero as $t$ goes to $t_0$. Now suppose that $\left \|\mcal{H}[\rho_{t}-\rho_{t_0}]\right \|_{W_n^{\infty}(\R)}\underset{t\rightarrow t_0}{\longrightarrow}0$ for some $n\geq0$. We have that:
		$$\left \|\mcal{H}[\rho_{t}^{(n+1)}-\rho_{t_0}^{(n+1)}]\right \|_{\infty}\leq C\left(\left \|\rho_{t}^{(n+1)}-\rho_{t_0}^{(n+1)}\right\|_{L^1(\R)}+\left \|\rho_{t}^{(n+2)}-\rho_{t_0}^{(n+2)}\right \|_{\infty}\right ).$$
		The second term on the RHS goes to zero by Theorem \ref{thm:locallimittheorem}. For the first term, since for all $k\in\N$, $\rho_{t_0}^{(k)}\in L^1(\R)$ and we have the following domination by Leibniz formula:
		$$|\rho_{t}^{(n+1)}(x)|\leq(1+\max_{s\in[0,1]}\|u_s\|_{W_{n+1}^\infty(\R)})\cdot\sum_{k=0}^{n+1}|\rho_{t_0}^{(k)}(x)|,$$
		we can conclude by the dominated convergence theorem that the $L^1$-norm goes to zero. The fact that $u_s$ is uniformly bounded comes from Corollary \ref{cor:contfixe}.
	\end{proof}
	
	Secondly, we set  for all $t\in[0,1]$, $M_{V_{G,\phi,t}}$  (see Lemma \ref{lem:M_V}) equal to any value $M$ such that: $$M>\max\left(1+\|\phi'\|_{\infty}+2P\max_{t\in[0,1]}\|\mcal{H}[\rho_{t}]\|_{\infty},2\left(\|\phi'\|_{\infty}+2P\max_{t\in[0,1]}\|\mcal{H}[\rho_{t}]\|_{\infty} \right)\right).$$It is well-defined because of Lemma \ref{lem:leshilbertssontbornéesent}. The following Lemma will be useful to show some controls on different quantities below. We recall that $\theta_t=(\log\rho_t)'$ and $\alpha_t=\dfrac{1}{\theta_t}$.
	
	\begin{lemma}\label{lem:M}
		For all $t\in[0,1]$, for all $i,j\in\N$, for all $|x|\geq  M$, 
		$$\bigg|\theta_t(x)\bigg|\geq 1 \hspace{0,5cm}\text{ and }\hspace{0,5cm}\Big|\alpha_t^{(i)}(x)\Big|^{j}\leq \delta_{i,0}\dfrac{C_{0,j}}{|x|^j}+\dfrac{C_{i,j}}{|x|^{2j}}$$
		for constants $C_{i,j}>0$  independent of $t$.
	\end{lemma}
	\begin{proof}
		Let $x\in\R$, $-\theta_t(x)=x+t\phi'(x)+2P\mcal{H}[\rho_{t}](x).$ Thus if $|x|\geq M$,
		$$\bigg|\dfrac{\rho_{t}'}{\rho_{t}}(x)\bigg|\geq1+\left (\|\phi'\|_{\infty}-t|\phi'(x)|\right )+2P\left (\max_{s\in[0,1]}\|\mcal{H}[\rho_{s}]\|_{\infty}-|\mcal{H}[\rho_{t}](x)|\right )\geq1.$$
		For the second point, one notices by differentiation and \eqref{deriv1} that there exists polynomials $\msf{P}_k$ with coefficients independent of $t$ such that
		\begin{equation}\label{lemb.2 preuve}
			\alpha_t^{(i)} (x)=\sum_{k=1}^{i}\dfrac{\msf{P}_k\left(t\phi'(x),\dots,t\phi^{(i+1)}(x),\mcal{H}[\rho_{t}](x),\dots,\mcal{H}\left [\rho_{t}^{(i)}\right ](x)\right) }{\left(x+t\phi'(x)+2P\mcal{H}[\rho_{t}](x)\right)^{k+1} }.
		\end{equation}
		Furthermore, since $|x|\geq2\left(\|\phi'\|_{\infty}+2P\displaystyle\max_{t\in[0,1]}\|\mcal{H}[\rho_{t}]\|_{\infty} \right)$, we have:
		\begin{equation*}
			\Big|x+t\phi'(x)+2P\mcal{H}[\rho_{t}](x)\Big|\geq \dfrac{|x|}{2}+\left(\dfrac{|x|}{2}-\|\phi'\|_{\infty}-2P\max_{s\in[0,1]}\|\mcal{H}[\rho_s]\|_{\infty}\right)\geq\dfrac{|x|}{2}.
		\end{equation*}
		Finally, the whole dependence in $t$ and $x$ of the numerator in \eqref{lemb.2 preuve} are in the entries which are bounded uniformly in $t$ and $x$, we can conclude that each numerator in the sum is bounded by a constant $C_k>0$. We can conclude that
		$$\Big|\left( \dfrac{1}{\alpha_t}\right) ^{(i)}(x)\Big|\leq4i\max_{k\leq i}C_k|x|^{-2},$$
		raising to the power $j$ leads to the conclusion.
	\end{proof}
	
	\begin{lemma}\label{lemappB:invrho}
		The map $t\mapsto\|\rho_{t}^{-1}\|_{L^\infty([-M,M])}$ is continuous.
	\end{lemma}
	\begin{proof}
		Let $x\in\R$, $t_0\in[0,1]$,
		\begin{equation}\label{eq:lemb3}
			|\rho_{t}(x)^{-1}-\rho_{t_0}(x)^{-1}|=\left|\dfrac{\delta t u_t(x)\rho_{t_0}(x)}{\rho_{t_0}(x)\rho_{t}(x)}\right| \leq|\delta t|\dfrac{\|u_t\|_{\infty}\|\rho_{t_0}^{-1}\|_{L^\infty([-M,M])}}{\left(1-|\delta t|\|u_t\|_{\infty}\right)}
		\end{equation}
		where $\delta t=t-t_0$ and $u_t$ is defined in \eqref{def:u_t} and is, by Corollary \ref{cor:contfixe}, continuous with respect to $t$. Taking the supremum over $x\in[-M,M]$ and let $t$ goes to $t_0$ establishes the result.
	\end{proof}
	
	Now it remains to bound the $L^2$ or $L^\infty$ norms of the functions $t\mapsto\mfrak{f}^{(i)}$ and $t\mapsto\mcal{I}_i$ used in the proofs of Theorems \ref{thm: cont inverse master}, \ref{thm: cont inverse master linfini} and \ref{thm:continversemaster thetaxi L2}. 
	
	\begin{lemma}\label{lem:I_a}
		The map $t\mapsto\|\mcal{I}_i^{t}\|_{\infty}$ for all $i\in\{1,2\}$ is continuous where $\mcal{I}_i^{t}\defi\mcal{I}_i$ is defined in \eqref{eq:defI_a} under the choice of potential $V_{G,\phi,t}$.
	\end{lemma}
	
	\begin{proof}The fact that these maps are well-defined can be seen by a tail estimate. To show the continuity, let $x>0$, $t,t_0\in[0,1]$, by the mean-value theorem and with $u_t$ defined in \eqref{def:u_t}, we get
		\begin{multline*}
			|\mcal{I}_1^{t}(x)-\mcal{I}_1^{t_0}(x)|=\Bigg|\dfrac{1}{\rho_{t}(x)}\int_x^{+\infty}\rho_{t}(s)\diff s -\dfrac{1}{\rho_{t_0}(x)}\int_x^{+\infty}\rho_{t_0}(s)\diff s \Bigg|
			\\\leq \dfrac{1}{\rho_{t_0}(x)}\int_x^{+\infty}\left| \dfrac{1+\delta tu_t(s)}{1+\delta tu_t(x)}-1\right|\rho_{t_0}(s)\diff s 
			\\\leq\dfrac{|\delta t|\|u_t'\|_{\infty}}{(1-|\delta t|\|u_t\|_{\infty})}\dfrac{1}{\rho_{t_0}(x)}\int_x^{+\infty}(s-x)\rho_{t_0}(s)\diff s .
		\end{multline*}
		One thus concludes by the fact that $\delta t$ goes to zero and that the following maps $$x\in]0,+\infty[\mapsto\dfrac{1}{\rho_{t_0}(x)}\displaystyle\int_x^{+\infty}(s-x)\rho_{t_0}(s)\diff s \hspace{1cm} \text{and}\hspace{1cm} t\mapsto\|u_t\|_{W_1^{\infty}(\R)}\hspace{1cm}\hspace{1cm}\hspace{1cm}$$ are bounded. For the second map, this follows from Corollary \ref{cor:contfixe}. For the first one, while it is obviously bounded on a compact, the boundedness at infinity can be proven by an integration by parts:
		\begin{multline*}
			\dfrac{1}{\rho_{t_0}(x)}\int_x^{+\infty}(s-x)\rho_{t_0}(s)\diff s =\dfrac{1}{\rho_{t_0}(x)}\left [\dfrac{\rho_{t_0}}{\rho_{t_0}'}(s)\rho_{t_0}(s)(s-x)\right ]_{x}^{+\infty}
			\\-\dfrac{1}{\rho_{t_0}(x)}\int_{x}^{+\infty}\left[\dfrac{\rho_{t_0}}{\rho_{t_0}'}(s)+\left(\dfrac{\rho_{t_0}}{\rho_{t_0}'}\right)'(s)(s-x)\right]\rho_{t_0}(s)\diff s .
		\end{multline*}
		
		The first term in the right hand side is zero, while by assumption \textit{\ref{assumption5}}, the last term is bounded by a tail-estimate shown by using Lemma \ref{lem:M}.
		Doing the same thing over $]-\infty,0]$ establishes that $t\mapsto\|\mcal{I}_1^{t}\|_{\infty}$ is continuous.
		
		Just as before, we get by the mean-value theorem,
		\begin{multline}
			|\mcal{I}_2^{t}(x)^2-\mcal{I}_2^{t_0}(x)^2|\leq \dfrac{1}{\rho_{t_0}^{2}(x)}\int_x^{+\infty}\left| \dfrac{(1+\delta tu_t(s))^2}{(1+\delta tu_t(x))^2}-1\right|\rho_{t_0}(s)^2\diff s
			\\\leq\dfrac{|\delta t|\|u_t'\|_{\infty}(2+\|u_t\|_\infty)}{(1-|\delta t|\|u_t\|_{\infty})^2}\dfrac{1}{\rho_{t_0}(x)^2}\int_x^{+\infty}(s-x)\rho_{t_0}(s)^2\diff s.
		\end{multline}
		We conclude by showing that $x\in[0,+\infty[\mapsto\dfrac{1}{\rho_{t_0}(x)^2}\displaystyle\int_x^{+\infty}(s-x)\rho_{t_0}(s)^2\diff s$ is bounded which can again be proven by the same integration by parts and by doing the exact same thing on $]-\infty,0]$. Therefore by the fact that
		$$0\leq\bigg|\|\mcal{I}_2^{t}\|_{\infty}^2-\|\mcal{I}_2^{t_0}\|^2_\infty\bigg|\leq\|\left(\mcal{I}_2^{t}\right) ^2-\left(\mcal{I}_2^{t_0}\right) ^2\|_{\infty}\tend{t\rightarrow t_0}0$$
		we conclude that $t\mapsto\|\mcal{I}_2^{t}\|_{\infty}$ is continuous.
	\end{proof}
	Next, we show that any polynomial in $\theta_t=(\log\rho_t)'$ and its derivatives yield a continuous dependence in $t$.
	\begin{lemma}\label{lemappB:norminfPoly}
		Let $\msf{Q}$ a polynomial in $\theta_t,\dots,\theta_t^{(k)}$ for some $k\geq0$ with coefficients independent of $t$, let $l\in\N$ then then the two following maps are continuous:
		\begin{enumerate}
			\item[(i)]$t\mapsto\|\msf{Q}\big(\theta_t,\dots,\theta_t^{(k)}\big)\|_{L^\infty([-M,M])}$
			\item[(ii)]$t\mapsto\|\sqrt{\rho_{t}^{(l)}}\msf{Q}\big(\theta_t,\dots,\theta_t^{(k)}\big)\|_{L^\infty(\R)}$
			\item[(iii)]$t\mapsto\|\rho_{t}^{(l)}\msf{Q}\big(\theta_t,\dots,\theta_t^{(k)}\big)\|_{L^\infty(\R)}.$	\end{enumerate}
	\end{lemma}
	
	\begin{proof}
		Proving these continuity results for any monomial in those variables is enough. Furthermore, since by continuity $x\in[-M,M]\mapsto\theta_t^{(i)}(x)$ is uniformly bounded in $t$ for all $i\leq k$, thus this monomial in $(\theta_t^{(i)})_{0\leq i \leq k}$ converges uniformly to the monomial in $(\theta_{t_0}^{(i)})_{0\leq i \leq k}$ as $t$ goes to $t_0$. It comes from the fact that the product of two bounded, uniformly converging sequences of functions converges to the product of the limits and that for all $i>0$, $\theta_t^{(i)}(x)-\theta_{t_0}^{(i)}(x)=-\delta t\phi^{(i+1)}(x)-2P\mcal{H}\left[\rho_{t}^{(i)}-\rho_{t_0}^{(i)}\right]$. The latter, when taking the supremum over $x\in[-M,M]$, goes to zero by Lemma \ref{lem:leshilbertssontbornéesent}. This establishes \textit{(i)}.
		Furthermore, notice that \textit{(iii)} implies \textit{(ii)} since $$\|\sqrt{\rho_{t}^{(l)}}\msf{Q}\left(\theta_t,\dots,\theta^{(k)}\right)\|_{L^\infty(\R)}=\|\rho_{t}^{(l)}\msf{Q}\left(\theta_t,\dots,\theta_t^{(k)}\right)^2\|_{L^\infty(\R)}^{1/2}$$ and $\msf{Q}$ is arbitrary so we only prove \textit{(iii)}. Moreover since $\rho_{t}^{(l)}=\exp(\log\rho_{t})^{(l)}$ can be written, by Faà di Bruno's formula,  as $\msf{R}\left(\theta_t,\dots,\theta_t^{(l)}\right)\rho_{t}$  where $\msf{R}$ is a polynomial with coefficients independent of $t$, it suffices to prove the result for $l=0$.
		
		For all $i\in\N$, we have $\theta_t^{(i)}(x)=-\delta_{i,0}x- t\phi^{(i+1)}(x)-2P\mcal{H}\left[\rho_{t}^{(i)}\right]$. Noticing, by Leibniz formula and the mean value theorem, that for all $j\in\N$ and for all $0<h<1$:
		\begin{multline*}
			\bigg|x^{j}\rho_{t}(x)^h-x^j\rho_{t_0}(x)^h\bigg|\leq\bigg|x^j\rho_{t_0}(x)^h\bigg|.\Big|\left (1+\delta tu_t(x)\right )^h-1\Big|
			\\\leq\dfrac{h|\delta t| \max_{s\in[0,1]}\|u_s\|_{W_l^{\infty}(\R)}}{\left(1-{|\delta t| \max_{s\in[0,1]}\|u_s\|_{W_l^{\infty}(\R)}}\right)^{1-h}}  \Big\|x\mapsto x^{j}\rho_{t_0}(x)^h\Big\|_{\infty}
		\end{multline*}
		where the existence of the max is justified by Corollary \ref{cor:contfixe}. Taking the supremum over $x\in\R$ and let $t$ goes to $t_0$ shows that $t\mapsto\left(x\mapsto x^j\rho_{t}(x)^h\right)\in L^\infty(\R)$ is continuous. By boundedness and continuity with respect to the $t$ parameter of $t\mapsto\phi^{(i+1)}(x)+2P\mcal{H}\left [\rho_{t}^{(i)}\right ]$ by Lemma \ref{lem:leshilbertssontbornéesent}, we deduce that for all $i\in\N$ and $\alpha>0$,
		$$\|\theta_t^{(i)}\rho_{t}(x)^{h}-\theta_{t_0}^{(i)}\rho_{t_0}(x)^{h}\|_{\infty}\tend{t\rightarrow t_0}0.$$
		From this last uniform convergence result, we show that by taking a monomial $\prod_{i=0}^{k}\left(\theta^{(i)}_t\right)^{l_i}$ such that $\sum_{i=0}^kl_i=m$, we deduce that $\rho_{t}\prod_{i=0}^{k}\left(\theta^{(i)}_t\right)^{l_i}=\prod_{i=0}^{k}\left(\theta^{(i)}_t\sqrt[m]{\rho_{t}}\right)^{l_i}$, as a product of bounded, uniformly converging $t$-sequences of functions, it converges uniformly. This concludes the proof.
	\end{proof}
	To state the following lemma, we recall that $\alpha_t= \dfrac{\rho_{t}}{\rho_{t}'}$.
	
	\begin{lemma}\label{lem:b6}[Continuity of $L^{\infty}$ norms]
		For all $j\in\{1,2,3,5,6\}$, $t\mapsto\|\mfrak{f}^{(j),t}_{l,\bm{i}}\|_{L^\infty([-M,M]^c)}$ is continuous where we recall that for $\bm{i}=(i_1,i_2)$:
		\begin{align*}
			&\mfrak{f}^{(1),t}_{l,\bm{i}}:x\mapsto Q_{i_1}^{l}\left(\theta_t,\dots,\theta_t^{(i_1)}\right)(x)\alpha_t(x) P_{i_2}^{l-i_1}(\alpha_t,\dots,\alpha_t^{(i_2)})(x),
			\\&\mfrak{f}_{l,\bm{i}}^{(2),t}:x\mapsto  \dfrac{Q_{i_1}^{l}\left(\theta_t,\dots,\theta_t^{(i_1)}\right)(x)}{\rho_{t}(x)}\cdot\int\limits_{x}^{\mrm{sgn}(x)\infty}\left[\alpha_t P_{i_2}^{l-i_1}\left (\alpha_t,\dots,\alpha_t^{(i_2)}\right)  \right] '(y)\rho_{t}(y)\diff y,
			\\&\mfrak{f}^{(3),t}_{l,\bm{i}}:x\mapsto \dfrac{Q_{i_1}^{l}\left(\theta_t,\dots,\theta_t^{(i_1)}\right)(x)}{\rho_{t}(x)}\cdot\Bigg(\int\limits_{x}^{\mrm{sgn}(x)\infty}\alpha_t(y)^2P_{i_2}^{l-i_1}\left (\alpha_t,\dots,\alpha_t^{(i_2)}\right )(y)^2\rho_{t}(y)^2\diff y\Bigg)^{1/2},
			\\&\mfrak{f}^{(5),t}_{l,\bm{i}}:x\mapsto \dfrac{Q_{i_1}^{l}\left(\theta_t,\dots,\theta_t^{(i_1)}\right)(x)}{\rho_{t}(x)}\cdot\int\limits_{x}^{\mrm{sgn}(x)\infty}\Big|P_{i_2}^{l-i_1}(\alpha_t,\dots,\alpha_t^{(i_2)})(y)\Big|\rho_{t}(y)\diff y,
			\\&\mfrak{f}^{(6),t}_{l,\bm{i}}:x\mapsto\dfrac{\Big|Q_{i_1}^{l}\left(\theta_t,\dots,\theta_t^{(i_1)}\right)(x)\Big|}{\rho_{t}(x)}\cdot\Bigg(\int\limits_{x}^{\mrm{sgn}(x)\infty}\Big|P_{i_2}^{l-i_1}(\alpha_t,\dots,\alpha_t^{(i_2)})(y)\Big|^2\rho_{t}(y)^2\diff y\Bigg)^{1/2}.
		\end{align*}
	\end{lemma}
	
	\begin{proof}
		First, one can check that, from Lemma \ref{lem:gammaetalpha}, $\mfrak{f}^{(1),t}_{l,\bm{i}}(x)=\underset{|x|\rightarrow\infty}{O}(x^{-1})$ for all $t\in[0,1]$. Recalling the definition of $\alpha_t(x)$, there exists $n>0$ and a polynomial expression $\msf{P}$ with coefficients independent of $t$ such that
		$$\mfrak{f}^{(1),t}_{l,\bm{i}}(x)=\dfrac{\msf{P}\left(x^{-1},t\phi',\dots,t\phi^{(k)},P\mcal{H}[\rho_{t}],\dots,P\mcal{H}[\rho_{t}^{(k-1)}]\right) }{(1+t\phi'(x)x^{-1}+2P\mcal{H}[\rho_{t}](x)x^{-1})^{n}}.$$
		where the numerator is a $\underset{|x|\rightarrow\infty}{O}(x^{-1})$.  Since, $\mfrak{f}^{(1),t}_{l,\bm{i}}$ is a bounded rational function such that the denominator is bounded from below uniformly in $t$ (see Lemma \ref{lem:M}) and such that both the numerator and denominator converges uniformly, we conclude from this closed form, that $\mfrak{f}^{(1),t}_{l,\bm{i}}$ converges uniformly to $\mfrak{f}^{(1),t_0}_{l,\bm{i}}$ when $t$ goes to $t_0$ on $[-M,M]^c$. Thus, $t\mapsto\|\mfrak{f}^{(1),t}_{l,\bm{i}}\|_{L^\infty([-M,M]^c)}$ is continuous.
		
		We only prove the continuity of $t\mapsto\|\mfrak{f}^{(j),t}_{l,\bm{i}}\|_{L^\infty([-M,M]^c)}$ in the case $j=5$, since the exact same arguments also prove the cases $j\in\{2,3,6\}$. Since $|\rho_{t}(x)-\rho_{t_0}(x)|\leq|\delta t |\max_{s\in[0,1]}\|u_s\|_{\infty}|\rho_{t_0}(x)$ by \eqref{def:u_t} and Corollary \ref{cor:contfixe}, and that the following map is uniformly bounded in $t\in[0,1]$ and $x>M$
		$$g_t:x\mapsto \dfrac{Q_{i_1}^{l}\left(\theta_t,\dots,\theta_t^{(i_1)}\right)(x)}{\rho_{t}(x)}\int\limits_{x}^{+\infty}\Big|P_{i_2}^{l-i_1}(\alpha_t,\dots,\alpha_t^{(i_2)})(y)\Big|\rho_{t_0}(y)\diff y,$$
		we can just show that $\|g_t-g_{t_0}\|_{\infty}\tend{t\rightarrow t_0}0$. Moreover by \eqref{eq:lemb3},
		it is enough to show the uniform convergence for $\dfrac{\rho_{t}}{\rho_{t_0}}g_t$. One can also notice that, for constants $C_{l,i_1}$ and $C_{l,\bm{i}}$ independents of $t$, $$|Q_{i_1}^l\left(\theta_t,\dots,\theta_t^{(i_1)}\right)(x)|\leq C_{l,i_1}|x|^{l-i_1}\hspace{0,5cm}\text{ and }\hspace{0,5cm}\Big|P_{i_2}^{l-i_1}(\alpha_t,\dots,\alpha_t^{(i_2)})(y)\Big|\leq C_{l,\bm{i}}|x|^{-(l-i_1)}.$$Finally, by writing $|x^{^{-(l-i_1)}}Q_{i_1}^l\left(\theta_t,\dots,\theta_t^{(i_1)}\right)(x)|$ as polynomial in $x^{-1}$, $t\phi^{(i+1)}(x)$ and $\mcal{H}[\rho_{t}^{(i)}]$ for $i\geq 0$ and $x^{l-i_1}P_{i_2}^{l-i_1}(\alpha_t,\dots,\alpha_t^{(i_2)})(x)$ as a rational function in those same variables, we conclude that these functions converge uniformly towards the same functions at $t_0$. Therefore
		$$x\mapsto x^{-(l-i_1)}Q_{i_1}^{l}\left(\theta_t,\dots,\theta_t^{(i_1)}\right)(x)\dfrac{x^{l-i_1}}{\rho_{t_0}(x)}\int\limits_{x}^{+\infty}\dfrac{\rho_{t_0}(y)\diff y}{y^{l-i_1}}\Big|y^{l-i_1}P_{i_2}^{l-i_1}(\alpha_t,\dots,\alpha_t^{(i_2)})(y)\Big|$$
		converges uniformly to the same functions at $t_0$. This establishes the proposition.
	\end{proof}	
	
	The final ingredient for showing the continuity with respect to $t$ of the constants in Theorems \ref{thm: cont inverse master}, \ref{thm: cont inverse master linfini} and \ref{thm:continversemaster thetaxi L2} is the continuity of these $L^2$-norms.
	
	\begin{lemma}\label{lem:appBl2}[Continuity $L^2$-norms]
		For all $j\in\llbracket1,4\rrbracket$, the maps $t\mapsto\|\mfrak{f}^{(j),t}_{l,\bm{i}}\|_{L^2([-M,M]^c)}$ are continuous where 
		\begin{equation}
			\mfrak{f}^{(4),t}_{l,\bm{i}}:x\mapsto\dfrac{Q_{i_1}^{l}\left(\theta_t,\dots,\theta_t^{(i_1)}\right)(x) }{\rho_{t}(x)}\cdot\Bigg|\int\limits_{x}^{\mrm{sgn}(x)\infty}\left[\alpha_t P_{i_2}^{l-i_1}\left (\alpha_t,\dots,\alpha_t^{(i_2)}\right)  \right] '(y)\rho_{t}(y)^2\diff y\Bigg|^{1/2}
		\end{equation}
		and where $\mfrak{f}^{(j),t}_{l,\bm{i}}$ for $j=1,2,3$ are defined in Lemma \ref{lem:b6}.
	\end{lemma}
	
	\begin{proof}
		For the continuity of $t\mapsto\|\mfrak{f}^{(1),t}_{l,\bm{i}}\|_{L^2([-M,M]^c)}$, we use dominated convergence theorem. Since, we showed uniform convergence and that $\|\mfrak{f}^{(1),t}_{l,\bm{i}}\|_{L^\infty([-M,M]^c)}<+\infty$, we conclude that for all $x\in[-M,M]^c$,
		$$\mfrak{f}^{(1),t}_{l,\bm{i}}(x)^2\tend{t\rightarrow t_0}\mfrak{f}^{(1),t_0}_{l,\bm{i}}(x)^2.$$
		The domination follows from the fact $\mfrak{f}^{(1),t}_{l,\bm{i}}(x)=\underset{|x|\rightarrow\infty}{O}(x^{-2})$ by Lemma \ref{lem:gammaetalpha}, and that all the dependence in $t$ is bounded, hence there exists a constant $C_{l,\bm{i}}>0$ independent of $t$ such that, for all $x\in[-M,M]^c$ and all $t\in[0,1]$,
		$$|\mfrak{f}^{(1),t}_{l,\bm{i}}(x)|\leq\dfrac{C}{x^2}.$$
		This establishes that $\|\mfrak{f}^{(1),t}_{l,\bm{i}}\|_{L^2([-M,M]^c)}\tend{t\rightarrow t_0}\|\mfrak{f}^{(1),t_0}_{l,\bm{i}}\|_{L^2([-M,M]^c)}$.
		
		We now establish the continuity of $t\mapsto\|\mfrak{f}^{(3),t}_{l,\bm{i}}\|_{L^2([-M,M]^c)}$, the case $j\in\{2,4\}$ can be obtained with the exact same arguments. We want to use dominated convergence theorem, for $\left({\mfrak{f}^{(3),t}_{l,\bm{i}}}\right) ^2$. The latter, when $t\rightarrow t_0\in[0,1]$, also converges uniformly since it is uniformly bounded and that we proved that $\mfrak{f}^{(3),t}_{l,\bm{i}}$ converges uniformly. It just remains to verify the domination hypothesis. By Lemma \ref{lem:gammaetalpha}, we know that  $\alpha_t(y)^2P_{i_2}^{l-i_1}\left (\alpha_t,\dots,\alpha_t^{(i_2)}\right )(y)^2=\underset{|y|\rightarrow\infty}{O}(y^{-2(l-i_1+1)}).$ We conclude by Lemma \ref{lem:M} that there exists a constant $C_{l,\bm{i}}>0$ such that for all $y>M$,
		$$\alpha_t(y)^2P_{i_2}^{l-i_1}\left (\alpha_t,\dots,\alpha_t^{(i_2)}\right )(y)^2\leq \dfrac{C_{l,\bm{i}}}{y^{2(l-i_1+1)}}.$$
		Similarly, $\Big|Q_{i_1}^{l}\left(\theta_t,\dots,\theta_t^{(i_1)}\right)(x)^2\Big|\leq C_{l,i_1}|x|^{2(l-i_1)}$ for all $x>M$ and for $C_{l,i_1}>0$ a constant independent of $t$ and $x$. Finally, we get the following domination for an arbitrary $t_0$ and all $x>M$:
		\begin{multline*}
			\Bigg|	\dfrac{Q_{i_1}^{l}\left(\theta_t,\dots,\theta_t^{(i_1)}\right)(x)^2}{\rho_{t}(x)^2}\int\limits_{x}^{+\infty}\alpha_t(y)^2P_{i_2}^{l-i_1}\left (\alpha_t,\dots,\alpha_t^{(i_2)}\right )(y)^2\rho_{t}(y)^2\diff y\Bigg|
			\\\leq\dfrac{C_{l,i_1}|x|^{2(l-i_1)}}{\rho_{t_0}(x)^{2}(1-\max_{s\in[0,1]}\|u_s\|_\infty)^{2}}\int_x^{+\infty}\dfrac{C_{l,\bm{i}}(1+\max_{s\in[0,1]}\|u_s\|_\infty)^{2}\rho_{t_0}(y)^{2}\diff y}{y^{2(l-i_1+1)}}.
		\end{multline*}
		The RHS is in $L^1([M,+\infty[)$ by integration by parts as it was done in the proof of Lemma \ref{lem:I_a}. We conclude by doing the same on $]-\infty,-M]$.
	\end{proof}
	
	\begin{proposition}\label{prop:constantescontinues}
		With the choice of potential $V_{\phi,t}$, the following map is continuous
		\begin{equation*}
			t\in[0,1]\mapsto\Big(C_{H^n}(\widetilde{\Xi_1^{-1}}),C_{W_n^{\infty}}(\widetilde{\Xi_1^{-1}}),C_{H^n}(\Theta^{(2)}\circ\widetilde{\Xi_1^{-1}})\Big)
		\end{equation*}
	\end{proposition}
	\begin{proof}
		By recalling the expression of those constants in \eqref{def:constanteXiHn}, \eqref{eq:C_9(V,l)}, \eqref{eq:C12}, it is clear that they are continuous expressions of the $t$-continuous quantities considered in Lemmas \ref{lemappB:invrho}, \ref{lem:I_a}, \ref{lemappB:norminfPoly}, \ref{lem:b6} and \ref{lem:appBl2}. The conclusion follows.
	\end{proof}
	\subsection{Parameter-continuity of $C_{\mcal{L}_t}$ and $K_{V_{G,\phi,t}}$}
	In \cite[App. A]{DwoMemin}, the authors showed that for a general potential $V$, the operator $\mcal{A}$ considered as an unbounded operator on $\msf{H}$ has the same spectrum as the Schrödinger operator $\mcal{S}:\mcal{D}(\mcal{S})\rightarrow L^2(\R)$, defined by $\mcal{S}\defi -\Delta+w_V$ with:
	\begin{equation*}
		\mcal{D}(\mcal{S})\defi \left\{u\in H^1(\R),uV'\in L^2(\R),-u''+w_V u\in L^2(\R), \int_\R u(x)\diff x=0\right\} 
	\end{equation*} and
	\begin{equation*}
		\label{eq:w_p}
		w_{V}\defi \frac{1}{2}\left( \frac{1}{2}V'^2-V'' + 2PV'\mathcal{H}[\rho_V]-2P\mathcal{H}[\rho_V']+2P^2\mathcal{H}[\rho_V]^2 \right)=\dfrac{1}{2}\Big[(\log\rho_V)''+\dfrac{1}{2}(\log\rho_V)'^2\Big]\,.
	\end{equation*}
	Since $C_{\mcal{L}}=\lambda_1(\mcal{A})=\lambda_1(\mcal{S})$ by Theorem \ref{thm inverse master}, we just have to show that when choosing the potential $V=V_{G,\phi,t}$, the $t$-dependent Schrödinger operator $\mcal{S}_t$ with potential $w_t\defi w_{V_{G,\phi,t}}$ produces a continuous smallest eigenvalue $\lambda_1(\mcal{S}_t)$. 
	\begin{proposition}\label{lem:contlambda(S)}
		The map $t\in[0,1]\mapsto C_{\mcal{L}_t}= \lambda_1(\mcal{S}_t)$ is continuous.
	\end{proposition}
	
	\begin{proof} First for all $t\in[0,1]$, $\lambda_1(\mcal{S}_t)>0$.
		Secondly, we have the following equalities:
		$$\lambda_1(\mcal{S}_t)=\min_{\substack{u\in\mcal{D}(\mcal{S}_t)
				\\\|u\|_{L^{2}(\R)}=1}}\Braket{u,\mcal{S}_t [u]}_{L^2(\R)}=\inf_{\substack{u\in\mcal{C}^{\infty}_c(\R)
				\\ \|u\|_{L^{2}(\R)}=1}}\int_{\R}(u(x)')^2\diff x+\int_\R u(x)^2w_{t}(x)\diff x.$$
		From the definition of $w_{t}$ and Lemma \ref{lem:leshilbertssontbornéesent}, it is easy to deduce that $$\|w_{t}-w_{t'}\|_{L^\infty(\R)}\tend{t\rightarrow t'}0.$$ Hence for all $t,t'\in[0,1]$, $u\in\mcal{C}^{\infty}_c(\R)$ with $\|u\|_{L^{2}(\R)}=1$, we have
		$$\left|\Braket{u,\mcal{S}_t [u]}_{L^2(\R)}-\Braket{u,\mcal{S}_{t'} [u]}_{L^2(\R)}\right|\leq\|w_{t}-w_{t'}\|_{L^\infty(\R)}  $$
		hence $\sup_{\substack{u\in\mcal{C}^{\infty}_c(\R)
				\\ \|u\|_2=1}}\left|\Braket{u,\mcal{S}_t [u]}_{L^2(\R)}-\Braket{u,\mcal{S}_{t'} [u]}_{L^2(\R)}\right|\tend{t\rightarrow t'}0.$ Since uniform convergence is enough to ensure convergence of infinimums, we get the result.
	\end{proof}
	
	We know prove the continuity of the constant $K_{V_{G,\phi,t}}$ introduced in Theorem \ref{thm:bounddensity}.
	\begin{lemma}\label{lem:contk_V}
		The following map is continuous $$t\mapsto K_{V_{G,\phi,t}}=2P\|\mcal{H}[\rho_{t}]\|_{\infty}+C+P\displaystyle\Big|\iint_{\R^2}\log|x-y|\rho_t(x)\rho_t(y)\diff x\diff y\Big|$$
		for $C$ some fixed constant (independent of $t$).
	\end{lemma}
	\begin{proof}
		We already proved the continuity of the map $t\mapsto\|\mcal{H}[\rho_{t}]\|_{L^\infty(\R)}$ in Lemma \ref{lem:leshilbertssontbornéesent} so it just remains to show that the double integral is continuous with respect to $t$. We prove this by dominated convergence theorem. The function $(x,y)\mapsto\log|x-y|\rho_{t}(x)\rho_{t}(y)$ converges almost everywhere to $$(x,y)\mapsto\log|x-y|\rho_{t_0}(x)\rho_{t_0}(y)$$ as $t$ goes to $t_0$. Furthermore we have the following domination $(x,y)$-almost everywhere $$\Big|\log|x-y|\Big|\rho_{t}(x)\rho_{t}(y)\leq\Big|\log|x-y|\Big|(1+\max_{s\in[0,1]}\|u_s\|_{\infty})^2\rho_{t_0}(x)\rho_{t_0}(y).$$
		This allows us to conclude on the continuity of $t\mapsto K_{V_{G,\phi,t}}$.
	\end{proof}
	
	\begin{proposition}\label{prop:integrabilityremainder}
		The map $t\mapsto C_{\mrm{rem}}(V_{G,\phi,t})$ where $C_{\mrm{rem}}(V)$ is defined in \eqref{eq:controleraminedercorrelator} is integrable on $[0,1]$.
	\end{proposition}
	\begin{proof}
		By the bounds on  $C_{H^n}(\widetilde{\Xi_1^{-1}})$, $C_{W_n^{\infty}}(\widetilde{\Xi_1^{-1}}),C_{H^n}(\Theta^{(2)}\circ\widetilde{\Xi_1^{-1}})$
		in Theorems \ref{thm: cont inverse master}, \ref{thm: cont inverse master linfini} and \ref{thm: cont theta inverse master linfini} in addition to the continuity results of Lemmas \ref{lem:contlambda(S)} and \ref{lem:contk_V}, we conclude that $C_{\mrm{rem}}(V_{G,\phi,t})\in L^{1}([0,1])$.
\end{proof}}

\bibliographystyle{alpha}
\bibliography{biblio}
		
\end{document}